# How to estimate the total number of citations of a researcher using his $h$-index and his $h-$core?

Romeo Meštrović[1*] and Branislav Dragović[1]
[1]Maritime Faculty Kotor, University of Montenegro,
Dobrota, 85330 Kotor, Montenegro; Tel./Fax: +382 303 184
E-mails: romeo@ucg.ac.me; branod1809@gmail.com

**Keywords**: Hirsch $h$-index, citations, $h$-core, partition of a positive integer, Durfee square, asymptotic normality of the $h$-index, $h$-defect of asymptotic normality of citations, confidence interval.

*Mathematics Subject Classification* **2020** (**MSC 2020**). 05A17, 11P81, 68-11, 62D05.

## Abstract

So far, many researchers have investigated the following question: Given total number of citations, what is the estimated range of the $h$-index? Here we consider the converse question. Namely, the aim of this paper is to estimate the total number of citations of a researcher using only his $h$-index, his $h$-core and perhaps a relatively small number of his citations from the tail. For these purposes, we use the asymptotic formula for the mode size of the Durfee square when $n \to \infty$, which was proved by Canfield, Corteel and Savage (1998), seven years before Hirsch (2005) defined the $h$-index. This formula confirms the asymptotic normality of the Hirsch citation $h$-index. Using this asymptotic formula, in Section 4 we propose five? estimates of a total number of citations of a researcher using his $h$-index and his $h-$core. These estimates are refined mainly using small additional citations from the $h-$ tail of a researcher. Related numerous computational results are given in Section 5. **Notice that the relative errors $|\delta_B|$ of the estimate $B$ of a total number of citations of a researcher are surprisingly close to zero for E. Garfield, H.D. White (Table 2), G. Andrews, L. Leydesdorf and C.D. Savage (Table 5).**

## 1. Introduction

In the past two decades, numerous Scientometrics and Bibliometrics indicators were proposed to evaluate the scientific impact of individuals, institutions, colleges, universities and research teams. To appraise the scientific impact of scholars, institutions and research areas among others, several publication-based indicators are used, such as the size-dependent indicators (total number of citations and number of highly cited papers) and size-independent indicators (average number of citations per paper and proportion of highly cited publications) (Waltman, 2016), as well as citation frequency (life cycle) of papers etc. Based on the limitations of these indicators, Hirsch (2005) proposed a new indicator called *h*-index, whilst Egghe (2006a, b, c) defined and studied an improvement of the $h$-index, called the $g$-index.

The $h$-index is defined as the highest number of publications of a scientist that received $h$ or more citations (Hirsch, 2005). Hirsch (2007) reported results of an empirical study of the predictive power of the *h*-index compared with other bibliometric indicators. Hirsch indicated that the $h$-index is better than other indicators considered (total citation count, citations per paper and the total paper count) in predicting future scientific achievement.

The $g$-index is defined as the highest rank such that the cumulative sum of the number of citations received is larger than or equal to the square of this rank. As noticed by Bihari et al. (2021), the $h$-index and the $g$-index got a lot of attention due to its simplicity, and several other indicators were proposed to extend the properties of the $h$-index and to overcome its shortcomings. Bornmann, Mutz, Hug and Daniel (2011) noticed that many studies have investigated many variants of the $h$-index, but these variants do not offer clear advantages over the original $h$-index (also see Anderson et al., 2008).

The $h$-index and its variants are extensively studied in Informetrics (both experimentally and theoretically). Numerous generalizations, extensions and variations of the $h$-index and the $g$-index were defined, studied and compared in the last two decades. The comprehensive recent review on $h$-index and its alternative indices was presented by Bihari et al. (2021) (see also Alonso et al., 2009).

So far, many researchers have investigated the following question: Given total number of citations, what is the estimated range of the $h$-index? Here we consider the converse question. Namely, the aim of this paper is to estimate the total number of citations of a researcher using only his $h$-index, his $h$-core and perhaps a small number of his citations from the tail. For these purposes, we use the asymptotic formula for the mode size of the Durfee square when $n \to \infty$, which was proved by Canfield, Corteel and Savage (1998), seven years before Hirsch (2005) defined the $h$-index.

This article is organized as follows. Section 2 contains the definitions of bibliometric indicators and new indices which are used in our computational results. Preliminary definitions on a partition of a positive integer $n$ are presented in Section 3. We present the asymptotic formula for the mode size of the Durfee square when $n \to \infty$, which was proved by Canfield, Corteel and Savage (1998). As it was observed by Ekhad and Zeilberger (2014), this formula confirms the asymptotic normality of the Hirsch citation index (alias size of Durfee square; see Subsection 3.1). Using this asymptotic formula, in Section 4 we propose five confidence intervals to estimate a total number of citations of a researcher using his $h$-index and his $h$-core. Related numerous computational results concerning 41 research profiles at Scopus and Google Scholar databases (which include 38 different researchers) are presented in Section 5. These computational results suggest the fact that generally, the means of these intervals present very good approximations of a total number of citations of a researcher. Concluding remarks are given in Section 6.

## 2. The $h$-, $g$-, $A$-, $R$-, $e$, $H$- and $D$-indices

Suppose that a scientist $S$ has published $p$ publications whose number of citations $cit_1, cit_2, cit_3, ..., cit_p$ in considered database are usually ranked in decreased order, i.e.,

$$cit_1 \geq cit_2 \geq \cdots \geq cit_p \geq 0. \qquad (1.1)$$

In the whole paper, we will use the following notations and related notions.

$p$ - the total number of publications of a scientist in considered database;

$N_p^+$ - the total number of publications having at least one citation;

$N_{cit} = \sum_{j=1}^{p} cit_j$ - the total number of citations of a scientist in considered database;

$N_{cit}(s) = \sum_{j=1}^{s} cit_j$ ( $s = 1,2,...,n$ ) - the total number of citations of a scientist in considered database up to the rank $s \in \{1,2,...,p\}$;

$N_{cit}(h) = \sum_{j=1}^{h} cit_j$ - the total number of citations of a scientist in considered database into $h$-core.

The Hirsch or $h$-*index* (Hirsch, 2005) is defined as the highest rank $m$ such that the first $m$ publications of a scientist received $m$ or greater than $m$ citations in considered database, or equivalently,

$$h = \max\{m : cit_m \geq m, 1 \leq m \leq p\}. \tag{1.2}$$

The Egghe's $g$-*index* (Egghe, 2006a, b, c) is defined as the largest positive integer $k$ ( $k = 1,2,...,n$ ) such that

$$\frac{1}{k}\sum_{j=1}^{k} cit_j \geq k. \tag{1.3}$$

Obviously, we have $h \leq g$.

Jin's $A$-*index* (the name was suggested by Rousseau, 2006), introduced by Jin (2006) is defined as the average number of citations received by the publications into $h$-core, i.e.,

$$A = \frac{1}{h}\sum_{j=1}^{h} cit_j, \tag{1.4}$$

which is under our notation equal to $N_{cit}(h)/h$.

Observe that the $A$-index achieves the same goal as the $g$-index, i.e. it takes into account the total number of citations included in the $h$-core. The $A$-index could potentially disadvantage prolific researchers, so the $A$-index was defined to remove this problem (Jin et al., 2007)

Another attempt to improve the insensitivity of the $h$-index to the number of citations to highly cited papers is the $R$-index, introduced in Jin et al. (2007. p. 857). The $R$-*index* as defined as

$$R = \sqrt{\sum_{j=1}^{h} cit_j}, \tag{1.5}$$

So as the $h$-index, also this measure takes into account the actual $cit_j$ - values in the $h$-core. Note that

$$R = \sqrt{hA} = \sqrt{N_{cit}(h)}. \tag{1.6}$$

It is easy to show that (see Proposition 1 and Corollary in Jin et al, 2007)

$$h \leq g \leq A.$$

Moreover, from $R = \sqrt{hA}$ and $A \geq h$ by Corollary in Jin et. al. (2007, p. 657) it follows that

$$h \leq R.$$

Aimed at the same goal as the $g$-index, the $e$-*index* was proposed in 2009 by (Zhang, 2009; also see Zhang, C.-T., 2013) with the aim of considering the contribution of excess citations, which are mostly from highly cited papers. It is defined as

$$e = \sqrt{\sum_{j=1}^{h} cit_j - h^2}. \tag{1.7}$$

Clearly, one has $e = 0$ if and only if $cit_1 = cit_2 = \cdots = cit_h = h$. Since by (1,5), $N_{cit}(h) = \sum_{j=1}^{h} cit_j = R^2 = hA$, we find that

$$e = \sqrt{hA - h^2} = \sqrt{R^2 - h^2} = \sqrt{N_{cit}(h) - h^2} \tag{1.8}$$

and hence,

$$\frac{e}{R} = \sqrt{1 - \left(\frac{h}{R}\right)^2} = \sqrt{1 - \frac{h^2}{N_{cit}(h)}}. \tag{1.9}$$

As already highlighted in some papers (see e.g., Jin et al., 2007 and Todeschini, 2011), the indices $g$-, $A$-, $R$- and $e$- are very sensitive to few highly-cited papers. In particular (see Todeschini, 2011), $g$-, $A$-, $R$- and $e$-indices appear to measure different characteristics of the author production because the information related to the number of citations prevails over that related to the number of papers in the $h$-core (e.g., the $g$-index). Burrell (2007) noticed that numerical investigations suggest that the $A$-index is a linear function of time and of the $h$-index, while the size of the Hirsch $h$-core has an approximate square-law relationship with time, and hence also with the $A$-index and the $h$-index. Notice that the $e$- index takes into account the contribution of excess citations.

Recently, Meštrović and Dragović (2023) defined the $H$- index as

$$H = h\sqrt{r}, \tag{1.10}$$

where (Corollary 2.15 in Meštrović and Dragović, 2023)

$$r = \left\lfloor \frac{2N_{cit}(h)}{h(h+1)} \right\rfloor - 1 = \left\lfloor \frac{2R^2}{h(h+1)} \right\rfloor - 1 = \left\lfloor \frac{2A}{h+1} \right\rfloor - 1, \tag{1.11}$$

where $\lfloor x \rfloor$ denotes the greatest integer which does not exceed a real number $x$.

The continuous analague of the value $r$ is the value $q$ defined as

$$q = \frac{2N_{cit}(h)}{h^2} - 1 = \frac{2R^2}{h^2} - 1. \tag{1.12}$$

Then the index associated to the above value $q$ is the $D$-index defined as (Meštrović and Dragović, 2023)

$$D = h\sqrt{q} = h\sqrt{\frac{2N_{cit}(h)}{h^2} - 1} = \sqrt{2N_{cit}(h) - h^2} = \sqrt{2R^2 - h^2}. \tag{1.13}$$

Since $r \leq q$, it follows that $H \leq D$. Moreover, it is easy to prove the following inequalities (Meštrović and Dragović, 2023):

$$h \leq R \leq D \leq A.$$

**Remark 2.1.** Notice that the ratio $q/e = (2N_{cit}(h)/h^2 - 1)/\sqrt{N_{cit}(h) - h^2}$ plays a crucial role in this paper (see Sections 4 and 5) in order to estimate the total number of citations of researchers using their $h$-core and mainly a relatively small number of citations from the tail.

**Remark 2.2.** If we consider the $h$-index as a Durfee square of length $h$ of a random partition of integer $N_{cit}$ which is defined in Section 3, then $N_{cit}(h)$, the indices $A$, $R$ and $e$ are compound like discrete random variables because the $h$-index is the integer random variable (see Feller, 1968, Chapter XII and Meštrović, 2015).

## 3. The asymptotic normality of the Hirsch citation index

Brito and Rodriguez Navarro (2021) noticed that many studies, including the original paper proposing the $h$-index, have claimed that the $h$-index can be calculated from the total number of citations (see Malesios, 2015). Specifically, van Raan (2006) found that a power law links the $h$-index with the sum of the number of citations ($N_{cit}$). Consistent with this observation, assumed that citation distributions are lognormal, Brito and Rodriguez (2021) observed good fits of their data for $N_{cit}$ and the $h$-index to a power law (not shown); the exponent 0.42 is not very different from the exponent given by van Raan (2006), which is 0.45. The dependence of $h$-index and $N_{cit}$ has been extensively studied by Egghe and Rousseau (2006), Glänzel (2006), Malesios (2015), Molinari and Molinari (2008) and Montazerian and al. (2019). Brito and Rodriguez Navarro (2021) observed that the size-independent indicator $h/p$ ($p$ is a number of publications) shows no correlation with the probability of publishing a paper exceeding any of the citation thresholds.

Yong (2014) considered the following question: Given $N_{cit}$, what is the estimated range of $h$? Yong (2014) noticed that taking only $N_{cit}$ as input hardly seems like sufficient information to obtain a meaningful answer. For computational purposes, Yong used the expression (3.3) established by Canfield, Corteel and Savage (1998) which is the asymptotic formula for the most likely Durfee square size for a partition in the set of all partition of a positive integer $n$. Yong (2014) reinterpret this result as the "rule of thumb" for $h$-index. Theoretically, Yong (2014) noticed that one expects the rule of thumb to be nearly correct for large $N_{cit}$. He noticed that this is empirically not true, even for pure mathematicians. However, Yong (2014) observed that something related to $0.54044\sqrt{N_{cit}}$ is higher than the actual $h$-index for almost every very highly cited ($N_{cit} > 10000$) scholar in Mathematics, Physics, Computer Science and Statistics.

Here we consider the converse question. Namely, as it is noticed in Abstract and Introduction, the aim of this paper is to estimate the total number of citations of a researcher using only his $h$-index, the $h$-core and a relatively small number of citations from the tail. For these purposes, here we use the asymptotic formula for the mode size of the Durfee square when $n \to \infty$, which was proved by Canfield, Corteel and Savage (1998) seven years before Hirsch (2005) defined the $h$-index (see the asymptotic formula (3.3) in Subsection 3.2).

### 3.1. Preliminary definitions on a partition of a positive integer $n$

A *partition* $\lambda$ of a positive integer $n$ is a non-increasing sequence of positive integers $\lambda_1, \lambda_2, ..., \lambda_k$ (i.e., $\lambda_1 \geq \lambda_2 \geq \cdots \geq \lambda_k$), each called a *part of the parttion*, such that

$$n = \lambda_1 + \lambda_2 + \cdots + \lambda_k. \tag{3.1}$$

We will also write $|\lambda| = n$. The terminology is extended to allow us to say that $n = \lambda_1 + \lambda_2 + \cdots + \lambda_k$ is a partition of $n$. So, for example, the all partitions of 5 are

$$5, \ 4+1, \ 3+2, \ 3+1+1, \ 2+2+1, \ 2+1+1+1 \text{ and } 1+1+1+1+1.$$

So, there are seven partitions of 5.

The *partition function* $p(n)$ equals the number of all partition of a positive integer $n$. For example, $p(5) = 7$. Notice that the sequence $\{p(n)\}$ (with the initial term $p(0) = 1$) is the OEIS A070177 sequence in Sloane, N.J.A. For example, $p(10) = 42$, $p(20) = 627$, $p(30) = 5604$,

$p(30) = 5604$, $p(100) = 190569292$, $p(1000) \approx 2.406147 \cdot 10^{30}$ and
$p(10000) \approx 3.616725 \cdot 10^{105}$.

No closed form expression for the partition function $p(n)$ is known, but it has both asymptotic expansion that accurately approximate it and recurrence relations by which it can be calculated exactly. It grows as an exponential function of the exponential function of the square root of its argument (Andrews, 1976, p. 69), as follows:

$$p(n) \approx \frac{1}{4n\sqrt{3}} \exp\left(\pi\sqrt{\frac{2n}{3}}\right) \text{ as } n \to \infty.$$

The above asymptotic formula was first obtained by Hardy and Ramanujan (1918) and independently by J.V. Uspensky in 1920 (https://en.wikipedia.org/wiki/Partition_function_(number_theory)).

To a partition $(\lambda_1, \lambda_2, ..., \lambda_k)$ with $\lambda_1 \geq \lambda_2 \geq \cdots \geq \lambda_k$, let us associate the so-called *Ferrers diagram*, which is a left-justified array of dots, with $\lambda_i$ dots in the $i$ th row $(i = 1, 2, ..., k)$. If the dots are replaced by squares, we obtain the *Young diagram* representation of a partition $\lambda$. This graphical representation contains a square in the upper left corner, and the largest such square is called *Durfee square* (see e.g., Andrews and Eriksson, 2004, p. 76 or Hardy and Wright, 1960, p. 281). Durfee squares are named after American mathematician William Pitt Durfee, a student of English mathematician J.J. Sylvester. This is attributed in Sylvesters' letter to Arthur Cayley in 1883 (Hunger, 1998, p. 224; https://en.wikipedia.org/wiki/Durfee_square#cite_note-4) in which he wrote: "*Durfee's square is a great invention of the importance of which its author has no conception.*".

It is known (see e.g., Andrews, 1976) that the generating function for partitions with Durfee square of length $d$ is given by

$$\frac{q^{d^2}}{\prod_{j=1}^{d}(1-q^j)^2}. \tag{3.2}$$

The above expression is closely related to the following famous Euler-Gauss identity (Andrews, 1976, pp. 27-28, Yong, 2014):

$$\prod_{i=1}^{\infty}\frac{1}{1-q^i} = \sum_{k=0}^{\infty}\frac{q^{d^2}}{\prod_{j=1}^{d}(1-q^j)^2}.$$

For example, for $d = 3$, the Maclaurin series expansion of the expression (3.2) is

$$\frac{q^9}{\prod_{j=1}^{3}(1-q^j)^2} = q^9 + 2q^{10} + 5q^{11} + 10q^{12} + 18q^{13} + 30q^{14} + 49q^{15} + 74q^{16} + 110q^{17} + 158q^{19} + \cdots.$$

Hence, as the coefficient of $x^{11}$ in the above expansion is 5, there are 5 unrestricted partitions of $n = 11$ whose Durfee squares is 3. These partitions are $5+3+3$, $4+4+3$, $4+3+3+1$, $3+3+3+2$ and $3+3+3+1+1$. Observe that every partition of 11 (whose total number is $p(11) = 56$) has Durfee square that is les than equal of 3.

For example, the Young diagram of the partition (5,3,2,1,1) of 12 is

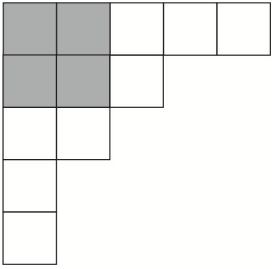

In terms of Ferrers diagram, the Durfee square can be defined as follows. The largest square in the Ferrers diagram of a partition is called *Durfee square*. If there are $\lambda_i$ dots in the $i$ th row ($i = 1, 2, ..., k$), then its size $d$ can be described as the number of parts $\lambda_i$ satisfying the condition $\lambda_i \geq i$. In the above example, the Durfee square has side of length 2 and area 4, and so $d = 2$.

The above standard definitions, as well as the one for a Durfee square, can be found for example, in Andrews and Eriksson (2004). Durfee squares in partitions have been previously studied with respect to the length of the side by Canfield, Corteel and Savage (1998), and independently by Mutafchiev (2002), where it was proved a local theorem for the length of the side of the Durfee square in a random partition of a positive integer $n$ as $n \to \infty$.

In Canfield, Corteel and Savage (1998) the emphasis was an asymptotic estimates for the length (see Subsection 3.2), while in Blecher (2012) some related interesting partition identities are derived.

### 3.2. The asymptotic normality of the $h$ - index

Ekhad and Zeilberger (2014) wrote: "When Rodney Canfield, Sylvie Corteel, and Carla Savage wrote their beautiful article (Canfield, Corteel and Savage, 1998), proving, rigorously, by a very deep and intricate analysis, the *asymptotic normality* of the random variable "size of Durfee square" defined on integer-partitions of $n$ (as $n \to \infty$), with precise asymptotics for the mean and variance, **they did not dream that one day their result should be of interest to everyone who has ever published a paper**."

Accordingly, the asymptotic normality of Hirsch $h$ - index was provided by Canfield, Corteel and Savage (1998) seven years before Hirsch (2005) defined the $h$-index.

Let $P(n)$ be the set of all partitions of a positive integer $n$, and let $P(n,d)$ be the set of all partitions of $n$ with Durfee square of size $d$. Let $F(n)$ be a family of partitions of a positive integer $n$ and let $F(n,d)$ denote the set of partitions in $F(n)$ with Durfee square of size $d$. For a fixed $n$, Canfield, Corteel and Savage (1998) investigated the sequences $\{F(n,d)\}$ ($0 \leq d \leq \sqrt{n}$). Denote by $m_F(n)$ and $a_F(n)$ the most likely size and the average size of the Durfee square of a partition in $F(n)$, respectively.

It was proved in Canfield, Corteel and Savage (1998, Section 5) that if $\{F(n,d)\}$ satisfies a recurrence of the type in Section 2 of the same paper and has a particular asymptotic form, then for fixed $n$, the most likely value of the Durfee square of a partition in $F(n)$ grows as $c_F \sqrt{n}$ for some constant $c_F$ depending only on the recurrence for $F(n,d)$. Theorem 1 in this paper gives the asymptotic formula for $P(n,d)$. Consequently, by Corollary 1, the numbers $P(n,d)$ are asymptotically normal as $n \to \infty$. It was observed in Canfield, Corteel and Savage (1998) that using the asymptotic formula for $P(n,d)$, it can be showed that the most likely Durfee square size for a partition in $P(n)$ (i.e., a mode of Durfee square size of $P(n)$) is (see the row (1) of Table 2 in Canfield, Corteel and Savage, 1998).

$$m_P(n) \approx \frac{\sqrt{6}\ln 2}{\pi}\sqrt{n} \approx 0.54044\sqrt{n}. \tag{3.3}$$

The above asymptotic formula presents the mode size of the Durfee square when $n \to \infty$. Moreover (Canfield, Corteel and Savage, 1998), it follows that for sufficiently large $n$,

$$|m_F(n) - a_F(n)| < 1/2 + o(1) \tag{3.4}$$

and that $\{P(n,d)\}$ ($\varepsilon n^{1/2} \leq d \leq (1-\varepsilon)n^{1/2}$) is log-concave.

**Remarks 3.1.** Inspired by Alexander Yong's recent critique (Yong, 2014) of the Hirsch citation index, Ekhad and Zeilberger (2014) gave an empirical (yet very convincing!) reproof of the asymptotic normality of the Hirsch citation index (alias size of Durfee square) with respect to the uniform distribution on the "sample space" of integer partitions of $n$ (see Subsection 3.1 in Ekhad and Zeilberger, 2014). It was noticed in the mentioned paper that the result of Canfield, Corteel and Savage (1998) was first proved rigorously (but with much greater effort!). Ekhad and Zeilberger (2014) confirmed the asymptotic formula (3.3) and estimated the variance numerically as $0.00811\sqrt{n} + O(1)$, and were obtained estimates extremely close to those of the standard Normal Distribution for the first 12 standardized moments. Namely, Ekhad and Zeilberger (2014) numerically estimated asymptotic expressions for the expectation $a_n$ and the variance $m_2(n)$ of Durfee square size for a partition in $P(n)$ as

$$a_n = 0.540446395\sqrt{n} + 0.085691 + 0.0374788\frac{1}{\sqrt{n}} + O\left(\frac{1}{n}\right), \tag{3.5}$$

and

$$m_2(n) = 0.081057\sqrt{n} + 0.018459 - 0.018015\frac{1}{\sqrt{n}} + O\left(\frac{1}{n}\right). \tag{3.6}$$

**Remarks 3.2.** Hirsch (2005) founded that the $h$-index is proportional to the square root of the total number of citations. Namely, Hirsch proposed the approximation $h \approx \sqrt{N_{cit}/a}$, where $a$ is a constant between 3 and 5. In Hirsch (2005) the constant 0.54044 corresponds to a value of 3.424 for the constant $a$. This asserts that the $h$-index is between $\sqrt{N_{cit}/5} \approx 0.45\sqrt{N_{cit}}$ and $\sqrt{N_{cit}/3} \approx 0.58\sqrt{N_{cit}}$ (cf. Ruch and Ball R., 2010).

If a researcher has $p$ publications and each of the $m$ most cited publications has at least $m$ citations, then the maximal value $m$ is the empirical $h$-index of this researcher, introduced by Hirsch (2005). In the papers Glänzel (2006) and Egghe and Rousseau (2006) some theoretical properties of the $h$-index were derived. Glänzel (2006) pointed out that the $h$-index is certainly an interesting indicator with interesting mathematical properties. The strength of this index lies in the potential application to the assessment of small paper sets were other, traditional bibliometric indicators often fail or at least were their application proved usually problematic.

Glänzel (2006) founded a strong relation between both number of papers published and citations received with $h$-index for practically all Paretian distributions.

More precisely, Glänzel (2006) gave a probabilistic approach to the definition of the theoretical $h$-index given as follows. Let $X$ be a discrete random variable which represents the citation rate of a paper. Let the probability distribution function of $X$ is given as $p_k = P(X = k)$ for all $k = 0,1,2,...$ and the cumulative distribution function defined as $F(k) = P(X < k)$. If we define $G_k = G(k) := 1 - F(k) = P(X \geq k)$, then Gumbel's $r$-th characteristic extreme value ($u_r$) is defined as $u_r = G^{-1}(r/n) = \max\{k : G(k) \geq r/n\}$, where $n$ is a given sample with distribution $F$ (i.e., the number of publications of researcher). The *theoretical $h$-index* ($H$) can consequently be defined as $H = \max\{r : u_r \geq r\} = \max\{r : \max\{k : G(k) \geq r/n\} \geq r\}$. If there exists such index $r$ so that $u_r = r$, then we have obviously $H := r$, and we can write $H := u_H$ i.e., the theoretical $h$-index $h = h_n$ of a researcher is defined (cf. Beirlant and Einmahl, 2007).

Glänzel (2006) considered two special cases, namely (1) discrete Paretian distributions and (2) the Price distribution with finite expectation. Glänzel (2006) pointed out that these two distributions cover most distributions used for modeling publication activity and citation processes. The distribution law of Paretian random variable satisfies the condition $p_k = P(X = k) \approx d(N + k)^{-(\alpha+1)}$ for sufficiently large $k$, where $\alpha > 1$, $N$ and $d$ are positive constants. Then the $h$-index is proportional to the $(\alpha+1)$ th root of the number of publications (Glänzel, Subsection 2.1).

In the case of Price distribution, introduced by Glänzel and Schubert (1985), Glänzel (2006, Subsection 2.2) proved that the $h$-index is proportional to the square root of the number of publications (cf. Price's square root law in Glänzel and Schubert, 1985).

Let $X_1,...,X_n$ be i.i.d. random variables with common distribution function $F$ as defined above. Assuming that the function $F$ is continuous and denoting (the right continuous) empirical distribution function of the random variables $X_1,...,X_n$, by $\hat{F}$, the empirical Hirsch index $\hat{H}$ is defined by Beirlant, and Einmahl (2007) as

$$1 - \hat{F}(\hat{H}) \leq \frac{\hat{H}}{n} \quad \text{and} \quad 1 - \lim_{x \uparrow \hat{H}} \hat{F}(\hat{H}) \geq \frac{\hat{H}}{n}.$$

Beirlant, J. and Einmahl, H.J. (2007) noticed that the above two inequalities have a unique solution. Furthermore, in case $\hat{F}$ puts only mass at non-negative integer values this definition coincides with the Hirsch definition of the $h$-index.

Notice that Beirlant and Einmahl (2007, Proposition 1 and Theorem 1) established the asymptotic normality of the empirical $h$-index. In case that the citations follow a Pareto-type or a Weibull-type distribution as defined in extreme value theory, a general result of this paper nicely specializes to results that are useful for constructing confidence intervals for the $h$-index on the basis of $\hat{H}$.

## 4. Auxiliary results

In this section we give useful definitions and propositions which will be necessary to obtain computational results presented in Section 5.

In accordance with Subsection 3.2 and the asymptotic formula (3.3) we give the following definition.

**Definition 4.1.** The *normal approximation* (*NA*) of the $h$-index, denoted as $h_{NA}$, is defined as

$$h_{NA} = \frac{\sqrt{6}\ln 2}{\pi}\sqrt{N_{cit}} \approx 0.54044463\,94667307\sqrt{N_{cit}}, \qquad (4.1)$$

where $N_{cit}$ is a total number of citations of a scientist in considered database.

*Examples.* Since $N_{cit} \geq N_{cit}(h) \geq h^2$, from (4.1) it follows that

$$\frac{h_{NA}}{h} \geq \frac{\sqrt{6}\ln 2}{\pi} \approx 0.54044463\,94667307,$$

where obviously, equality holds if and only if $cit_1 = cit_2 = \cdots = cit_h = h$ and $cit_i = 0$ for all $i = h+1,\ldots, p$.

For any positive integer $k$, put $cit_i = k+1-i$ for all $i = 1,\ldots, k$. Then clearly, $p = N_p^+ = k$, $N_{cit} = k(k+1)/2$ and $h = \lfloor (k+1)/2 \rfloor$. Therefore, we have

$$\frac{h_{NA}}{h} = \frac{\sqrt{6}\ln 2}{\pi} \cdot \frac{\sqrt{k(k+1)/2}}{\lfloor (k+1)/2 \rfloor} \to \frac{2\sqrt{3}\ln 2}{\pi} = 0.76430413\,88456883 \text{ as } k \to \infty.$$

More generally, for positive integers $k$ and $d$ and a nonnegative integer $a$, put $cit_i = a + (k-i)d$ for all $i = 1,\ldots, k$. Then clearly, $p = N_p^+ = k$, $N_{cit} = k(2a+(k-1)d)/2$ and $h = \lfloor (a+kd)/(d+1) \rfloor$. Therefore, we have

$$\frac{h_{NA}}{h} = \frac{\sqrt{6}\ln 2}{\pi} \cdot \frac{\sqrt{k(2a+(k-1)d)/2}}{\lfloor (a+kd)/(d+1) \rfloor}.$$

Letting $k \to \infty$ for fixed $d$ and $a$, from the above equality it follows that

$$\frac{h_{NA}}{h} \to \frac{\sqrt{6}\ln 2}{\pi} \cdot \frac{d+1}{\sqrt{2d}} := f(d) \text{ as } k \to \infty.$$

For example, $f(1) = 0.764$, $f(2) = 0.811$, $f(3) = 0.883$, $f(4) = 0.955$ $f(5) = 1.025$, $f(6) = 1.093$, $f(7) = 1.156$, $f(8) = 1.216$, $f(9) = 1.274$, $f(10) = 1.329$ and $f(26) = 2.024$. Hence, the arithmetic progressions with differences $d \in \{4,5,6\}$ give a good approximation $h_{NA}/h \approx 1$ for large integer values $k$.

**Definition 4.2.** Let a scientist $S$ has published $p$ publications, each of which received at least one citation. As usually, these citations $cit_1, cit_2, cit_3\ldots, cit_p$ are ranked in decreased order, i.e.,

$$cit_1 \geq cit_2 \geq \cdots \geq cit_p \geq 1. \qquad (4.2)$$

Then the *shifted h-index of order k* ($k = 0,1,...,p-1$), denoted as $h_k$, is defined as the $h$-index concerning the citations $cit_{k+1}, cit_{k+2},...,cit_p$. Observe that $h_0 = h$ and $h_{p-1} = 1$.

We also denote the total number of citations into $h_k$-core as $N_h^{(k)}$ ($k = 0,1,...,p-1$), i.e.,

$$N_h^{(k)} = \sum_{j=k+1}^{k+h_k} cit_j \text{ - } N_h^{(k)} = \sum_{j=k+1}^{k+h_k} cit_j . \tag{4.3}$$

Furthermore, for all $k = 0,1,...,p-1$ we also define $N_{cit}^{(k)}$, $e_k$ and $q_k$ as

$$N_{cit}^{(k)} = \sum_{j=k+1}^{p} cit_j \text{ - } N_{cit}^{(k)} = \sum_{j=k+1}^{p} cit_j , \tag{4.4}$$

$$e_k = \sqrt{N_h^{(k)} - (h_k)^2} = \sqrt{\sum_{j=k+1}^{k+h_k} cit_j - (h_k)^2} \tag{4.5}$$

and

$$q_k = \frac{2 N_h^{(k)}}{(h_k)^2} - 1. \tag{4.6}$$

Notice that $N_h^{(0)} = N_{cit}(h) = \sum_{j=1}^{h} cit_j$, $N_{cit}^{(0)} = N_{cit} = \sum_{j=1}^{np} cit_j$, $e_0 = e = \sqrt{N_{citt}(h) - h^2} = \sqrt{\sum_{j=1}^{h} cit_j - h^2}$ and $q_0 = q = 2 N_{cit}(h)/h^2 - 1$.

**Definition 4.3.** Under notations of Definition 4.2, for all $k = 0,1,...,p-1$ we define the $h_{NA}^{(k)}$-*index* as

$$h_{NA}^{(k)} = \frac{\sqrt{6} \ln 2}{\pi} \sqrt{N_{cit}^{(k)}} = 0.54044463\ 94667307 \sqrt{N_{cit}^{(k)}} . \tag{4.7}$$

Observe that $h_{NA}^{(0)} = h_{NA}$, where $h_{NA}$ is the normal approximation of the $h$-index defined by the expression (4.1) of Definition 4.1.

For all $k = 0,1,...,p-1$ we also define the intervals $I_k = (a_k, b_k)$ and $J_k = (c_k, f_k)$ as

$$I_k = \left(\left(\frac{\pi}{\sqrt{6}\ln 2}h_k\left(1-\frac{q_k}{e_k}\right)\right)^2, \left(\frac{\pi}{\sqrt{6}\ln 2}h_k\left(1+\frac{q_k}{e_k}\right)\right)^2\right) = \left(\left(1.85032828h_k\left(1-\frac{q_k}{e_k}\right)\right)^2, \left(1.85032828h_k\left(1+\frac{q_k}{e_k}\right)\right)^2\right) \quad (4.8)$$

and

$$J_k = \left(\left(\frac{\pi}{\sqrt{6}\ln 2}h_k\left(1-\frac{q_k}{e_k}\right)\right)^2 + N_{cit} - N_{cit}^{(k)}, \left(\frac{\pi}{\sqrt{6}\ln 2}h_k\left(1+\frac{q_k}{e_k}\right)\right)^2 + N_{cit} - N_{cit}^{(k)}\right)$$

$$= \left(\left(1.85032828h_k\left(1-\frac{q_k}{e_k}\right)\right)^2 + N_{cit} - N_{cit}^{(k)}, \left(1.85032828h_k\left(1+\frac{q_k}{e_k}\right)\right)^2 + N_{cit} - N_{cit}^{(k)}\right). \quad (4.9)$$

**Definition 4.4.** Under notions and notations of Definition 4.2, we define the $h$-*defect of asymptotic normality of citations* $cit_1, cit_2, ..., cit_p$, denoted as $d$ ($0 \leq d \leq h$), as follows.

1) If $e_k \geq h_k$ for all $k = 0, 1, ..., h$, then $d = 0$.
2) If $e_0 := e \geq h$ and $e_k < h_k$ for some $k \in \{1, ..., h\}$, then $d$ ($1 \leq d \leq h$) is the smallest integer $k \in \{1, ..., h-1\}$ such that $e_i \geq h_i$ for all $i = 0, 1, ..., k$ and $e_{k+1} < h_{k+1}$.
3) If $e_0 := e < h$, we put $d = 0$.
4) If $e_0 := e \leq h$ and $e_k > h_k$ for some $k \in \{1, ..., h\}$, then $d$ ($1 \leq d \leq h$) is the smallest integer $k \in \{1, ..., h-1\}$ such that $e_i \leq h_i$ for all $i = 0, 1, ..., k$ and $e_{k+1} > h_{k+1}$.

Following Definition 4.4, we will call the multiset $D := \{cit_1, cit_2, ..., cit_d\}$ the $h$-*defect-core of asymptotic normality of citations* $cit_1, cit_2, ..., cit_p$ and we will call the multiset $A := \{cit_{d+1}, cit_{d+2}, ..., cit_{d+h_d}\}$ the *domain of asymptotic normality of citations* $cit_1, cit_2, ..., cit_p$. Notice that if $d = 0$, then $D$ is the empty set.

Our computational results given in Subsection 5.2 and some additional computations that are not presented here, show that that the case 2) of Definition 4.4. refers to the majority of scientists whose $h$-indices are greater than equal to 15.

**Definition 4.5.** Under notations of Definitions 4.3 and 4.4, we define

$$\bar{I}_d = \frac{a_d + b_d}{2}, \quad \bar{I}_{d+1} = \frac{a_{d+1} + b_{d+1}}{2}, \quad A' = \frac{\bar{I}_d + \bar{I}_{d+1}}{2} \quad (4.10)$$

and

$$\bar{J}_d = \frac{c_d + f_d}{2}, \quad \bar{J}_{d+1} = \frac{c_{d+1} + f_{d+1}}{2}, \quad A = \frac{\bar{J}_d + \bar{J}_{d+1}}{2}. \quad (4.11)$$

**Definition 4.6.** Considering four cases of Definition 4.4, we define the numbers $\alpha_d$, $\beta_d$, $\alpha_{d+1}$, $\beta_{d+1}$ ($0 \leq d \leq p-1$) and the estimations $B'$, $B''$ and $B$ of $N_{cit}$ as follows.

1a) If $e > h+1$ in the case 1) of Definition 4.4, we define

$$B = \frac{1}{2}\left(\left(\frac{\pi}{\sqrt{6}\ln 2}h\left(1+\frac{q}{e}\right)\right)^2 + \left(\frac{\pi}{\sqrt{6}\ln 2}h_1\left(1+\frac{q_1}{e_1}\right)\right)^2\right). \tag{4.12}$$

1b) If $h \leq e \leq h+1$ in the case 1) of Definition 4.4, we define $\beta_0 = e - \lfloor e \rfloor$, $\alpha_0 = 1 - \beta_0$ and

$$B = \frac{1}{2}\left(\alpha_0\left(\frac{\pi}{\sqrt{6}\ln 2}h\left(1-\frac{q}{e}\right)\right)^2 + \beta_0\left(\frac{\pi}{\sqrt{6}\ln 2}h\left(1+\frac{q}{e}\right)\right)^2\right) \tag{4.13}$$

2) In the case 2) of Definition 4.4 it holds $e_d \geq h_d$ and $e_{d+1} < h_{d+1}$. Then we consider the following four subcases.

2a) If $e_d > h_d + 1$, we define $\beta_d = 1$, $\alpha_d = 0$ and

$$B' = \left(\frac{\pi}{\sqrt{6}\ln 2}h\left(1+\frac{q_d}{e_d}\right)\right)^2 + N_{cit} - N_{cit}^{(d)} = \left(\frac{\pi}{\sqrt{6}\ln 2}h\left(1+\frac{q_d}{e_d}\right)\right)^2 + \sum_{i=1}^{d}cit_i. \tag{4.14}$$

2b) If $e_d \leq h_d + 1$, we define $\beta_d = e_d - \lfloor e_d \rfloor$, $\alpha_d = 1 - \beta_d$ and

$$B' = \frac{1}{2}\left(\alpha_d\left(\frac{\pi}{\sqrt{6}\ln 2}h\left(1-\frac{q_d}{e_d}\right)\right)^2 + \beta_d\left(\frac{\pi}{\sqrt{6}\ln 2}h\left(1+\frac{q_d}{e_d}\right)\right)^2\right) + N_{cit} - N_{cit}^{(d)}$$

$$= \frac{1}{2}\left(\alpha_d\left(\frac{\pi}{\sqrt{6}\ln 2}h\left(1-\frac{q_d}{e_d}\right)\right)^2 + \beta_d\left(\frac{\pi}{\sqrt{6}\ln 2}h\left(1+\frac{q_d}{e_d}\right)\right)^2\right) + \sum_{i=1}^{d}cit_i. \tag{4.15}$$

2c) If $h_{d+1} > e_{d+1} + 1$, we define $\beta_{d+1} = 0$, $\alpha_{d+1} = 1$ and

$$B'' = \left(\frac{\pi}{\sqrt{6}\ln 2}h\left(1-\frac{q_{d+1}}{e_{d+1}}\right)\right)^2 + N_{cit} - N_{cit}^{(d+1)}$$

$$= \left(\frac{\pi}{\sqrt{6}\ln 2}h\left(1-\frac{q_{d+1}}{e_{d+1}}\right)\right)^2 + \sum_{i=1}^{d+1}cit_i. \tag{4.16}$$

2d) If $h_{d+1} \leq e_{d+1} + 1$, we define $\beta_{d+1} = e_{d+1} - \lfloor e_{d+1} \rfloor$, $\alpha_{d+1} = 1 - \beta_{d+1}$ and

$$B'' = \frac{1}{2}\left(\beta_{d+1}\left(\frac{\pi}{\sqrt{6}\ln 2}h\left(1-\frac{q_{d+1}}{e_{d+1}}\right)\right)^2 + \alpha_{d+1}\left(\frac{\pi}{\sqrt{6}\ln 2}h\left(1+\frac{q_d}{e_d}\right)\right)^2\right) + N_{cit} - N_{cit}^{(d+1)} \qquad (4.17)$$

$$= \frac{1}{2}\left(\beta_{d+1}\left(\frac{\pi}{\sqrt{6}\ln 2}h\left(1-\frac{q_{d+1}}{e_{d+1}}\right)\right)^2 + \alpha_{d+1}\left(\frac{\pi}{\sqrt{6}\ln 2}h\left(1+\frac{q_d}{e_d}\right)\right)^2\right) + \sum_{i`=1}^{d+1} cit_i.$$

Then in all subcases 2a) - 2d) of Case 2) we put

$$B = \frac{B'+B''}{2}. \qquad (4.18)$$

3) In the case 3) of Definition 4.4 we consider the following two subcases.

3a) If $h - e \leq 1$, then $d = 0$ and we define $\beta_0 = e - \lfloor e \rfloor$, $\alpha_0 = 1 - \beta_0$ and

$$B = \frac{1}{2}\left(\beta_0\left(\frac{\pi}{\sqrt{6}\ln 2}h\left(1-\frac{q}{e}\right)\right)^2 + \alpha_0\left(\frac{\pi}{\sqrt{6}\ln 2}h\left(1+\frac{q}{e}\right)\right)^2\right). \qquad (4.19)$$

3b) If $h - e > 1$, then $d = 0$ and we define $\beta_0 = e - \lfloor e \rfloor$ and

$$B = \beta_0\left(\frac{\pi}{\sqrt{6}\ln 2}h\left(1-\frac{q}{e}\right)\right)^2. \qquad (4.20)$$

4) In the case 4) of Definition 4.4, we define the values $B'$ and $B''$ by the expressions (4.14)-(4.18) in which $\beta_d$ is replaced by $\alpha_d$, $\alpha_d$ is replaced by $\beta_d$, $\beta_{d+1}$ is replaced by $\alpha_{d+1}$ and $\alpha_{d+1}$ is replaced by $\beta_{d+1}$.

**Proposition 4.7.** *The value $A$ defined by* (4.11) *of Definition 4.5 is given by the expression*

$$A = \frac{\overline{J}_d + \overline{J}_{d+1}}{2} = \frac{\pi^2}{12\ln^2 2}\left(h_d^2\left(1+\left(\frac{q_d}{e_d}\right)^2\right) + h_{d+1}^2\left(1+\left(\frac{q_{d+1}}{e_{d+1}}\right)^2\right)\right) + \sum_{i`=1}^{d} cit_i + \frac{cit_{d+1}}{2}. \qquad (4.21)$$

*Proof.* The expression (4.21) immediately follows from the formula (4.9) of Definition 4.3 and the formula (4.4) of Definition 4.2.

Notice that in the case 2) of Definition 4.4 which refers to the most researchers in all databases, the numerous computational results show that $0 \leq h_d - e_d \leq 2$ and $0 \leq e_{d+1} - h_{d+1} \leq 2$, and so we can approximately assume that $e_d \approx h_d$, $e_{d+1} \approx h_{d+1}$, and hence, $q_d = 2N_h^{(d)}/h_d^2 - 1 \approx 4h_d^2/h_d^2 - 1 = 3$, $q_d/e_d \approx 3/h_d$, $q_{d+1} = 2N_h^{(d+1)}/h_{d+1}^2 - 1 \approx 4h_{d+1}^2/h_{d+1}^2 - 1 = 3$ and $q_{d+1}/e_{d+1} \approx 3/h_{d+1}$. Substituting this into (4.20), we obtain

$$A \approx \frac{\pi^2}{12\ln^2 2}\left(h_d^2 + h_{d+1}^2\right) + 30.8134 + \sum_{i=1}^{d+1} cit_i - \frac{cit_{d+2}}{2} = 1.71186\left(h_d^2 + h_{d+1}^2\right) + 30.8134 + \sum_{i=1}^{d} cit_i + \frac{cit_{d+1}}{2}. \quad (4.22)$$

Since our computational results given in Section 5 show that the value $A$ defined by (4.11) and the values $B$ defined by (4.12), (4.13) and (4.18) are close to $N_{cit}$, i.e., $A \approx N_{cit}$ and $B \approx N_{cit}$, we propose these two values as good approximations of $N_{cit}$ (see Tables 2 and 5). Related computational results for 14 experts in Bibliometrics and Infometrics (Tables 1-4 and Table 8) and 27 researcher Google Scholar's and Scopus' profiles, including 24 researchers (Tables 5-7 and Table 8) are given in the next section.

Put $N_{cit}(h) = (x+1)h^2$ with a real number $x \geq 0$. Then $q = 2N_{cit}(h)/h^2 - 1 = 2x+1$, $e = h\sqrt{x}$. Define the function $f$ as $f(x) := hq/e = (2x+1)/\sqrt{x}$ for $x \in (0, +\infty)$. The function $f$ attains the minimum value $2\sqrt{2}$ at $x = 1/2$ and so, $N_{cit}(h) = 2h^2$, $e = h\sqrt{2}$ and hence, $q/e \geq 2\sqrt{2}/h$ for each $x > 0$. Note also that by the inequality $a + b \leq \sqrt{2(a^2 + b^2)}$ ($a, b \geq 0$) it follows that $h + e = h + \sqrt{N_{cit}(h) - h^2} \leq \sqrt{2N_{cit}(h)}$ and equality holds when $N_{cit}(h) = 2h^2$, and then $e = h$.

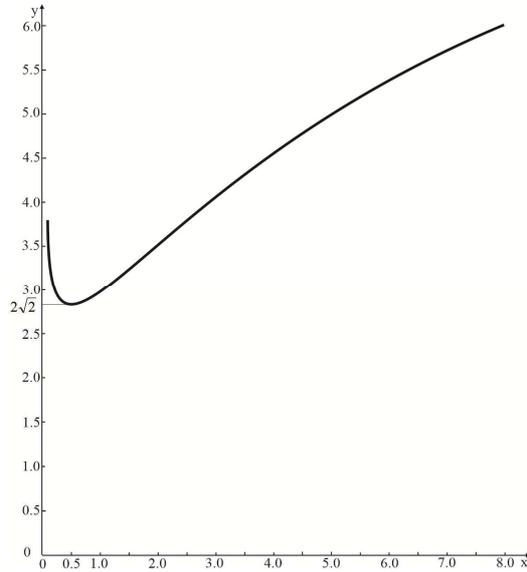

**Figure 1.** The graphic of the function $f(x) = \dfrac{2x+1}{\sqrt{x}}$ ($x > 0$).

Now under above notations, we define $q' = N_{cit}(h)/h^2 = x+1$. Note that $q' \leq q$ and equality holds if and only if $cit_1 = cit_2 = \cdots = cit_h = h$. Define the function $g$ as $g(x) := hq'/e = (x+1)/\sqrt{x}$ for $x \in (0, +\infty)$. The function $g$ attains the minimum value 2 at $x=1$, and then $N_{cit}(h) = 2h^2$, $e = h$ and $q'/e \geq 2/h$ for each $x > 0$.

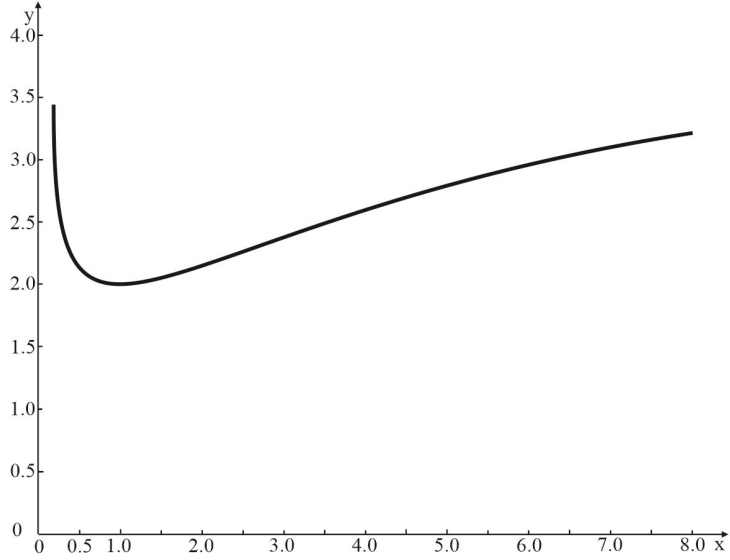

**Figure 2.** The graphic of the function $g(x) = \dfrac{x+1}{\sqrt{x}}$ ($x > 0$).

Observe that $f(x)/g(x) = (2x+1)/(x+1) = 1 + x/(x+1)$, and hence, $1 < f(x)/g(x) < 2$ for all $x > 0$, while $f(x)/g(x) \to 2$ as $x \to \infty$.

**Proposition 4.7.** *For each $k \in \{0, 1, ..., p-1\}$ we have*

$$h_k - 1 \leq h_{k+1} \leq h_k. \tag{4.23}$$

*Furthermore, $h_{k+1} = h_k$ if and only if $cit_{k+1+h_k} = h_k$, and $h_{k+1} = h_k - 1$ if and only if $cit_{k+1+h_k} < h_k$.*

*Proof.* The inequalities (4.23) immediately follow from Definition 4.2.

If $cit_{k+1+h_k} = h_k$, then by Definition 4.2, $h_{k+1} = h_k$. Conversely, if $h_{k+1} = h_k$, then by Definition 4.2 must be $cit_{k+1+h_k} \geq h_k$. If we suppose that $cit_{k+1+h_k} \geq h_k + 1$, it follows that $h_k \geq h_k + 1$. This contradiction shows that $cit_{k+1+h_k} = h_k$.

If $cit_{k+1+h_k} < h_k$, then by Definition 4.2, $h_{k+1} = h_k - 1$. Conversely, if $h_{k+1} = h_k - 1$, then by Definition 4.2 must be $cit_{k+1+h_k} < h_k$. This completes the proof.

**Proposition 4.8.** *Under above notations, for each $k \in \{0, 1, ..., p-1\}$, the following three assertions are true.*

(i) *If $h_{k+1} = h_k$ (or equivalently, by Proposition 4.7, if and only if $cit_{k+1+h_k} = h_k$), then $e_k \geq e_{k+1}$.*

(ii) *If $h_{k+1} = h_k - 1$ and $cit_{k+1} \geq 2h_k - 1$, then $e_k \geq e_{k+1}$.*

(iii) If $h_{k+1} = h_k - 1$ and $cit_{k+1} < 2h_k - 1$, then $e_k < e_{k+1}$.

*Proof.* (i) Suppose that $h_{k+1} = h_k$. Then by Proposition 4.7, $cit_{k+1+h_k} = h_k$ and hence,

$$e_{k+1} = \sqrt{N_h^{(k+1)} - (h_{k+1})^2} = \sqrt{\sum_{j=k+2}^{h_{k+1}+k+1} cit_j - (h_{k+1})^2} = \sqrt{\sum_{j=k+2}^{h_k+k+1} cit_j - (h_k)^2}$$

$$= \sqrt{\left(\sum_{j=k+1}^{h_k+k} cit_j - (h_k)^2\right) - (cit_{k+1} - cit_{h_k+k+1})} = \sqrt{(e_k)^2 - (cit_{k+1} - cit_{h_k+k+1})} \leq e_k,$$

which proves the assertion (i)

*Proof* of (ii) and (iii). Assume that $h_{k+1} = h_k - 1$. Then

$$e_{k+1} = \sqrt{N_h^{(k+1)} - (h_{k+1})^2} = \sqrt{\sum_{j=k+2}^{h_k+k} cit_j - (h_k - 1)^2} = \sqrt{\left(\sum_{j=k+1}^{h_k+k} cit_j - (h_k)^2\right) - (cit_{k+1} - 2h_k + 1)}$$

$$= \sqrt{(e_k)^2 - (cit_{k+1} - 2h_k + 1)}.$$

Hence, if $cit_{k+1} \geq 2h_k - 1$, then $e_k \geq e_{k+1}$. Conversely, if $cit_{k+1} < 2h_k - 1$, then $e_k < e_{k+1}$. This completes the proof of proposition.

As an immediate consequence of Proposition 4.8, we obtain the following result.

**Proposition 4.9.** *Suppose that $e_k \geq h_k$ for some $k \in \{0,1,...,p-1\}$. Then the following assertions are true.*

(i) *If $h_{k+1} = h_k$ (or equivalently, by Proposition 4.7, if $cit_{k+1+h_k} = h_k$), then $e_{k+1} < h_{k+1}$ if and only if $e_k < \sqrt{(h_k)^2 - h_k + cit_{k+1} - 1}$.*

(ii) *If $h_{k+1} = h_k - 1$ (or equivalently, by Proposition 4.7, if $cit_{k+1+h_k} < h_k$), then $e_{k+1} < h_{k+1}$ if and only if $e_k \leq \sqrt{(h_k)^2 - 4h_k - 1 + cit_{k+1}}$.*

*Proof.* (i) Suppose that $h_{k+1} = h_k$. Then from the proof of (i) of Proposition 4.8 and the fact that by Proposition 4.7, $cit_{k+1+h_k} = h_k$, we obtain

$$e_{k+1} = \sqrt{(e_k)^2 - (cit_{k+1} - cit_{h_k+k+1})},$$

whence it follows that the inequality $e_{k+1} < h_{k+1} = h_k$ is equivalent to the following one:

$$(e_k)^2 \leq (h_k)^2 - h_k + cit_{k+1} - 1,$$

i.e.,

$$e_k \leq \sqrt{(h_k)^2 - h_k + cit_{k+1} - 1} = \sqrt{(cit_{k+1+h_k})^2 - cit_{k+1+h_k} + cit_{k+1} - 1}.$$

(ii) Assume that $h_{k+1} = h_k - 1$. Then from the proofs of (ii) and (iii) of Proposition 4.8 we find that

$$e_{k+1} = \sqrt{(e_k)^2 - cit_{k+1} + 2h_k - 1},$$

whence it follows that the inequality $e_{k+1} < h_{k+1} = h_k - 1$ is equivalent to the following one:

$$(e_k)^2 \leq (h_k)^2 - 4h_k + cit_{k+1},$$

i.e.,

$$e_k \leq \sqrt{(h_k)^2 - 4h_k - 1 + cit_{k+1}} = \sqrt{(cit_{k+1+h_k})^2 - 4cit_{k+1+h_k} - 1 + cit_{k+1}}.$$

This completes the proof of proposition.

The "converse analogue" of Proposition 4.9 is given as follows. Since his proof is quite similar to those of of Proposition 4.9, it can be omitted.

**Proposition 4.10.** *Suppose that $e_k < h_k$ for some $k \in \{0,1,...,p-1\}$. Then the following assertions are true.*

(i) *If $h_{k+1} = h_k$ (or equivalently, by Proposition 4.7 if $cit_{k+1+h_k} = h_k$), then $e_{k+1} \geq h_{k+1}$ if and only if $e_k \geq \sqrt{(h_k)^2 - h_k + cit_{k+1}}$.*

(ii) *If $h_{k+1} = h_k - 1$ (or equivalently, by Proposition 4.5, if $cit_{k+1+h_k} < h_k$), then $e_{k+1} \geq h_{k+1}$ if and only if $e_k \geq \sqrt{(h_k)^2 - 4h_k - 2 + cit_{k+1}}$.*

**Remarks 4.11.** From Propositions 4.9 and 4.10 we can get the recursive computational algorithm for the determination of $h$-defect of asymptotic normality of citations $cit_1, cit_2,...,cit_p$ ($d$), defined by Definition 4.4. Namely, if $e \geq h$ (which is equivalent to $N_{cit}(h) := N_h^{(0)} \geq 2h^2$) and $e_{k+1} < h_{k+1}$ for at least one $k \in \{0,1,...,h-1\}$, then by the Case 2) of Definition 4.4, $d$ ($0 \leq d \leq h$) is the smallest integer $k$ such that $e_i \geq h_i$ for all $i = 0,1,...,k$ and $e_{k+1} < h_{k+1}$. Such a value $d$ is the smallest positive integer $k$ satisfying the estimates (i) and (ii) of Proposition 4.8. Notice that the Case 2) of Definition 4.4 is satisfied for 35 researcher profiles at Scopus and Google Scholar databases which are considered in our computational results (13 among 14 researchers in Table 2 and 22 among 27 researchers in Table 5), except six researchers concerning the cases 3a) and 3b) (namely, one researcher from Table 2 and five researchers from Table 5).

The all values from Tables 2 and 5 in Subsection 5.2 can be calculated using the following two recurrence relations (see Propositions 4.9 and 4.10):

$$h_0 = h, N_h^{(0)} = N_{cit}(h); h_{k+1} = h_k \text{ if } cit_{h_k+k+1} = h_k \text{ and } h_{k+1} = h_k - 1 \text{ if } cit_{h_k+k+1} < h_k \quad (4.24)$$

and

$$N_h^{(k+1)} = N_h^{(k)} - cit_{k+1} + \delta_k cit_{h_k+k+1}, \text{ where } \delta_k = 1 \text{ if } cit_{h_k+k+1} = h_k \text{ and } \delta_k = 0 \text{ if } cit_{h_k+k+1} < h_k. \quad (4.25)$$

**Remarks 4.12.** S. Brown (2018) in Abstract of his paper wrote "The $h$-index has been considered in terms of both unrestricted integer partitions and fuzzy integrals, and the expected value of $h$ for a given number of citations has been estimated. However, the distribution of $h$ as a function of both the number of citations and the number of papers has not been considered explicitly. Using Durfee squares determined from restricted integer partitions, it is shown that for a small number of papers the expected value of $h$-index estimated from the unrestricted partitions of the number of citations is unreliable. Despite this, it is confirmed that the distribution of his asymptotically normal. This means that $h$-indices should be considered in the context of the number of publications unless that number is large."

Brown (2018) pointed out that for any combination of $p$ published papers and $N_{cit}$ citations, a range of values of the $h$-index is possible. Brown (2018) noticed that if, as is often the case, $p$ is much smaller than $N_{cit}$, then the use of unrestricted partitions ($p \leq N_{cit}$) to estimate the distribution of $h$ is inappropriate and this discrepancy grows as ($N_{cit}/p$) increases. Brown (2018) proposed the empirical approximation of the 95% confidence interval for the $h$-index as

$$(0.54\sqrt{N_{cit}} - 1.96(0.57 + 0.045\sqrt{N_{cit}}), 0.54\sqrt{N_{cit}} + 1.96(0.57 + 0.045\sqrt{N_{cit}})). \quad (4.26)$$

Brown (2018) observed that for small value of $N_{cit}$, the distribution of $h$ is not normal, but as $N_{cit}$ increases (for moderate $p = 100$) from 200 a normal distribution provides an increasingly good approximation (based on the Anderson-Darling normality test). The distributions of $h$ for $N_{cit} = 100$ with $p = 6, 8, 10, 20$ and for $N_{cit} = 200, 400, 600, 1000$ with $p = 100$ are presented in Figures 5 and 6 in Brown (2018), respectively.

Note that the confidence interval given by (4.26) is used for giving our computational results presented in Table 8.

**Remarks 4.13.** Yong (2014) observed that the $h$-index is not a good predictor of winners of the Fields Medal in Mathematics. In particular, Yong (2014) considered Terence Tao's Google Scholar profile. Since he has 30053 citations, the "rule of thumb" predicts his $h$-index is 93.6. This is far from his actual $h$-index of 65 (cf. with recent Tao's Google Scholar's and Scopus' indices, which are equal $h = 106$, $h_{NA} = 162.3$ in Google Scholar database and $h = 70$ and $h_{NA} = 116.981$ in Scopus database, which is presented in Table 5 of Subsection 5.2 and in Tables A2 and A4 in Appendix). Notice that the ratios of these three pairs of Tao's indices are 93.6/65=1.44, 162.3/106=1.531 and 116.981/70=1.67, respectively. Yong (2014) observed that his top five citations (joint with E. Candes on compressed sensing) are applied. Removing the papers on this topic leaves 13942 citations his new estimate is therefore $h_{NA} = 63.7$ and his revised $h$-index is 61.16, and so $h_{NA}/h = 63.7/61.16 = 1.042$. Notice that Tao's Google Scholar's and Scopus' $h$-defects of asymptotic normality of his citations (by Definition 4.4) are recently equal to 14 and 7, respectively (see Table 5).

But in a subfield of physics, Lehmann et al. (2006) concluded that the $h$-index was inferior to the mean number of citations per paper as indicator of scientific quality.

Recently, Mahmoudi et al. (2021) proposed a statistical approach to model the $h$-index based on the total number of citations and the duration from the publishing of the first article. Namely, this research deals with a practical approach to model the $h$-index based on the total number of citations ($N_{cit}$) and the duration from the publishing of the first article ($D_1$). The obtained computational results indicated that both $N_{cit}$ and $D_1$ had a significant effect on $h$ (with $p < 0.001$) and the authors of this paper estimated the $h$-index as a function of $N_{cit}$ and $D_1$ by the following approximate estimates: $h \approx h(N_{cit}) = 0.600(N_{cit})^{0.476}$ and $h \approx h(D_1) \approx 0.667(D_1)^{1.041}$, which is similar to the approximate estimates (3.3), $h \approx h(N_{cit}) = 0.544(N_{cit})^{0.5}$ proved by Canfield, Corteel and Savage (1998). Using multiple nonlinear regression, Mahmoudi et al. (2021) also estimated the $h$-index as a function of two variables $N_{cit}$ and $D_1$.

Notice also that by Glänzel-Schubert model (Schubert and Glänzel, 2007) there is a strong correlation between the $h$-index and $(N_{cit}/N_p)^{a/(1+a)}$, where $N_p$ is the total number of publications and the constant $a$ is related to the Lotka exponent in the Lotkaian framework (see Radicchi and Castellano, 2013 and Schubert and Glänzel, 2007). Furthermore, Radicchi and Castellano (2013) proposed the approximations $h \approx (N_{cit})^{0.42}$ and $h \approx (N_{cit})^{0.41}(N_p)^{0.18}$. Spruit (2012) proposed the approximation $h \approx 0.5(\sqrt{N_{cit}} + 1)$. Considering the dataset of 255 physicists, Redner (2012) proposed the approximation $\sqrt{N_{cit}}/(2h) \approx 1.045$. Furthermore, considering the dataset of 147 university chemistry research groups in the Netherlands, van Raan (2005) found the approximation $h \approx 0.42(N_{cit})^{0.45}$.

## 5. Computational results and discussion

In order to estimate the total number of citations of a researcher using his $h$-index and his $h$-core, we propose the methods given in Subsection 5.1. In Subsection 5.2 we refine these methods by using additionally a relatively small number of citations from the tail. Computational results are presented for 14 leading experts in Scientometrics (Tables 1-4) and for some other 27 scientists (Tables 5-7).

The researchers whose bibliometric indicators and indices are calculated in Subsections 5.1 and 5.2, using their Scopus (S) and Google Scholar profiles (GS) are:

1) (Tables 1-4 and Table 8). Scopus profiles of 14 Price awardees/leading experts in Scientometrics (February 21 2023).

Leot Leydesdorff (S), Wolfgang Glänzel (S), Henk Moed (S), Anthony Van Raan (S), Ronald Rousseau (S), Andras Schubert (S), Ben Martin (S), Francis Narin (S), Eugene Garfield (S), Tibor Braun (S), Henry Small (S), Leo Egghe (S), Peter Ingwersen (S), Howard D. White (S).

2) (Tables 5-7 and Table 8). 27 research profiles at Google Scholar (GS) and Scopus (S) databases of 24 researchers from various scientific fields, on April-June 2023 (which include 24 different researchers, with being considered both Scopus and Google Scholar profiles for Paul Erdös, Terence Tao and George Andrews).

Sigmund Freud (GS, June 26, Neurology, Psyhology, Psyhopatology), H. Jihn Kim (GS, May 26, High energy physics, astro-particle physics, scintillator), Ronald C. Kessler (GS, May 24, Psychiatric Epidemiology), Albert Einstein (GS, May 21, Physics), Paul Erdös (GS, April 20, Number Theory, Combinatorics, Probability, Set Theory, Mathematical Analysis), Terence Tao (GS, April 20, Mathematical analysis, Combinatorics, Random Matrix Theory, PDE, Number theory), L. Leydesdorf (GS, May 24, Bibliometrics and Scientometrics), Christian G. Meyer (GS, Infectious diseases), Terence Tao (S, May 20, Mathematical analysis, Combinatorics, Random matrix theory, PDE, Number theory), George Andrews (GS, May 05, Mathematical analysis, Combinatorics, Number theory), Michael McAleer (GS, May 27, Economics, Econometrics, Financial economics, Financial econometrics), Jorge E. Hirsch (S, May 23, Physics, Supercodutivity, Ferromagnetism), Paul Erdös (S, April 21, Number theory, Combinatorics, Probability, Set theory, Mathematical analysis), Frank. T. Edelmann (May 27, GS, Inorganic chemistry), Carl Friedrich Gauss (GS, May 28, Mathematics), George Andrews (S, May 23, Mathematical analysis, Combinatorics, Number theory), Stratos Papadimitrou (GS, May 27, Transportation, Maritime logistics), Doron Zeilberger (S, May 23, Combinatorics, Hypergeometric identities, Symbolic computation), Irena Orović (GS, May 27, Signal processing, Multimedia signals and systems, Compressive sensing, Time-frequency signal analysis), Carla. D. Savage (GS, June 3, Algorithms, Combinatorics), Velibor Spalević (GS, May 25, Sediment yield, Soil erosion modeling, Erosion, Watershed management, IntErO), Parisa Ziarati (GS, June 26, Toxicology, Green chemistry, Enviromental chemistry, Inorganic chemistry), Alexander Yong (GS, May 23, Combinatorics, Algebraic geometry, Representation theory), David Kalaj, (GS, May 23, Complex analysis), Petar Vukoslavčević (GS, May 28, Mechanical engineering, Measurement, Turbulent flow, Boundary layer, Applied Mathematics), Katarina Monkova (S, May 23, Mechanical engineering), Ljuben R. Mutafchiev (GS, June 3, Random combinatorial structures).

Notice that H. Jihn Kim (*h*= **338**, GS, High energy physics, Astro-particle physics, Scintillator) and Ronald C. Kessler (*h*= **328**, GS, Psychiatric Epidemiology) (Table 5) have recently the highest $h$-index in the Google Scholar database. Moreover, in the Google Scholar database Sigmund Freud (Neurology, Psyhology, Psyhopatology) has the $h$-index equal to **288** and the total number of his citations is **643730**.

### 5.1. The estimates of total number of citations of a researcher using only his $h$-index and his $h$-core

We recall the following notions, notations and expressions presented in Table 1.

$p$-the total number of publications; $N_{cit} = \sum_{j=1}^{p} cit_j$, - the total number of citations of a considered researcher; the $h$-index (the expression (1.2); $N_{cit}(h) = \sum_{j=1}^{h} cit_j$ - the total number of citations into $h$-core; $h_{NA} = 0.5404446 \sqrt{N_{cit}}$ - the normal approximation of $h$-index (the expression (3.3)); $r = \lfloor 2N_{cit}(h)/(h(h+1)) \rfloor - 1$ (the expression (1.11)); $e = \sqrt{N_{cit}(h) - h^2} = \sqrt{\sum_{j=1}^{h} cit_j - h^2}$ (the expression (1.7)); $q = 2N_{cit}(h)/h^2 - 1$ (the expression (4.6) with $k = 0$); the interval $I = I_0 = (a,b) = \left((1.85032828(1 - q/e)h)^2, (1.85032828(1 + q/e)h)^2\right)$ (the expression (4.8) with $k = 0$). Then the mean of the interval $I$ is equal to

$$\overline{I} = 3.42371438\left(1 + \left(\frac{q}{e}\right)^2\right)h^2. \tag{5.1}$$

Here we also consider the interval $I(q)$ defined as

$$I(q) = (a',b') = \left(\left(1.85032828\left(1 - \frac{q}{e}\right)(h+q)\right)^2, \left(1.85032828\left(1 + \frac{q}{e}\right)(h+q)\right)^2\right), \tag{5.2}$$

whose mean is equal to

$$\overline{I(q)} = 3.42371438\left(1 + (q/e)^2\right)(h+q)^2, \tag{5.3}$$

and the interval $I(r)$ defined as

$$I(r) = (a'',b'') = \left(\left(3.42371438\left(1 - \frac{r}{e}\right)(h+r)\right)^2, \left(3.42371438\left(1 + \frac{r}{e}\right)(h+r)\right)^2\right), \tag{5.4}$$

whose mean is equal to

$$\overline{I(r)} = 3.42371438\left(1 + (r/e)^2\right)(h+r)^2. \tag{5.5}$$

Recall that after Figure 1 in Section 4 we defined the ratio $q'$ as

$$q' = \frac{N_{cit}(h)}{h^2}. \tag{5.6}$$

Note that $q' \leq q$.

Define the shorter interval $I(q')$ than the interval $I(q)$ as

$$I(q') = (a'',b'') = \left(\left(1.85032828\left(1 - \frac{q'}{e}\right)(h+q')\right)^2, \left(1.85032828\left(1 + \frac{q'}{e}\right)(h+q')\right)^2\right), \tag{5.7}$$

whose mean is equal to

$$\overline{I(q')} = 3.42371438\left(1 + (\frac{q'}{e})^2\right)(h+q')^2. \tag{5.8}$$

We also denote the relative errors of $\overline{I}$, $\overline{I(r)}$, $\overline{I(q)}$ and $\overline{I(q')}$ as $\delta_1 = (\overline{I} - N_{cit})/N_{cit}$, $\delta_2 = (\overline{I(r)} - N_{cit})/N_{cit}$, $\delta_3 = (\overline{I(q)} - N_{it})/N_{cit}$ and $\delta_4 = (\overline{I(q')} - N_{cit})/N_{cit}$, respectively.

Recall that the all calculations in this paper are made using the software *Wolfram Mathematica 11*.

Using data of Table A1 from Appendix (the reduced Table A2 in Meštrović and Dragović, 2023) and the above notions, notations and the expressions, we obtain the following Table 1.

**Table 1**. Price awardees data (bases on Scopus, February 21 2023). The bibliometric indicators and indices.

| Researcher | L. Leydesdorff | W. Glänzel | H.F. Moed | A.F.J. Van Raan | R. Rousseau | A. Schubert | B. Martin |
|---|---|---|---|---|---|---|---|
| $N_p^+$ | 406 | 258 | 127 | 124 | 295 | 141 | 85 |
| $N_{cit}$ | 25005 | 11766 | 7606 | 8308 | 8053 | 7587 | 7598 |
| $h$ | 79 | 61 | 49 | 48 | 43 | 42 | 38 |
| $h/N_p^+$ | 0.195 | 0.236 | 0.386 | 0.387 | 0.146 | 0.298 | 0.447 |
| $N_{cit}(h)$ | 17360 | 8049 | 6351 | 6833 | 5203 | 6359 | 7048 |
| $N_{cit}(h)/h$ | 219.747 | 131.95 | 129.612 | 142.235 | 121 | 151.140 | 185.48 |
| $N_{cit}(h)/N_{cit}$ | 0.694 | 0.684 | 0.835 | 0.822 | 0.646 | 0.838 | 0.928 |
| $h_{NA}$ | 85.461 | 58.623 | 47.133 | 49.260 | 48.499 | 47.074 | 47.108 |
| $h_{NA}/h$ | **1.0818** | **0.9610** | **0.9619** | **1.0263** | **1.1279** | **1.1208** | **1.2397** |
| $q$ | 4.563 | 3.326 | 4.290 | 4.931 | 4.628 | 6.210 | 8.762 |
| $e$ | 105.447 | 65.788 | 62.849 | 67.298 | 57.914 | 67.786 | 74.860 |
| $e/h$ | 1.335 | 1.079 | 1.283 | 1.402 | 1.347 | 1.614 | 1.970 |
| $q/e$ | 0.043 | 0.051 | 0.0683 | 0.073 | 0.080 | 0.093 | 0.130 |
| $I$ | (19558.16; 23256.68) | (11484.07; 14060.34) | (7136.42; 9380.86) | (6774.63; 9086.55) | (5359.12; 7382.63) | (4983.55; 7196.69) | (3854.27; 6168.88) |
| $N_{cit} \in I$ | No | Yes | Yes | Yes | No | No | No |
| $\overline{I}$ | 21407.42 | 12772.20 | 8258.64 | 7930.59 | 6370.88 | 6090.12 | 5011.58 |
| $\delta_1$ | **-0.144** | **0.086** | **0.086** | **0.045** | **-0.2088** | **-0.197** | **-0.340** |
| $h_q := h+q$ | 83.563 | 64.326 | 53.290 | 52.931 | 47.628 | 48.210 | 46.762 |
| $h_{NA}/h_q$ | **1.023** | **0.911** | **0.884** | **0.931** | **1.018** | **0.976** | **1.007** |
| $I(q)$ | (21882.74; 26020.85) | (12770.53; 15635.41) | (8440.72; 11095.38) | (8238.03; 11049.35) | (6574.78; 9057.30) | (6566.21; 9482.18) | (5836.61; 9341.68) |
| $N_{cit} \in I(q)$ | Yes | No | No | Yes | Yes | Yes | Yes |
| $\overline{I(q)}$ | 23951.80 | 14202.97 | 9768.05 | 9643.69 | 7816.04 | 8024.20 | 7589.15 |
| $\delta_2$ | **-0.070** | **0.207** | **0.284** | **0.161** | **-0.029** | **0.058** | **-0.001** |
| $r$ | 4 | 3 | 4 | 4 | 4 | 6 | 8 |
| $h+r$ | 83 | 64 | 53 | 52 | 47 | 48 | 46 |
| $h_{NA}/(h+r)$ | 1.030 | 0.916 | 0.889 | 0.947 | 1.032 | 0.981 | 1.024 |
| $r/e$ | 0.038 | 0.046 | 0.064 | 0.060 | 0.069 | 0.089 | 0.107 |
| $I(r)$ | (21827.50; 25412.56) | (12763.04; 15343.37) | (8425.60; 10887.61) | (8180.13; 10401.98) | (6555.30; 8642.69) | (6546.62; 9354.83) | (5777.18; 8877.86) |
| $N_{cit} \in I(r)$ | Yes | No | No | Yes | Yes | Yes | Yes |
| $\overline{I(r)}$ | 23620.03 | 14053.21 | 9656.61 | 9291.06 | 7599.00 | 7950.73 | 7327.52 |
| $\delta_3$ | **-0.055** | **0.194** | **0.270** | **0.118** | **-0.057** | **0.048** | **-0.036** |
| $q' := N_{cit}(h)/h^2$ | 2.782 | 2.163 | 3.16785 | 3.6059 | 4.35533 | 4.30102 | 5.26177 |
| $q'/e$ | 0.0264 | 0.03288 | 0.05040414 | 0.05358109 | 0.075203405 | 0.06345 | 0.0698127 |
| $h_q' := h+q'$ | 81.782 | 63.163 | 52.168 | 51.606 | 47.356 | 46.301 | 43.262 |
| $h_{NA}/h_q'$ | **0.957** | **0.928** | **0.903** | **0.955** | **1.024** | **1.017** | **1.089** |
| $I(q')$ | (21706.44; 24122.62) | (12775.67; 14572.12) | (8402.02; 10280.61) | (8167.04; 10121.24) | (6566.59; 8876.24) | (6734.56; 8683.24) | (5544.36; 7333.75) |
| $N_{cit} \in \overline{I(q')}$ | No | No | No | Yes | Yes | Yes | No |
| $\overline{I(q')}$ | 22914.53 | 13673.90 | 9341.315 | 9144.14 | 7721.42 | 7708.9 | 6439.06 |
| $\delta_4$ | **-0.084** | **0.162** | **0.228** | **0.101** | **-0.041** | **0.016** | **-0.153** |

**Table 1 -Continued**. Price awardees data (bases on Scopus, February 2023). The bibliometric indicators and indices.

| Researcher | F. Narin | E. Garfield | T. Braun | H. Small | L. Egghe | P. Ingwersen | H.D. White |
|---|---|---|---|---|---|---|---|
| $N_p^+$ | 64 | 106 | 216 | 57 | 211 | 88 | 37 |
| $N_{cit}$ | 7209 | 11515 | 5680 | 7693 | 5640 | 3606 | 2399 |
| $h$ | 38 | 37 | 37 | 34 | 30 | 27 | 19 |
| $h/N_p^+$ | 0.594 | 0.349 | 0.171 | 0.597 | 0.142 | 0.307 | 0.514 |
| $N_{cit}(h)$ | 6823 | 10509 | 3566 | 7471 | 3995 | 2952 | 2332 |
| $N_{cit}(h)/h$ | 179.553 | 284.027 | 153.514 | 219.735 | 133.167 | 109.333 | 122.737 |
| $N_{cit}(h)/N_{cit}$ | 0.946 | 0.913 | 0.628 | 0.971 | 0.708 | 0.819 | 0.972 |
| $h_{NA}$ | 45.887 | 57.994 | 40.731 | 47.402 | 40.587 | 32.454 | 26.471 |
| $h_{NA}/h$ | **1.208** | **1.568** | **1.101** | **1.394** | **1.353** | **1.211** | **1.393** |
| $q$ | 8.450 | 14.353 | 4.210 | 11.926 | 7.878 | 7.099 | 11.920 |
| $e$ | 73.342 | 95.603 | 46.872 | 79.467 | 55.633 | 47.149 | 44.396 |
| $h$ | 38 | 37 | 37 | 34 | 30 | 27 | 19 |
| $e/h$ | 1.930 | 2.584 | 1.267 | 2.337 | 1.854 | 1.746 | 2.337 |
| $q/e$ | 0.115 | 0.150 | 0.090 | 0.150 | 0.142 | 0.151 | 0.268 |
| $I$ | (3870.2; 6148.67) | (3385.36; 6200.06) | (3882.90; 5566.85) | (2859.02; 5234.89) | (4015.81; 2270.45) | (1800.88; 3304.06) | (669.45; 1974.80) |
| $N_{cit} \in I$ | No | No | No | No | No | No | No |
| $\overline{I}$ | 5009.47 | 4792.71 | 4724.88 | 4046.96 | 3143.13 | 2552.47 | 1322.13 |
| $\delta_1$ | **-0.305** | **-0.584** | **-0.168** | **- 0.474** | **-0.443** | **-0.292** | **-0.449** |
| $h_q := h+q$ | 46.450 | 51.353 | 41.219 | 45.926 | 37.878 | 34.099 | 30.920 |
| $h_{NA}/h_q$ | **0.988** | **1.129** | **0.988** | **1.032** | **1.072** | **0.954** | **0.856** |
| $I(q)$ | (5782.90; 9187.24) | (6521.28; 11943.29) | (4818.90; 6908.78) | (5216.46; 9551.40) | (3619.47; 6401.84) | (2872.37; 5269.91) | (1751.52; 5266.87) |
| $N_{cit} \in I(q)$ | Yes | Yes | Yes | Yes | Yes | Yes | Yes |
| $\overline{I(q)}$ | 7485.07 | 9232.29 | 5863.84 | 7383.93 | 5010.66 | 4071.14 | 3509.19 |
| $\delta_2$ | **0.038** | **-0.198** | **0.032** | **-0.040** | **-0.112** | **0.129** | **0.463** |
| $r$ | 8 | 13 | 4 | 11 | 7 | 6 | 11 |
| $e$ | 73.342 | 95.603 | 46.872 | 79.467 | 55.633 | 47.149 | 44.396 |
| $h_{NA}$ | 45.887 | 57.994 | 40.731 | 47.402 | 40.587 | 32.454 | 26.471 |
| $h+r$ | 46 | 50 | 41 | 45 | 37 | 33 | 30 |
| $h_{NA}/(h+r)$ | 0.998 | 1.160 | 0.993 | 1.053 | 1.097 | 0.983 | 0.882 |
| $r/e$ | 0.109 | 0.136 | 0.085 | 0.138 | 0.126 | 0.127 | 0.248 |
| $I(r)$ | (5750.33; 8911.26) | (6389.78; 11045.32) | (4814.88; 6779.47) | (5146.50; 8985.23) | (3581.77; 5940.77) | (2839.87; 4737.73) | (1743.58; 4797.44) |
| $N_{cit} \in I(r)$ | Yes | No | Yes | Yes | Yes | Yes | Yes |
| $\overline{I(r)}$ | 7330.80 | 8717.55 | 5797.18 | 7065.87 | 4761.27 | 3788.80 | 3270.51 |
| $\delta_3$ | **0.017** | **-0.243** | **0.021** | **-0.082** | **-0.156** | **0.051** | **0.363** |
| $q':=N_{cit}(h)/h^2$ | 4.725 | 7.676 | 2.605 | 6.463 | 4.439 | 4.050 | 6.460 |
| $e$ | 73.342 | 95.603 | 46.872 | 79.467 | 55.633 | 47.149 | 44.396 |
| $q'/e$ | 0.064 | 0.080 | 0.056 | 0.081 | 0.080 | 0.086 | 0.146 |
| $h_q' := h+q'$ | 42.725 | 44.676 | 39.605 | 40.463 | 34.439 | 31.049 | 25.460 |
| $h_{NA}/h_q'$ | **1.074** | **1.298** | **1.028** | **1.171** | **1.179** | **1.045** | **1.040** |
| $I(q')$ | (5470.41; 7080.98) | (5780.31; 7975.14) | (4788.34; 5985.49) | (4730.77; 6554.26) | (3438.51; 4734.49) | (2758.07; 3891.98) | (1620.42; 2912.08) |
| $N_{cit} \in \overline{I(q')}$ | No | No | Yes | No | No | Yes | Yes |
| $\overline{I(q')}$ | 6275.70 | 6877.73 | 5386.92 | 5642.52 | 4086.5 | 3325.03 | 2266.25 |
| $\delta_4$ | **-0.129** | **-0.403** | **-0.052** | **-0.267** | **-0.275** | **-0.078** | **-0.055** |

From Table 1 we see that the estimates $|h_{NA}/h_q - 1| < 0.1$ and $|h_{NA}/h_q' - 1| < 0.1$ are true for 11 scientists. Moreover, the values of $|\delta_1|$ are greater than 0.1 for 11 researchers (L. Leydesdorff, R. Rousseau, A. Schubert, B. Martin, F. Narin, E. Garfield, T. Braun, H. Small, L. Egghe, P. Ingwersen and H.D. White). Also, the values $|\delta_2|$ are greater than 0.1 for seven researchers (W. Glänzel, H.F. Moed, A.F.J. Van Raan, E. Garfield, L. Egghe, P. Ingwersen and H.D. White). Moreover, the values $|\delta_3|$ are greater than 0.1 for nine researchers (W. Glänzel, H.F. Moed, A.F.J. Van Raan, E. Garfield, L. Egghe, P. Ingwersen and H.D. White). The values $|\delta_4|$ are greater than 0.1 for seven researchers (W. Glänzel, H.F. Moed, A.F.J. Van Raan, F. Narin, E. Garfield, H. Small and L. Egghe). The smallest values of $|\delta_1|$, $|\delta_2|$, $|\delta_3|$ and $|\delta_4|$ of Table 1 are 0.045 (A.F.J. Van Raan), 0.001 (B. Martin), 0.017 (F. Narin) and 0.016 (A. Schubert), respectively. The ratio $h_{NA}/h$ is the biggest for E. Garfield ($h_{NA}/h = 1.568$), who published his last work in 2016, one year before his death (with the title "Citation indexes from the Science Citation Index to the Web of Science").

## 5.2. The estimates of total number of citations of a researcher using his $h$–core and a relatively small number of citations from the tail

We recall the following notions, notations and expressions presented in Tables 2 and 3.

$N_{cit}$ - the total number of citations; $N_{cit}(h) = \sum_{j=1}^{h} cit_j$ - the total number of citations into $h$-core; $h_{NA} = 0.5404446\sqrt{N_{cit}}$ - the normal approximation of $h$-index (the expression (3.3)); $d$ the $h$-defect of asymptotic normality of citations $cit_1, cit_2, ..., cit_p$ (Definition 4.4); $h_d$ -and $h_{d+1}$ - the shifted $h$-indices of order $d$ and $d+1$, respectively (Definition 4.2); $N_h^{(0)} = N_{cit}(h) = \sum_{j=1}^{h} cit_j$ ((4.3) of Definition 4.2); $N_h^{(d)} = \sum_{j=d+1}^{d+h_d} cit_j$ ((4.3) of Definition 4.2); $N_h^{(d+1)} = \sum_{j=d+2}^{d+1+h_{d+1}} cit_j$ ((4.3) of Definition 4.2); $N_{cit}^{(0)} = N_{cit} = \sum_{j=1}^{p} cit_j$ ((4.4) of Definition 4.2); $N_{cit}^{(d)} = \sum_{j=d+1}^{p} cit_j$ ((4.4) of Definition 4.2); $N_{cit}^{(d+1)} = \sum_{j=d+2}^{p} cit_j$ ((4.4) of Definition 4.2);

$h_{NA}^{(0)} = h_{NA} = \sqrt{6}\ln 2/\pi \sqrt{N_{cit}^{(k)}} = 0.54044463 94667307 \sqrt{N_{cit}}$ (the expression (3.3));

$h_{NA}^{(d)} = \sqrt{6}\ln 2/\pi \sqrt{N_{cit}^{(d)}} = 0.54044463 94667307 \sqrt{N_{cit}^{(d)}}$ (the expression (3.3));

$h_{NA}^{(d+1)} = \sqrt{6}\ln 2/\pi \sqrt{N_{cit}^{(d+1)}} = 0.54044463 94667307 \sqrt{N_{cit}^{(d+1)}}$ (the expression (3.3));

$e_d = \sqrt{N_h^{(d)} - (h_d)^2} = \sqrt{\sum_{j=d+1}^{d+h_k} cit_j - (h_d)^2}$, $e_{d+1} = \sqrt{N_h^{(d+1)} - (h_{d+1})^2} = \sqrt{\sum_{j=d+2}^{d+1+h_{d+1}} cit_j - (h_{d+1})^2}$ ((4.5)

from Definition 4.2);

$q_d = \dfrac{2N_h^{(d)}}{(h_d)^2} - 1$, $q_{d+1} = \dfrac{2N_h^{(d+1)}}{(h_{d+1})^2} - 1$ ((4.6) from Definition 4.2);

$J_d$ (the interval (4.9) from Definition 4.3)); $\overline{J_d}$ (the mean of interval $J_d$, (4.11) from Definition 4.5);

$J_{d+1}$ (the interval (4.9) from Definition 4.3)); $\overline{J_{d+1}}$ (the mean of interval $J_{d+1}$, (4.11) from Definition 4.3);

$A := (\overline{J_d} + \overline{J_{d+1}})/2$ (4.11) from Definition 4.5);

$\alpha_d, \alpha_{d+1}, \beta_d, \beta_{d+1}, B', B'', B$ (Definition 4.6);

the interval $(B', B'')$ if $B' \leq B''$ or $(B'', B')$ if $B'' \leq B'$;

the absolute errors $\Delta_A := N_{cit} - A$ and $\Delta_B = N_{cit} - B$;

the relative errors $\delta_d := (N_{cit} - \overline{J_d})/N_{cit}$, $\delta_A = (N_{cit} - A)/N_{cit}$ and $\delta_B = (N_{cit} - B)/N_{cit}$.

Using data of Table A1 and the above notions, notations and expressions, we obtain the following Table 2.

**Table 2**. Price awardees data (bases on Scopus, February 2023). The extended bibliometric indicators and indices

| Researcher | L. Leydesdorff | W. Glänzel | H.F. Moed | A.F.J. Van Raan | R. Rousseau | A. Schubert | B. Martin |
|---|---|---|---|---|---|---|---|
| $N_p^+$ | 406 | 258 | 127 | 124 | 295 | 141 | 85 |
| $N_{cit}$ | 25005 | 11766 | 7606 | 8308 | 8053 | 7587 | 7598 |
| $h$ | 79 | 61 | 49 | 48 | 43 | 42 | 38 |
| Case of Definition 4.4 | 2) | 2) | 2) | 2) | 2) | 2) | 2) |
| $d$ (the $h$ - defect of AN) | 3 | 2 | 5 | 14 | 2 | 9 | 20 |
| $h_d$ | 78 | 59 | 44 | 39 | 42 | 36 | 23 |
| $h_d + d$ | **81** | **61** | **50** | **53** | **44** | **45** | **43** |
| $(h_d + d)/h$ | **1.025** | **1.000** | **1.020** | **1.104** | **1.023** | **1.071** | **1.132** |
| $(h_d + d)/N_p^+$ | 0.200 | 0.240 | 0.394 | 0.427 | 0.149 | 0.319 | 0.506 |
| $N_{cit}(h)$ | 17360 | 8049 | 6351 | 6833 | 5203 | 6359 | 7048 |
| $N_h^{(d)}$ | 12284 | 7066 | 3998 | 3049 | 3711 | 2592 | 1093 |
| $N_{cit}^{(d)}$ | 19850 | 10783 | 5161 | 4463 | 6519 | 3707 | 1496 |
| $h_{NA}^{(d)}$ | 76.1433 | 56.12 | 38.826 | 35.485 | 43.636 | 32.905 | 20.903 |
| $h_{NA}^{(d)}/h_d$ | **0.976** | **0.920** | **0.863** | **0.910** | **1.039** | **0.914** | **0.909** |
| $e_d$ | 78.740 | 59.875 | 45.409 | 39.090 | 44.125 | 36.000 | 23.785 |
| $q_d$ | 3.038 | 3.060 | 3.130 | 3.009 | 3.160 | 3.000 | 3.132 |
| $q_d/e_d$ | 0.039 | 0.051 | 0.069 | 0.077 | 0.072 | 0.083 | 0.132 |
| $J_d$ | (24405.9, 27623.3) | (11714, 14150.1) | (8191.0, 10018.6) | (8747.6, 10351.1) | (6739.4, 10003.4) | (7608.4, 9087.5) | (7446.9, 8422.4) |
| $N_{cit} \in J_d$ | Yes | Yes | No | No | Yes | No | Yes |
| $\overline{J_d}$ | 26015.9 | 12932.1 | 9104.8 | 9549.4 | 8371.4 | 8348.0 | 7944.7 |
| $\beta_d$ | 0.740 | 0.875 | 1 | 0.090 | 1 | 0 | 0.749 |
| $\alpha_d$ | 0.260 | 0.125 | 0 | 0.910 | 0 | 1 | 0.251 |
| $\delta_d = (N_{cit} - \overline{J_d})/N_{cit}$ | **-0.040** | **-0.058** | **-0.197** | **-0.149** | **0.056** | **0.0028206** | **0.077** |
| $h_{d+1}$ | 78 | 59 | 44 | 39 | 42 | 36 | 23 |
| $N_h^{(d+1)}$ | 11882 | 6708 | 3794 | 2938 | 3275 | 2477 | 1048 |
| $N_{cit}^{(d+1)}$ | 19371 | 10366 | 4913 | 4161 | 6041 | 3556 | 1428 |
| $e_{d+1}$ | 76.151 | 56.807 | 42.942 | 37.643 | 38.8716 | 34.366 | 22.782 |
| $q_{d+1}$ | 2.906 | 2.85406 | 2.919 | 2.863 | 2.71315 | 2.822 | 2.962 |
| $q_{d+1}/e_{d+1}$ | 0.038 | 0.050 | 0.068 | 0.076 | 0.070 | 0.082 | 0.130 |
| $J_{d+1}$ | (24904.3, 28084.2) | (12150.5, 14545.6) | (8453.9; 10249.6) | (8592.4; 10176.8) | (7231.8, 8923.9) | (7769.2; 9226.9) | (7540.8; 8482.8) |
| $N_{cit} \in J_{d+1}$ | Yes | No | No | No | Yes | No | Yes |
| $\overline{J_{d+1}}$ | 26496.3 | 13348.0 | 9351.8 | 9384.6 | 8077.9 | 8498.0 | 8011.8 |
| $\beta_{d+1}$ | 0 | 0 | 0.105 | 1 | 1 | 1 | 0.782 |
| $\alpha_{d+1}$ | 1 | 1 | 0.895 | 0 | 0 | 0 | 0.218 |
| $A := (\overline{J_d} + \overline{J_{d+1}})/2$ | 26255.1 | 12900.6 | 9228.3 | 9467.0 | 8224.7 | 8423 | 7978.2 |
| $\Delta_A := N_{cit} - A$ | -1250.1 | -1134.55 | -1622.3 | -1159.0 | -171.7 | -836 | -380.2 |
| $\delta_A = (N_{cit} - A)/N_{cit}$ | **-0.050** | **-0.096** | **-0.213** | **-0.140** | **-0.021** | **-0.110** | **0.050** |
| $(B'', B')$ or $(B', B'')$ | **(24894.3, 26.787.7)** | **(12150.5, 13847.1)** | (9213.9, 10061.9) | **(8592.4, 8891.3)** | (7237.8, 10003.4) | (7608.4, 7769.2) | (7746.5, 8177.6) |
| $N_{cit} \in (B', B'')$ Or $N_{cit} \in (B'', B')$ | Yes | No | No | No | Yes | No | No |
| $B := (B + B')/2$ | 25846.0 | 12998.8 | 9637.9 | 9056.2 | 8620.6 | 7688.8 | 7964.4 |
| $\Delta_B := N_{cit} - B$ | -841 | -1232.8 | -2031.9 | -748.2 | -567.6 | -101.8 | -366.4 |
| $\delta_B = (N_{cit} - B)/N_{cit}$ | **0.034** | **-0.105** | **-0.267** | **-0.090** | **-0.070** | **-0.013** | **-0.048** |

**Table 2 -Continued**. Price awardees data (bases on Scopus, February 2023). The extended bibliometric indicators and indices.

| Researcher | F. Narin | E. Garfield | T. Braun | H. Small | L. Egghe | P. Ingwersen | H.D. White |
|---|---|---|---|---|---|---|---|
| $N_p^+$ | 64 | 106 | 216 | 57 | 211 | 88 | 37 |
| $N_{cit}$ | 7209 | 11515 | 5680 | 7693 | 5640 | 3606 | 2399 |
| $h$ | 38 | 37 | 37 | 34 | 30 | 27 | 19 |
| Case of Definition 4.4 | 2) | 2) | 2) | 2) | 2) | 2) | 1a) |
| $d$ (the $h$ - defect of AN) | 23 | 25 | 4 | 23 | 4 | 6 | 0 |
| $h_d$ | 21 | 21 | 34 | 17 | 29 | 25 | 19 |
| $h_d + d$ | **44** | **46** | **38** | **40** | **33** | **31** | **19** |
| $(h_d + d)/h$ | **1.158** | **1.243** | **1.027** | **1.176** | **1.222** | **1.148** | **1.000** |
| $(h_d + d)/N_p^+$ | 0.688 | 0.434 | 0.176 | 0.702 | 0.375 | 0.352 | 0.051 |
| $N_{cit}(h)$ | 6823 | 10509 | 3566 | 7471 | 3995 | 2952 | 2332 |
| $N_h^{(d)}$ | 921 | 901 | 2378 | 599 | 1702 | 1259 | 2332 |
| $N_{cit}^{(d)}$ | 1132 | 1682 | 4456 | 705 | 3259 | 1798 | 2399 |
| $h_{NA}^{(d)}$ | 18.183 | 22.165 | 36.077 | 14.350 | 30.853 | 22.916 | 26.471 |
| $h_{NA}^{(d)}/h_d$ | **0.866** | **1.055** | **1.061** | **0.844** | **1.064** | **0.888** | **1.393** |
| $e_d$ | 21.909 | 21.448 | 34.957 | 17.607 | 29.343 | 25.179 | 44.396 |
| $q_d$ | 3.177 | 3.086 | 3.114 | 3.145 | 3.048 | 3.029 | 11.920 |
| $q_d/e_d$ | 0.145 | 0.144 | 0.089 | 0.179 | 0.104 | 0.120 | 0.268 |
| $J_d$ | (7180.7, 8056.5) | (10939.5, 11808.6) | (4508.1, 5918.4) | (7655.5, 8362.6) | (4693.3, 5889.5) | (3454.0, 4483.6) | (661.4, **1988.8**) |
| $N_{cit} \in J_d$ | Yes | Yes | Yes | Yes | Yes | Yes | No |
| $\overline{J_d}$ | 7618.6 | 11374.1 | 5213.3 | 8009.1 | 5291.4 | 3968.8 | 1325.1 |
| $\beta_d$ | 0.909 | 0.448 | 0.957 | 0.607 | 0.657 | 0.179 | - |
| $\alpha_d$ | 0.091 | 0.552 | 0.043 | 0.393 | 0.343 | 0.821 | - |
| $\delta_d = (N_{cit} - \overline{J_d})/N_{cit}$ | -0.057 | 0.012 | 0.082 | -0.041 | 0.062 | -0.101 | 0.448 |
| $h_{d+1}$ | 21 | 21 | 34 | 17 | 28 | 24 | 19 |
| $N_h^{(d+1)}$ | 874 | 849 | 2261 | 559 | 1549 | 1104 | 1322 |
| $N_{cit}^{(d+1)}$ | 1064 | 1609 | 4305 | 648 | 3106 | 1653 | 1303 |
| $e_{d+1}$ | 20.809 | 20.199 | 33.242 | 16.432 | 27.659 | 22.978 | 31.000 |
| $q_{d+1}$ | 2.964 | 2.850 | 2.912 | 2.869 | 5.503 | 2.833 | 6.324 |
| $q_{d+1}/e_{d+1}$ | 0.142 | 0.141 | 0.087 | 0.175 | 0.199 | 0.123 | 0.204 |
| $J_{d+1}$ | (7255.4; 8115.6) | (11019.8; 11872.0) | (4659.1; 6069.4) | (7719.1; 8410.1) | (4256.4; 6392.5) | (3468.7; 4441.4) | (1818.8; **2810.7**) |
| $N_{cit} \in J_{d+1}$ | No | Yes | Yes | No | Yes | Yes | Yes |
| $\overline{J_{d+1}}$ | 7685.5 | 11445.9 | 5364.2 | 8064.6 | 5324.4 | 3955.1 | 2314.8 |
| $\beta_{d+1}$ | 0.809 | 0.199 | 0.242 | 0.432 | 0.659 | 1 | - |
| $\alpha_{d+1}$ | 0.191 | 0.801 | 0.758 | 0.568 | 0.341 | 0 | - |
| $A := (\overline{J_d} + \overline{J_{d+1}})/2$ | 7652.1 | 11410.0 | 5288.8 | 8036.9 | 5307.9 | 3962.0 | 1820.0 |
| $\Delta_A := N_{cit} - A$ | -443.05 | 105.0 | 391.3 | -343.3 | 332.08 | -356 | -579.0 |
| $\delta_A = (N_{cit} - A)/N_{cit}$ | **-0.061** | **0.009** | **0.069** | **0.045** | **0.059** | **-0.099** | **-0.241** |
| $(B'',B')$ or $(B',B'')$ | **(7420.0, 7976.7)** | **(11328.5, 11700.5)** | **(5728.1, 5857.9)** | (8046.6, 8111.8) | **(5479.2, 5662.3)** | (3638.7, 3468.7) | - |
| $N_{cit} \in (B',B'')$ or $N_{cit} \in (B'',B')$ | **No** | **Yes** | **No** | **No** | **Yes** | **Yes** | - |
| $B := (B+B')/2$ | 7698.4 | 11515.45 | 5793.0 | 8098.2 | 5570.8 | 3553.7 | 2399.75 |
| $\Delta_B := N_{cit} - B$ | -489.4 | **-0.45** | -113 | -405.2 | 69.2 | 52.3 | **-0.75** |
| $\delta_B = (N_{cit} - B)/N_{cit}$ | **-0.068** | **-0.000004** | **-0.020** | **-0.053** | **0.012** | **-0.015** | **-0.000313** |

From Table 2 we see that the value $(h_d + d)/N_p^+$ is the smallest for H.D. Whte (0.051) and it is largest for H. Small (0.737). From Table 2 we also see that the case 2) of Definition 4.4 refers to 13 Price awardees, except H.D White for which refers the case 1a) of Definition 4.4. Moreover, the $h$ - defect of asymptotic normality of citations $cit_1, cit_2,..., cit_p$, denoted as $d$, is greatest for for E. Garfield ($d = 25$), and then for H. Small ($d = 23$), F. Narin ($d = 23$), B. Martin ($d = 20$) and A.F.J. Van Raan ($d = 14$).

The case 2) of Definition 4.4 refers to 13 researchers except H.D. White for which $e_k > h_k$ for all $k = 0, 1,..., h - 1 = 18$ (the case 1a) of Definition 4.4). Also, the equality $h_d + d = h$ holds for W. Glänzel and H.D. White. The estimate $|h_{NA}^{(d)}/h - h_d| < 0.1$ is satisfied for 9 researchers. The relative errors $|\delta_d| = |(N_{cit} - \overline{J_d})/N_{cit}|$ are less than 0.1 for 10 researchers, while for A. Schubert this value is the smallest and equal to **0.0028206**. Also, the relative errors $|\delta_A| = |(N_{cit} - A)/N_{cit}|$ are less than 0.1 for 10 researchers, while for E. Garfield this value is the smallest and equal to **0.009**. The relative errors $|\delta_B| = |(N_{cit} - B)/N_{cit}|$ are les than 0.1 for 12 researchers (except W. Glänzel and H.F. Moed), while for H.D. White this value is the smallest and equal to 0.002. **The smallest value of $|\delta_B|$ has E. Garfield, which is equal to the incredibly value 0.000004 close to 0 (with the surprisingly absolute value** $\Delta_B := N_{cit} - B = -0.45$), and after him follows H.D. White whose $|\delta_B|$ value is **0.000313 (with the surprisingly small absolute value** $\Delta_B := N_{cit} - B = -0.75$).

Generally, it seems that among all the estimates of $N_{cit}$ from Tables 1 and 2, the estimate $B := (B' + B'')$ from Table 2 could be considered as the best approximation of $N_{cit}$. Namely, from Table 2 we see that the inequality $|\delta_B| < |\delta_A|$ is satisfied for 11 researchers.

Supposing that each partition of $N_{cit}$ is chosen with equal probability, Yong (2014) presented the confidence intervals for $h$-index, as a function of $N_{cit}$ whose value vary between 50 and 10000 (Table 1 in Yong, 2014). Notice that theory concerning the asymptotics of these uniform random partitions was largely developed by Vershik et al. (see Su, 2008, pp. 44–79).

Using data of Table 1 in Yong (2014) and the above Table 2, we obtain the following two tables. Recall that the first two rows of Tables 3 and 4 are taken from Table 1 in Yong (2014).

From Table 3 we see that $h \in C$ is satisfied for six Price awardees. The estimates $1 \leq h/\overline{C} < 1.15$ is true for three Price awardees (with the highest $h$-indices), while the estimates $0 < 1 - h_d/\overline{C}_d < 0.11$ is true for three Price awardees.

**Table 3**. The confidence intervals for the $h$-indices (Table 1 in Yong, 2014) for 14 Price awardees (bases on Scopus, February 21 2023).

| $N_{cit}$ (Table 1 in Yong, 2014) | 10000 | 7500 | 8000 | 8000 | 7500 | 7500 | 7000 | 10000 | 5500 | 7500 | 5500 | 3500 | 2500 |
|---|---|---|---|---|---|---|---|---|---|---|---|---|---|
| The confidence interval for $h$-index, $C$ | [47, 60] | [40, 52] | [42, 54] | [42, 54] | [40, 52] | [40, 52] | [39, 51] | [47, 60] | [35, 45] | [40, 52] | [35, 45] | [27, 36] | [22, 31] |
| The mean of $C$, $\overline{C}$ | 53.5 | 46 | 48 | 48 | 46 | 46 | 45 | 53.5 | 40 | 46 | 40 | 31.5 | 26.5 |
| Researcher | W. Glänzel | H.F. Moed | A.F.J. Van Raan | R. Rousseau | A. Schubert | B. Martin | F. Narin | E. Garfield | T. Braun | H. Small | L. Egghe | P. Ingwersen | H.D. White |
| $h$-index | 61 | 49 | 48 | 43 | 42 | 38 | 38 | 37 | 37 | 34 | 30 | 27 | 19 |
| $N_{cit}$ of researcher | 11766 | 7606 | 8308 | 8053 | 7587 | 7598 | 7209 | 11515 | 5680 | 7693 | 5640 | 3606 | 2399 |
| $h \in C$ | No | Yes | Yes | Yes | Yes | No | No | No | Yes | No | No | Yes | No |
| $h/\overline{C}$ | **1.140** | **1.021** | **1.000** | **0.896** | **0.913** | **0.826** | **0.844** | **0.692** | **0.925** | **0.740** | **0.750** | **0.856** | **0.717** |

**Table 4**. The estimations of $h_d$-indices, based on related confidence intervals (Table 1 in Yong, 2014) for 14 Price awardees (bases on Scopus, February 2023).

| $N_{cit}$ (Table 1 in Yong, 2014) | 10000 | 5500 | 4500 | 6500 | 3500 | 1500 | 1250 | 1500 | 4500 | 500 | 3500 | 2000 | 2500 |
|---|---|---|---|---|---|---|---|---|---|---|---|---|---|
| The confidence interval for $h_d$-index, $C_d$ | [47, 60] | [35, 45] | [31, 41] | [37, 49] | [27, 36] | [17, 24] | [15, 22] | [17, 24] | [31, 41] | [9, 14] | [27, 36] | [20, 28] | [22, 31] |
| The mean of $C_d$, $\overline{C}_d$ | 53.5 | 40 | 36 | 43 | 31.5 | 20.5 | 18.5 | 20.5 | 36 | 11.5 | 31.5 | 24 | 26.5 |
| Researcher | W. Glänzel | H.F. Moed | A.F.J. Van Raan | R. Rousseau | A. Schubert | B. Martin | F. Narin | E. Garfield | T. Braun | H. Small | L. Egghe | P. Ingwersen | H.D. White |
| $h_d$-index | 59 | 45 | 39 | 42 | 36 | 23 | 21 | 21 | 34 | 14 | 29 | 25 | 19 |
| $N_{cit}^{(d)}$ of researcher | 10783 | 5739 | 4463 | 6519 | 3707 | 1496 | 1132 | 1682 | 4456 | 457 | 3259 | 1798 | 2399 |
| $h_d \in C_d$ | Yes | Yes | Yes | Yes | Yes | Yes | Yes | Yes | Yes | Yes | Yes | Yes | No |
| $h_d+d \in C_d$ | No | No | No | Yes | No | No | No | No | Yes | No | Yes | No | No |
| $h_d/\overline{C}_d$ | **1.103** | **1.125** | **1.083** | **0.977** | **1.143** | **1.122.** | **1.135** | **1.024** | **0.944** | **1.217** | **0.921** | **1.042** | **0.717** |
| $(h_d+d)/\overline{C}_d$ | **1.140** | **1.250** | **1.444** | **1.023** | **1.429** | **2.098** | **2.378** | **2.244** | **1.056** | **3.652** | **1.048** | **1.292** | **0.717** |

From Table 4 we see that $h_d \in C_d$ is satisfied for 12 Price awardees, while this is not true for H.D. White. The estimate $1 < h_d/\overline{C}_d < 1.15$ is true for eight Price awardees, while the estimates $0 < 1 - h_d/\overline{C}_d < 0.08$ is true for three Price awardees.

Using data of Tables A3, A4, A5 and A6 and the above notions, notations and expressions, we obtain the following Table 5.

**Table 5.** Researcher's data (bases on Google Scholar and Scopus, April 2023-June 2023) - Bibliometric indicators and indices.

| Researcher | S Freud (GS) | H.J. Kim (GS) | R.C. Kessler (GS) | A. Einstein (GS) | P. Erdös (GS) | T. Tao (GS) | L. Leydesdorf (GS) | C.G. Meyer (GS) | T. Tao (S) |
|---|---|---|---|---|---|---|---|---|---|
| $p$ | 2991 | 3000 | 1900 | 1044 | 1608 | 584 | 912 | 1115 | 297 |
| $N_p^+$ | 2378 | 2613 | 1559 | 617 | 1263 | 448 | 635 | 701 | 278 |
| $N_{cit}$ | 643730 | 518589 | 515591 | 161009 | 99866 | 90963 | 70821 | 49110 | 46852 |
| $h$ | 288 | 338 | 329 | 123 | 128 | 106 | 117 | 95 | 70 |
| $N_p^+/h$ | 8.257 | 7.731 | 4.739 | 5.016 | 9.867 | 4.226 | 5.427 | 7.379 | 3.971 |
| $N_{cit}(h)$ | 559724 | 312483 | 432109 | 150565 | 72123 | 73134 | 54898 | 36684 | 42192 |
| $N_{cit}(h)/h$ | 1943.5 | 924.5 | 1313.4 | 1224.1 | 563.5 | 689.9 | 469.2 | 386.1 | 602.7 |
| $N_{cit}(h)/N_{cit}$ | **0.869** | 0.603 | 0.838 | 0.935 | 0.722 | 0.804 | 0.775 | 0.747 | 0.900 |
| $h_{NA}$ | 433.614 | 389.165 | 388.064 | 216.858 | 170.447 | 162.998 | 143.824 | 119.767 | 116.981 |
| $h_{NA}/h$ | **1.5056** | 1.151 | 1.180 | 1.763 | 1.332 | 1.538 | 1.229 | 1.288 | 1.671 |
| $q$ | 12.496 | 4.470 | 7.082 | 18.904 | 7.804 | 12.018 | 7.021 | 7.129 | 16.221 |
| $e$ | 690.493 | 445.24 | 570.246 | 368.016 | 236.091 | 248.793 | 203.000 | 166.310 | 193.111 |
| $e/h$ | 2.398 | 1.317 | 1.733 | **2.992** | 1.844 | 2.347 | 1.735 | 1.751 | **2.759** |
| $h_q := h+q$ | 300.4964 | 342.470 | 336.082 | 141.904 | 135.804 | 118.018 | 124.021 | 102.129 | 86.221 |
| $h_{NA}/h_q$ | 1.44299 | 1.136 | 1.1547 | 1.528 | 1.255 | 1.486 | 1.160 | 1.173 | 1.357 |
| Case of Definition 4.4 | 2) | 2) | 2) | 2) | 2) | 2) | 2) | 2) | 2) |
| $d$ (the $h$-defect of AN) | **97** | **13** | **94** | **59** | **21** | **14** | **14** | **9** | **7** |
| $h_d$ | 236 | 332 | 270 | 87 | 120 | 101 | 111 | 91 | 66 |
| $e_d$ | 236.692 | 334.798 | 271.134 | 88.272 | 120.959 | 102.587 | 111.786 | 92.423 | 66.227 |
| $h_{d+1}$ | 234.902 | 332 | 270 | 87 | 120 | 101 | 110 | 91 | 55.964 |
| $e_{d+1}$ | 236 | 331.053 | 269.596 | 86.238 | 119.126 | 100.479 | 109.827 | 90.167 | 65 |
| $h_d + d$ | 333 | 345 | 364 | 146 | 141 | 115 | 125 | 100 | 73 |
| $(h_d+d)/h$ | **1.156** | **1.021** | **1.106** | **1.187** | **1.102** | **1.085** | **1.068** | **1.053** | **1.043** |
| $(h_d+d)/N_p^+$ | 0.140 | 0.132 | 0.233 | 0.241 | 0.112 | 0.257 | 0.197 | 0.143 | 0.263 |
| $(h_d+d)/p$ | 0.011 | 0.115 | 0.192 | 0.140 | 0.088 | 0.197 | 0.137 | 0.090 | 0.246 |
| $N_{cit}^{(d)}$ | **184240** | 426070 | 218790 | 23255 | 55147 | **37093** | 40049 | 28788 | **13201** |
| $h_{NA}^{(d)}$ | 231.976 | 352.770 | 252.793 | 82.420 | 126.268 | 104.087 | 108.155 | 77.0432 | 62.095 |
| $h_{NA}^{(d)}/h_d$ | **0.983** | **1.063** | **0.936** | **0.947** | **1.052** | **1.031** | **0.974** | **0.847** | **0.887** |
| $J_d$ | (645355.4, 655060.8)) | (463086, 476765) | (540866.3, 551975.2)) | (161903.2, 165495.2) | (90938.9, 97149.9) | (86740.7, **90912.3**) | (67196.04, 791137.7) | (47024.7, 50372.2()) | (**46762**, 50377)) |
| $N_{cit} \in J_d$ | No | No | No | No | No | Yes | Yes | Yes | Yes |
| $A := (\overline{J_d} + \overline{J_{d+1}})/2$ | 650208.1 | 471172.0 | 546977.7 | 163876.5 | 94331.8 | 89604.7 | 73025.81 | 48907.5 | 48596.3 |
| $\Delta_A := N_{cit} - A$ | -6478.09 | 47416.5 | -31386.7 | -2867.5 | 5534.2 | 1358.3 | -2204.8 | 202.5 | -1744.3 |
| $\delta_A = (N_{cit}-A)/N_{cit}$ | **-0.011** | **0.091** | **-0.061** | **-0.018** | **0.055** | **0.015** | **-0.031** | **0.004** | **-0.037** |
| $(B'',B')$ or $(B',B'')$ | (**646435.0, 652069.0**) | (476765.0, 478513.0) | (**542608.9, 546470.3**) | (162881.3, 164987.5) | (95889.1, 96111.2) | (**88384.3, 90912.3**) | (69976.1, 71407.7) | (**50372.2, 50359.4**) | (**44006.8, 47854.2**) |
| $N_{cit} \in (B',B'')$ | No | No | No | No | No | No | Yes | No | Yes |
| $B := (B'+B'')/2$ | 649252.0 | 477639.0 | 544539.6 | 163934.4 | 96000.2 | 89648.3 | 70691.9 | 50365.8 | 45930.5 |
| $\Delta_B := N_{cit} - B$ | -5522 | 40950 | -28948.6 | -2925.4 | 3465.8 | 1314.7 | -411.5 | -1225.8 | 921.5 |
| $\delta_B = (N_{cit}-B)/N_{cit}$ | **-0.000858** | **0.079** | **-0.056** | **-0.018** | **0.035** | **0.022** | **-0.006** | **-0.026** | **0.020** |

**Table 5 -Continued**. Researcher's data (bases on Google Scholar and Scopus, April 2023-June 2023) - Bibliometric indicators and indices.

| Researcher | G. Andrews (GS) | M. McAleer (GS) | J. Hirsch (S) | P. Erdös (S) | F.T. Edelmann (GS) | C.F. Gauss (GS) | G. Andrews (S) | S. Papadimitrou (GS) | D. Zeilberger (S) |
|---|---|---|---|---|---|---|---|---|---|
| $p$ | 730 | 1251 | 289 | 672 | 606 | 828 | 259 | 92 | 172 |
| $N_p^+$ | 434 | 761 | 256 | 613 | 431 | 224 | 229 | 64 | 146 |
| $N_{cit}$ | 32386 | 26266 | 24451 | 18830 | 13417 | 11606 | 6567 | 5407 | 3158 |
| $h$ | 66 | 79 | 60 | 62 | 54 | 44 | 41 | 26 | 25 |
| $N_p^+/h$ | 6.576 | 9.633 | 4.267 | 9.887 | 7.981 | 5.091 | 5.585 | 2.462 | 5.84 |
| $N_{cit}(h)$ | 25557 | 15879 | 20447 | 11257 | 7338 | 9930 | 3918 | 5060 | 2222 |
| $N_{cit}(h)/h$ | 387.227 | 201.0 | 340.8 | 1181.6 | 135,9 | 225.7 | 95.6 | 194.6 | 88.9 |
| $N_{cit}(h)/N_{cit}$ | 0.789 | 0.605 | 0.836 | 0.598 | 0.547 | 0.856 | 0.597 | 0.936 | 0.704 |
| $h_{NA}$ | 97.259 | 87.589 | 84.508 | **74.161** | 62.601 | 58.223 | 43.796 | 39.740 | 30.371 |
| $h_{NA}/h$ | **1.474** | **1.109** | **1.408** | **1.196** | **1.159** | **1.323** | **1.068** | **1.528** | **1.215** |
| $q$ | 10.734 | 4.089 | 10.393 | 4.857 | 4.033 | 9.258 | 3.662 | 13.970 | 6.111 |
| $e$ | 145.606 | 98.173 | 130.027 | 86.099 | 66.498 | 89.409 | 47.297 | 66.212 | 39.963 |
| $e/h$ | 2.206 | 1.243 | 2.167 | 1.389 | 1.231 | 2.032 | 1.154 | 2.547 | 1.599 |
| $h_q := h+q$ | 76.734 | 83.089 | 70.393 | 66.857 | 58.033 | 53.258 | 44.662 | 39.970 | 31.111 |
| $h_{NA}/h_q$ | 1.267 | 1.054 | 1.201 | **1.109** | 1.079 | 1.093 | 0.981 | 0.994 | 0.976 |
| Case of Definition 4.4 | 2) | 2) | 2) | 2) | 2) | 2) | 2) | 2) | 2) |
| $d$ (the $h$-defect of AN) | **7** | **10** | **12** | **11** | **3** | **18** | **1** | **13** | **10** |
| $h_d$ | 63 | 74 | 55 | 58 | 54 | 34 | 40 | 19 | 20 |
| $e_d$ | 64.861 | 74.135 | 55.723 | 58.481 | 54.580 | 35.071 | 40.633 | 19.596 | 20.976 |
| $h_{d+1}$ | 62 | 73 | 55 | 58 | 53 | 34 | 39 | 19 | 20 |
| $e_{d+1}$ | 61.213 | 72.719 | 53.824 | 56.577 | 51.137 | 33.734 | 38.000 | 17.378 | 18.841 |
| $h_d + d$ | 70 | 84 | 67 | 69 | 57 | 52 | 41 | 32 | 30 |
| $(h_d + d)/h$ | 1.061 | 1.063 | 1.117 | 1.113 | 1.056 | 1.182 | **1.000** | 1.231 | 1.200 |
| $(h_d + d)/N_p^+$ | 1.447 | 0.110 | 0.262 | 0.113 | 0.132 | 0.232 | 0.179 | 0.500 | 0.205 |
| $(h_d + d)/p$ | 0.096 | 0.067 | 0.239 | 0.103 | 0.094 | 0.063 | 0.158 | 0.348 | 0.174 |
| $N_{cit}^{(d)}$ | 14750 | 21132 | 9665 | 13949 | 11812 | 3765 | 5900 | 964 | 3158 |
| $h_{NA}^{(d)}$ | 65.637 | 78.5636 | 53.1315 | 63.830 | 58.737 | 33.162 | 21.866 | 16.780 | 21.866 |
| $h_{NA}^{(d)}/h_d$ | **1.042** | **1.062** | 0.966 | **1.101** | **1.088** | 0.975 | **1.093** | 0.883 | **1.093** |
| $J_d$ | (29948.9, **32563.4**) | (21553.3, **26388.4**) | (24039,0, 26038.1) | (15243.6, 17633.1) | (10506.3, 12732.9) | (11124.3, 12536.3) | (5350.0, 7002.0) | (5315.9, 6105.0) | (2493.3, 3355.5) |
| $N_{cit} \in J_d$ | Yes | Yes | Yes | No | No | Yes | Yes | Yes | Yes |
| $A := (\overline{J_d} + \overline{J_{d+1}})/2$ | 31334.2 | 23747.8 | 25173.5 | 16566.7 | 11670.1 | 11826.5 | 6184.39 | 5752.3 | 2964.1 |
| $\Delta_A := N_{cit} - A$ | 1051.8 | 2518.2 | -722.5 | 2263.3 | 1746.9 | 220.5 | 382.61 | -345.3 | 193.9 |
| $\delta_A = (N_{cit} - A)/N_{cit}$ | **0.032** | **0.096** | **-0.030** | **0.120** | 0.1302 | **0.019** | **0.058** | **-0.064** | **0.061** |
| $(B'', B')$ or $(B', B'')$ | **(32141.0, 32563.4)** | (24288.6, 24513.1) | **(24051.0, 25678.9)** | (16392.5, 16511.3) | **(10697.3, 11798.2)** | **(12536.3, 12281.9)** | (6394.9, 6985.3) | **(5529.7, 5785.9)** | **(2583.0, 2834.8)** |
| $N_{cit} \in (B', B'')$ or $N_{cit} \in (B'', B')$ | **Yes** | **No** | **Yes** | **No** | **No** | **No** | **Yes** | **No** | **No** |
| $B := (B' + B'')/2$ | 32367.2 | 24400.9 | 24865.0 | 16451.9 | 11247.8 | 12409.1 | 6690.1 | 5657.8 | 2708.9 |
| $\Delta_B := N_{cit} - B$ | 18.8 | 1865.1 | -414.0 | 2378.1 | 2169.2 | -803.1 | -123.1 | -250.8 | 449.1 |
| $\delta_B = (N_{cit} - B)/N_{cit}$ | **0.000580** | **0.071** | **-0.017** | **0.126** | **0.162** | **0.069** | **-0.019** | **-0.046** | **0.142** |

**Table 5 -Continued**. Researcher's data (bases on Google Scholar and Scopus, April 2023-June 2023) - Bibliometric indicators and indices.

| Researcher | I. Orović (GS) | C.D. Savage (GS) | V. Spalević (GS) | P. Ziarati (GS) | A. Yong (GS) | D. Kalaj (GS) | P. Vukoslavčević (GS) | K. Monkova (GS) | Lj.R. Mutafchiev (GS) |
|---|---|---|---|---|---|---|---|---|---|
| $p$ | 177 | 104 | 412 | 180 | 100 | 162 | 63 | 114 | 59 |
| $N_p^+$ | 148 | 96 | 233 | 140 | 72 | 122 | 37 | 77 | 38 |
| $N_{cit}$ | 3082 | 3031 | 2832 | 2123 | 1631 | 1577 | 912 | 741 | 312 |
| $h$ | 31 | 29 | 29 | 26 | 24 | 22 | 14 | 16 | 10 |
| $N_p^+/h$ | 4.774 | 3.310 | 8.035 | 6.923 | 3.000 | 5.546 | 2.643 | 4.813 | 3.800 |
| $N_{cit}(h)$ | 1993 | 2296 | 1199 | 1034 | 1206 | 947 | 814 | 440 | 172 |
| $N_{cit}(h)/h$ | 64.3 | 79.2 | 41.3 | 39.8 | 50.3 | 43.046 | 58.1 | 27.5 | 17.2 |
| $N_{cit}(h)/N_{cit}$ | 0.647 | 0.758 | 0.423 | 0.487 | 0.739 | 0.601 | 0.893 | 0.594 | 0.551 |
| $h_{NA}$ | 30.003 | 29.754 | 28.761 | 24.902 | 21.826 | 21.462 | 16.321 | 14.712 | 9.546 |
| $h_{NA}/h$ | **0.968** | **1.026** | **0.991** | **0.958** | **0.909** | **0.976** | **1.166** | **0.920** | **0.955** |
| $q$ | 3.148 | 4.461 | 1.851 | 2.059 | 3.188 | 2.913 | 7.306 | 2.438 | 2.440 |
| $e$ | 32.143 | 38.1445 | 18.921 | 18.921 | 25.100 | 21.517 | 24.860 | 13.565 | 8.485 |
| $e/h$ | 1.037 | 1.315 | 0.652 | 0.728 | 1.046 | 0.978 | 1.776 | 0.848 | 0.849 |
| $h_q := h+q$ | 34.148 | 33.460 | 30.851 | 28.059 | 27.188 | 24.913 | 21.306 | 18.438 | 12.440 |
| $h_{NA}/h_q$ | 0.879 | 0.889 | 0.932 | 0.668 | 0.803 | 0.861 | 0.766 | 0.798 | 0.767 |
| Case of Definition 4.4 | 2) | 2) | 3b) | 3b) | 2) | 3a) | 2) | 3b) | 3b) |
| $d$ (the $h$ - defect of AN) | **0** | **2** | **0** | **0** | **1** | **0** | **6** | **0** | **0** |
| $h_d$ | 31 | 28 | 29 | 26 | 23 | 22 | 9 | 16 | 10 |
| $e_d$ | 32.143 | 28.089 | 18.921 | 18.921 | 23.791 | 21.517 | 10.198 | 13.565 | 8.485 |
| $h_{d+1}$ | 31 | 28 | 28 | 26 | 23 | 21 | 9 | 16 | 10 |
| $e_{d+1}$ | 29.189 | 26.268 | 17.635 | 17.720 | 22.428 | 20.273 | 8.888. | 11.747 | 7.289 |
| $h_d + d$ | 31 | 30 | 29 | 26 | 24 | 22 | 15 | 16 | 10 |
| $(h_d+d)/h$ | **1.000** | **1.034** | **1.000** | **1.000** | **1.000** | **1.000** | **1.071** | **1.000** | **1.000** |
| $(h_d+d)/N_p^+$ | 0.209 | 0.313 | 0.124 | 0.186 | 0.333 | 0.180 | 0.405 | 0.208 | 0.263 |
| $(h_d+d)/p$ | 0.175 | 0.288 | 0.070 | 0.144 | 0.240 | 0.136 | 0.143 | 0.140 | 0.169 |
| $N_{cit}^{(d)}$ | 3082 | 25800 | 2832 | 2123 | 1520 | 1577 | 272 | 741 | 312 |
| $h_{NA}^{(d)}$ | 30.003 | 25.800 | 28.761 | 24.902 | 21.070 | 21.462 | 8.913 | 14.712 | 9.546 |
| $h_{NA}^{(d)}/h_d$ | **0.968** | **0.921** | **0.992** | **0.958** | **0.916** | **0.976** | **0.990** | **0.920** | **0.946** |
| $J_d$ | (2677, 3966.56) | (2891.3, 4042.9) | (2343.4, 3470.0) | (1838.1, 2341.8)) | (1475.6, 2431.8) | (1238.8, 2136.2) | (757.2, 1149.3)) | (589.8, 1219.8)) | (173.8, 567.6) |
| $N_{cit} \in J_d$ | Yes | Yes | Yes | Yes | Yes | Yes | Yes | Yes | Yes |
| $A := (\overline{J_d} + \overline{J_{d+1}})/2$ | 3223.8 | 3530.0 | 2927.1 | 2341.8 | 1995.2 | 1687.5 | 969.3 | 904.8 | 370.7 |
| $\Delta_A := N_{cit} - A$ | -141.8 | -498.9 | -95.1 | -218.8 | -364.2 | -110.5 | -57.3 | -163.8 | -58.7 |
| $\delta_A = (N_{cit}-A)/N_{cit}$ | **-0.046** | **-0.165** | **-0.034** | **-0.103** | **-0.223** | **-0.070** | **-0.063** | **-0.221** | **-0.189** |
| $(B'',B')$ or $(B',B'')$ | (2837.9, 2905.7) | (2993.9, 3082.7) | - | - | **(2102.5, 2120.7)** | - | **(834.61, 834.85)** | - | - |
| $N_{cit} \in (B',B'')$ | No | Yes | - | - | **No** | - | No | - | - |
| $B := (B'+B'')/2$ | 2871.8 | 3011.3 | 3195.9 | 2156.8 | 2111.6 | 1671.8 | 834.73 | 689.2 | 275.4 |
| $\Delta_B := N_{cit} - B$ | 210.2 | **19.7** | -363.9 | -33.8 | -480.6 | -94.8 | 77.27 | 51.80 | 36.6 |
| $\delta_B = (N_{cit}-B)/N_{cit}$ | **0.068** | **0.0065** | **-0.128** | **-0.016** | **-0.295** | **-0.060** | **-0.085** | **0.070** | **0.117** |

From Table 5 we see that the $h$-defect of asymptotic normality of citations $cit_1, cit_2, ..., cit_p$ ($d$) is the greatest for S. Freud ($d = 97$), whose Google Scholar's $h$-index is 288, and after him follows R.C. Kessler ($d = 94$), whose Google Scholar's $h$-index is 329, and then A. Einstein ($d = 59$), whose Google Scholar's $h$-index is 123. The case 2) of Definition 4.4 refers to 22 researchers, the case 3a) of Definition 4.4 refers to two one researcher (D. Kalaj, whose GS $h$-index is equal to 22) and the case 3b) of Definition 4.4 refers to four researchers (whose $h$-indices are 29, 26, 16 and 10). The equality $h_d + d = h$ holds for eight researchers (G. Andrews (S), I. Orović (GS), V. Spalević (GS), P. Ziarati (GS), A. Yong (GS), D. Kalaj (GS), K. Monkova (S) and Lj.R. Mutafchiev (GS)). Moreover, the ratio $(h_d + d)/p$ varies between 0.011 (for S. Freud) and 0.348 (for S. Papadimitrou). The estimates $|h_{NA}^{(d)}/h - h_d| < 0.1$ is satisfied for 22 researchers. The relative errors $|\delta_A| = |(N_{cit} - A)/N_{cit}|$ are less than 0.1 for 20 researchers, while this value is the smallest for C.G. Meyer from his profile at Google Scholar, and it is equal to **0.004**. Notice also that the relative error $\delta(\overline{J_d}) = \delta(\overline{J_8})$ from Google Scholar profile of C.G. Meyer is quite close to zero, namely, this value is equal to **0.000130523** (with the absolute error equals to $\Delta(\overline{J_d}) := \Delta(\overline{J_8}) = N_{cit} - B = 6.41$). The relative errors $|\delta_B| = |(N_{cit} - B)/N_{cit}|$ are les than 0.1 for 22 researchers, except five reseachers whose $h$-indices are 62, 54, 29, 24 and 10. **The smallest value of $|\delta_B|$ has G. Andrews from his profile at Google Scholar database, which is equal to the incredibly value 0.000580 close to 0 (with the surprisingly small absolute value $\Delta_B := N_{cit} - B = 18.8$), and after him follows L. Leydesdorf, whose $|\delta_B|$ value in his profile at Google Scholar database is 0.0058 (with the absolute error value $\Delta_B := N_{cit} - B = -411.5$), and then C.D. Savage, whose $|\delta_B|$ value in her profile at Google Scholar database is 0.0065 (with the absolute error value $\Delta_B := N_{cit} - B = 19.7$).** Moreover, the inequality $|\delta_B| < |\delta_A|$ is satisfied for 18 researchers. **Observe also that the length of the confidence interval $(B'', B')$ of P. Vukoslavčević from his profile at Google Scholar database is equal $B'' - B' = 0.24$.**

Finally, by data from from Tables 1, 2 and 5, we can conclude that among all the estimates of $N_{cit}$, the estimate $B := (B' + B'')$ could be considered as the best approximation of $N_{cit}$.

Using data of Table 1 in Yong (2014) and the above Table 5, we obtain the following two tables. Recall that the first two rows of Tables 6 and 7 are taken from Table 1 in Yong (2014).

**Table 6.** The confidence intervals for the $h$-indices (Table 1 in Yong, 2014) for for 13 researchers (based on their Google Scholar and Scopus profiles, April 2023-June 2023).

| $N_{cit}$ | 10000 | 6500 | 5500 | 3000 | 3000 | 3000 | 3000 | 2000 | 1500 | 1500 | 1000 | 750 | 300 |
|---|---|---|---|---|---|---|---|---|---|---|---|---|---|
| The confidence interval for $h$-index, $C$ | [47, 60] | [37, 49] | [35, 45] | [25, 34] | [25, 34] | [25, 34] | [25, 34] | [20, 28] | [17, 24] | [17, 24] | [13, 20] | [11, 17] | [7, 11] |
| The mean of $C$, $\overline{C}$ | 53.5 | 43 | 40 | 29.5 | 29.5 | 29.5 | 29.5 | 24 | 20.5 | 20.5 | 16.5 | 14 | 9 |
| Researcher | C.F. Gauss (GS) | G. Andrews (S) | S. Papadimitrou (GS) | D. Zeilberger (S) | I. Orović (GS) | C.D. Savage (GS) | V. Spalević (GS) | P. Ziarati (GS) | A. Yong (GS) | D. Kalaj (GS) | P. Vukoslavčević (GS) | K. Monkova | Lj.R. Mutafchiev (GS) |
| $h$-index | 44 | 41 | 26 | 25 | 31 | 29 | 29 | 26 | 24 | 22 | 14 | 16 | 10 |
| $N_{cit}$ of researcher | 11606 | 6567 | 5407 | 3158 | 3082 | 3031 | 2832 | 2123 | 1631 | 1577 | 912 | 741 | 312 |
| $h \in C$ | No | Yes | No | Yes | Yes | Yes | Yes | Yes | Yes | Yes | Yes | Yes | Yes |
| $h/\overline{C}$ | **0.822** | **0.953** | **0.650** | **0.847** | **1.051** | **0.983** | **0.983** | **1.083** | **1.171** | **1.073** | **0.845** | **1.143** | **1.111** |

From Table 6 we see that $h \in C$ is satisfied for 11 researchers. The estimates $1 < h/\overline{C} < 1.18$ is true for six researchers, while the estimates $0 < 1 - h_d/\overline{C}_d < 0.16$ is true for four researchers.

**Table 7.** The estimation of $h_d$-indices, based on related confidence intervals for 13 researchers (based on their Google Scholar and Scopus profiles, April 2023-June 2023).

| $N_{cit}$ | 4000 | 6000 | 1000 | 1500 | 3000 | 3000 | 3000 | 2000 | 1500 | 1500 | 1000 | 750 | 300 |
|---|---|---|---|---|---|---|---|---|---|---|---|---|---|
| The confidence interval for $h_d$-index, $C_d$ | [29, 39] | [36, 47] | [13, 20] | [17, 24] | [25, 34] | [25, 34] | [25, 34] | [20, 28] | [17, 24] | [17, 24] | [13, 20] | [11, 17] | [7, 11] |
| The mean of $C_d$, $\overline{C}_d$ | 34 | 41.5 | 16.5 | 20.5 | 29.5 | 29.5 | 29.5 | 24 | 20.5 | 20.5 | 16.5 | 14 | 9 |
| Researcher | C.F. Gauss (GS) | G. Andrews (S) | S. Papadimitrou (GS) | D. Zeilberger (S) | I. Orović (GS) | C.D. Savage (GS) | V. Spalević (GS) | P. Ziarati (GS) | A. Yong (GS) | D. Kalaj (GS) | P. Vukoslavčević (GS) | K. Monkova | Lj.R. Mutafchi (GS) |
| $h_d$-index | 34 | 40 | 19 | 20 | 31 | 28 | 29 | 26 | 23 | 22 | 9 | 16 | 10 |
| $h_d + d$ | 52 | 41 | 32 | 30 | 31 | 30 | 29 | 26 | 24 | 22 | 15 | 16 | 10 |
| $N_{cit}^{(d)}$ of researcher | 3765 | 5900 | 964 | 1637 | 3082 | 2279 | 2832 | 2123 | 1520 | 1577 | 272 | 741 | 314 |
| $h_d \in C_d$ | Yes | Yes | Yes | Yes | No | Yes | Yes | Yes | Yes | Yes | No | Yes | Yes |
| $h_d + d \in C_d$ | No | Yes | No | No | Yes | Yes | Yes | Yes | Yes | Yes | Yes | Yes | Yes |
| $h_d/\overline{C}_d$ | **1.000** | **0.964** | **1.152** | **0.976** | **1.051** | **0.949** | **0.983** | **0.923** | **0.854** | **1.073** | **0.545** | **1.143** | **1.111** |
| $(h_d + d)/\overline{C}_d$ | **1.529** | **0.988** | **1.939** | **1.463** | **1.051** | **1.017** | **0.983** | **1.083** | **1.171** | **1.073** | **0.909** | **1.143** | **1.111** |

From Table 7 we see that $h_d \in C_d$ is satisfied for 11 researchers. The estimate $1 \leq h_d/\overline{C}_d < 1.16$ is true for six researchers, while the estimate $0 < 1 - h_d/\overline{C}_d < 0.15$ is also true for six researchers.

As it was noticed in Section 4, Brown (2018) proposed the empirical approximation of the 95% confidence interval for $h$ (the interval (4.24) of Section 4), give as

$$I' = (0.54\sqrt{N_{cit}} - 1.96(0.57 + 0.045\sqrt{N_{cit}}), 0.54\sqrt{N_{cit}} + 1.96(0.57 + 0.045\sqrt{N_{cit}})).$$

Using the expression for the above interval and data from Table 1, we obtain the following table.

**Table 8**. The estimations of $h_d$-indices by using the confidence interval $I'$ (Brown, 2018) for 41 researchers (based on their Google Scholar and Scopus profiles, February 2023-June 2023)

| Researcher | L. Leydes-dorff | W. Glänzel | H.F. Moed | A.F.J. Van Raan | R. Rousseau | A. Schubert | B. Martin | F. Narin | E. Garfield | T. Braun | H. Small | L. Egghe | P. Ingwersen | H.D. White |
|---|---|---|---|---|---|---|---|---|---|---|---|---|---|---|
| $N_{cit}$ | 25005 | 11766 | 7606 | 8308 | 8053 | 7587 | 7598 | 7209 | 11515 | 5680 | 7693 | 5640 | 3606 | 2399 |
| P | 406 | 258 | 127 | 124 | 295 | 141 | 85 | 64 | 106 | 216 | 57 | 211 | 88 | 37 |
| $N_{cit}/p$ | 61,59 | 45.605 | 59.890 | 67.000 | 27.298 | 53.809 | 89.388 | 112.641 | 108.632 | 26.296 | 134.965 | 26.730 | 40.966 | 64.838 |
| $h$ | 79 | 61 | 49 | 48 | 43 | 42 | 38 | 38 | 37 | 37 | 34 | 30 | 27 | 19 |
| $I'$ | (82.9, 87.9) | (56.5, 60.6) | (45.2, 49.0) | (40..1, 58.4) | (46.6, 50.4) | (45.2, 48.9) | (37.9,56.2.7) | (39.554.5) | (55.9, 60. | (38.9, 42.5) | (45.5, 49.3) | (38.8, 42.3) | (28..2, 36.6) | (24.9,0, 28.0) |
| $\overline{I'}$ | 85.4 | 58.6 | 47.1 | 49.2 | 48.5 | 47.1 | 47.1 | 47.0 | 58.0 | 40.7 | 47.4 | 40.6 | 32.4 | 26.5 |
| $h \in I'$ | No | No | Yes | Yes | No | No | Yes | No | No | No | No | No | No | No |

| Researcher | S Freud (GS) | H.J. Kim (GS) | R.C. Kessler (GS) | A. Einstein (GS) | P. Erdös (GS) | T. Tao (GS) | L. Leydesdorf (GS) | C.G. Meyer (GS) | T. Tao (S) | G. Andrews (GS) | M. McAleer (GS) | J. Hirsch (S) | P. Erdös (S) | F.T. Edelmann (GS) |
|---|---|---|---|---|---|---|---|---|---|---|---|---|---|---|
| $N_{cit}$ | 643730 | 518589 | 515591 | 161009 | 99866 | 90963 | 70821 | 49110 | 46852 | 32386 | 26266 | 24451 | 18830 | 13417 |
| $N_{cit}/p$ | 215.222 | 172.863 | 271.364 | 154.223 | 62.106 | 155.759 | 77.655 | 44.045 | 157.508 | 44.364 | 212.284 | 84.606 | 28.021 | 22.140 |
| $h$ | 288 | 338 | 329 | 123 | 128 | 106 | 117 | 95 | 70 | 66 | 79 | 60 | 62 | 54 |
| $I'$ | (361.4, 505.1) | (324.2, 453.5) | (323.3, 452.2) | (180.2, 253.2) | (141.7, 199.6) | (135.1, 190.6) | (119.1, 168.3) | (99.0, 140.3) | (96.7, 137.1) | (80. 2, 114.2) | (72.1, 102.9) | (69.5, 99.3) | (60.9, 87.3) | (51.2, 73.9) |
| $\overline{I'}$ | 433.3 | 388.9 | 387.75 | 216.70 | 170.7 | 162.9 | 143.7 | 119.65 | 116.9 | 194.4 | 87.5 | 84.4 | 74.1 | 62.55 |
| $h \in I'$ | No | Yes | Yes | No | No | No | No | No | No | No | Yes | No | Yes | Yes |

| Researcher | C.F. Gauss (GS) | G. Andrews (S) | S. Papadimitrou (GS) | D. Zeilberger (S) | I. Orović (GS) | C.D. Savage (GS) | V. Spalević (GS) | P. Ziarati (GS) | A. Yong (GS) | D. Kalaj (GS) | P. Vukoslavčević (GS) | K. Monkova | Lj.R. Mutafchiev (GS) |
|---|---|---|---|---|---|---|---|---|---|---|---|---|---|
| $N_{cit}$ | 11606 | 6567 | 5407 | 3158 | 3082 | 3031 | 2832 | 2123 | 1631 | 1577 | 912 | 741 | 312 |
| $N_{cit}/p$ | 140.169 | 25.355 | 58.772 | 18.360 | 17.412 | 29.144 | 6.874 | 11.794 | 16.310 | 9.735 | 14.476 | 6.500 | 5.288 |
| $h$ | 44 | 41 | 26 | 25 | 31 | 29 | 29 | 26 | 24 | 22 | 14 | 16 | 10 |
| $I'$ | (49.8, 66.6) | (37.7, 43.8) | (32.1, 47.3) | (24.3, 36.4) | (24.0, 36.0) | (23.7, 35.7) | (22.7, 34.7) | (19.7, 30.0) | (16.6, 27) | (16.8 26.1) | (12.5, 20.0) | (11.2, 18.2) | (6.9, 12.2) |
| $\overline{I'}$ | 58.2 | 40.8 | 39.7 | 30.4 | 30 | 29.7 | 28.7 | 29.9 | 21.8 | 21.5 | 16.3 | 14.7 | 9.6 |
| $h \in I'$ | No | Yes | No | Yes | Yes | Yes | Yes | Yes | Yes | Yes | Yes | Yes | Yes |

Notice that the values $\overline{I'}$ presented in Table 8 are in fact the values $m_P(n) \approx \sqrt{6}\ln 2\sqrt{N_{cit}}/\pi \approx 0.54044\sqrt{N_{cit}}$, given by the expression (3.3), which are in terms of random partitions the most likely Durfee square sizes for partitions in $P(N_{cit})$ (i.e., the modes of Durfee square size of $P(N_{cit})$) (see Section 3). The only difference between them is that in formula (3.3) for $m_P(n)$, the constant 0.540444639 is used instead of the constant 0.54; namely, $\overline{I'} = 0.54\sqrt{N_{cit}}$.

Moreover, From Table 8 we see that only three Price Awardees have the $h$-index that belongs to the confidence interval $I'$. Among the other 27 researchers, 16 researchers have the $h$-index that belongs to the confidence interval $I'$. Thus, the percentage of the total number of researchers whose $h$-index belongs to the confidence interval $I'$ is equal to $46.3\%$. Observe that the $h$-index do not belong to the interval $I'$ for 9 researchers whose numbers of citations are very large, namely, each of them greater than 2400 (S. Freud (GS), A. Einstein (GS), P. Erdös (GS), T. Tao (GS), L. Leydesdorf (GS), C.G. Meyer (GS), T. Tao (S), G. Andrews (GS) and J. Hirsch (S)). From Table 8 we see that the values $N_{cit}/p$ of all these 9 researchers are relatively large values; namely, they are equal to 215.222, 154.223, 62.106, 155.759, 77.655, 44.045, 157.508, 44.364 and 84.606, respectively. This is consistent with the observation of Brown (2018) that the use of unrestricted partitions (with $p \leq N_{cit}$) to estimate the distribution of $h$-index is inappropriate and this discrepancy grows as ($N_{cit}/p$) increases.

Notice that from Table 5 we see that the values $q := 2N_{cit}(h) - h^2$ of all these 9 researchers are relatively large values; namely, they are equal to 12.496, 18.904, 7.804, 12.018, 7.021, 7.129, 16.221 and 10.734, respectively.

## 6. Concluding Remarks

It is noticed I Abstract of this article, that so far, many researchers have investigated the following question: Given total number of citations, what is the estimated range of the $h$-index? Here we consider the converse question. Namely, the aim of this paper is to estimate the total number of citations of a researcher using only his $h$-index, his $h$-core and mainly a relatively small number of his citations from the tail. For these purposes, we use the asymptotic formula for the mode size of the Durfee square when $n \to \infty$, **which was proved by Canfield, Corteel and Savage (1998), seven years before Hirsch (2005) defined the $h$-index.**

It will be very interesting to see if one could utilize the ideas in this paper to confirm our computational results in terms of random partitions of a positive integer $n := N_{cit}$ as $n \to \infty$.

# Appendix

1) (Tables 1-4 and Table 8). Scopus profiles of 14 Price awardees/leading experts in Scientometrics (February 21 2023; the reduced Table A2 in Meštrović and Dragović, 2023).

Leot Leydesdorff, Wolfgang Glänzel, Henk Moed, Anthony Van Raan, Ronald Rousseau, Andras Schubert, Ben Martin, Francis Narin, Eugene Garfield, Tibor Braun, Henry Small, Leo Egghe, Peter Ingwersen and Howard D. White.

**Table A1.** Price awardees data: the all citations in domain *D* of asymptotic normality of citations (bases on Scopus, 21st February 2023)

| | Price awardees | | | | | | | | | | | | |
|---|---|---|---|---|---|---|---|---|---|---|---|---|---|
| | L. Leyde-sdorff | W. Glänzel | H. Moed | A. Van Raan | R. Rousseau | A. Schubert | B. Martin | F. Narin | E. Garfield | T. Braun | H. Small | L. Egghe | P. Ingwe-Rsen | H.D. White |
| 1 | 3964 | 534 | 417 | 522 | 1050 | 1858 | 1828 | 803 | 1888 | 449 | 2846 | 1455 | 482 | 1029 |
| 2 | 647 | 449 | 417 | 455 | 484 | 449 | 728 | 780 | 1711 | 362 | 592 | 478 | 414 | 213 |
| 3 | 544 | 417 | 380 | 381 | 478 | 362 | 411 | 519 | 1635 | 255 | 540 | 242 | 348 | 178 |
| 4 | 479 | 364 | 350 | 350 | 265 | 273 | 391 | 481 | 597 | 158 | 410 | 206 | 232 | 148 |
| 5 | 424 | 273 | 303 | 303 | 242 | 255 | 371 | 443 | 589 | 151 | 242 | 153 | 166 | 124 |
| 6 | 422 | 255 | 290 | 290 | 175 | 211 | 362 | 354 | 430 | 147 | 197 | 124 | 156 | 123 |
| 7 | 410 | 249 | 288 | 284 | 153 | 159 | 282 | 339 | 264 | 138 | 193 | 110 | 155 | 104 |
| 8 | 401 | 220 | 248 | 283 | 126 | 158 | 244 | 306 | 255 | 135 | 190 | 107 | 133 | 72 |
| 9 | 376 | 211 | 197 | 276 | 118 | 155 | 223 | 283 | 252 | 127 | 186 | 88 | 113 | 54 |
| 10 | 352 | 201 | 183 | 206 | 111 | 151 | 196 | 226 | 238 | 120 | 157 | 81 | 70 | 50 |
| 11 | 352 | 188 | 182 | 172 | 110 | 147 | 179 | 193 | 238 | 119 | 154 | 71 | 57 | 35 |
| 12 | 301 | 163 | 174 | 166 | 97 | 133 | 131 | 183 | 178 | 96 | 151 | 70 | 56 | 32 |
| 13 | 291 | 159 | 166 | 157 | 96 | 127 | 128 | 159 | 155 | 77 | 144 | 70 | 49 | 29 |
| 14 | 266 | 155 | 152 | 152 | 88 | 120 | 122 | 151 | 154 | 75 | 129 | 68 | 46 | 28 |
| 15 | 252 | 147 | 146 | 150 | 87 | 119 | 114 | 138 | 153 | 73 | 113 | 65 | 45 | 25 |
| 16 | 242 | 145 | 122 | 146 | 75 | 105 | 101 | 121 | 151 | 69 | 112 | 50 | 42 | 24 |
| 17 | 241 | 141 | 118 | 142 | 74 | 98 | 77 | 114 | 139 | 67 | 111 | 48 | 42 | 23 |
| 18 | 232 | 135 | 118 | 136 | 72 | 92 | 73 | 94 | 132 | 61 | 102 | 47 | 41 | 21 |
| 19 | 225 | 133 | 117 | 130 | 70 | 89 | 72 | 87 | 114 | 59 | 98 | 46 | 39 | 20 |
| 20 | 221 | 127 | 103 | 111 | 67 | 75 | 69 | 81 | 113 | 57 | 90 | 45 | 39 | 19 |
| 21 | 184 | 126 | 95 | 110 | 66 | 73 | 68 | 76 | 111 | 55 | 81 | 42 | 38 | |
| 22 | 179 | 126 | 92 | 103 | 64 | 73 | 66 | 73 | 91 | 55 | 79 | 41 | 35 | |
| 23 | 171 | 121 | 88 | 98 | 63 | 72 | 65 | 73 | 88 | 52 | 71 | 40 | 35 | |
| 24 | 166 | 111 | 80 | 94 | 61 | 71 | 64 | 68 | 83 | 52 | 57 | 40 | 33 | |
| 25 | 165 | 109 | 76 | 92 | 58 | 70 | 63 | 67 | 74 | 51 | 51 | 39 | 29 | |
| 26 | 157 | 108 | 75 | 86 | 56 | 69 | 58 | 61 | 73 | 47 | 49 | 35 | 29 | |
| 27 | 154 | 96 | 75 | 82 | 56 | 69 | 57 | 61 | 69 | 46 | 47 | 34 | 28 | |
| 28 | 147 | 94 | 73 | 81 | 55 | 64 | 55 | 55 | 65 | 45 | 44 | 34 | 27 | |
| 29 | 147 | 92 | 71 | 81 | 51 | 63 | 50 | 54 | 64 | 45 | 43 | 34 | 27 | |
| 30 | 147 | 89 | 70 | 77 | 49 | 62 | 48 | 50 | 63 | 45 | 41 | 32 | 26 | |
| 31 | 141 | 89 | 69 | 76 | 48 | 62 | 48 | 43 | 63 | 43 | 40 | 30 | 25 | |
| 32 | 138 | 89 | 69 | 75 | 48 | 57 | 47 | 43 | 53 | 42 | 39 | 29 | 24 | |
| 33 | 135 | 85 | 63 | 73 | 47 | 53 | 46 | 43 | 51 | 40 | 38 | 29 | | |
| 34 | 133 | 85 | 62 | 72 | 46 | 52 | 46 | 42 | 49 | 40 | 34 | 28 | | |
| 35 | 133 | 84 | 61 | 72 | 45 | 47 | 45 | 41 | 45 | 39 | 21 | | | |
| 36 | 129 | 84 | 60 | 68 | 45 | 46 | 40 | 40 | 42 | 37 | 21 | | | |
| 37 | 129 | 83 | 60 | 65 | 45 | 46 | 40 | 40 | 39 | 37 | 20 | | | |
| 38 | 127 | 80 | 58 | 63 | 45 | 45 | 40 | 38 | 30 | 36 | 19 | | | |
| 39 | 127 | 80 | 57 | 62 | 44 | 45 | **32** | 37 | 28 | 34 | 18 | | | |
| 40 | 126 | 78 | 57 | 60 | 44 | 43 | 32 | 34 | 27 | | 17 | | | |
| 41 | 125 | 76 | 56 | 59 | 43 | 41 | 31 | 32 | 27 | | 17 | | | |
| 42 | 125 | 76 | 55 | 59 | 43 | 42 | 28 | 27 | 25 | | 15 | | | |
| 43 | 125 | 75 | 53 | 59 | 43 | 40 | 24 | 24 | 23 | | 14 | | | |
| 44 | 125 | 75 | 53 | 55 | 42 | 37 | 23 | 21 | 22 | | | | | |
| 45 | 124 | 75 | 52 | 52 | 42 | 36 | | 21 | 22 | | | | | |
| 46 | 123 | 74 | 52 | 51 | | 36 | | | 21 | | | | | |
| 47 | 122 | 73 | 50 | 48 | | | | | 21 | | | | | |
| 48 | 122 | 73 | 49 | 48 | | | | | | | | | | |
| 49 | 118 | 73 | 49 | 47 | | | | | | | | | | |
| 50 | 108 | 72 | 47 | 45 | | | | | | | | | | |
| 51 | 106 | 72 | 45 | 41 | | | | | | | | | | |
| 52 | 102 | 71 | | 41 | | | | | | | | | | |
| 53 | 100 | 71 | | 39 | | | | | | | | | | |
| 54 | 98 | 70 | | | | | | | | | | | | |
| 55 | 97 | 69 | | | | | | | | | | | | |
| 56 | 94 | 66 | | | | | | | | | | | | |
| 57 | 94 | 64 | | | | | | | | | | | | |
| 58 | 93 | 63 | | | | | | | | | | | | |
| 59 | 93 | 62 | | | | | | | | | | | | |
| 60 | 92 | 62 | | | | | | | | | | | | |
| 61 | 91 | 62 | | | | | | | | | | | | |
| 62 | 90 | 59 | | | | | | | | | | | | |

**Table A1** – Continued (only for Leydesdorf). Price awardees data: the all citations in domain $D$ of asymptotic normality of citations (bases on Scopus, 21$^{st}$ February 2023)

| L. Leyde-sdorff | 63 | 64 | 65 | 66 | 67 | 68 | 69 | 70 | 71 | 72 | 73 | 74 | 75 | 76 | 77 | 78 | 79 | 80 | 81 | 82 | 83 |
|---|---|---|---|---|---|---|---|---|---|---|---|---|---|---|---|---|---|---|---|---|---|
|  | 89 | 87 | 87 | 87 | 87 | 86 | 83 | 83 | 83 | 82 | 82 | 82 | 80 | 80 | 80 | 80 | **79** | 79 | 78 | 78 | 77 |

Note:

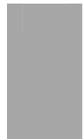 $h$-index  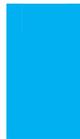 The bounds of domain $D$ of asymptotic normality of citations  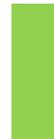 The values for which $h_d + d = h$

2) (Tables 5-7 and Table 8). 27 research profiles at Google Scholar (GS) and Scopus (S) databases of 24 researchers from various scientific fields, on April-June 2023 (which include 24 different researchers, with being considered both Scopus and Google Scholar profiles for Paul Erdös, Terence Tao and George Andrews).

Sigmund Freud (GS, June 26, Neurology, Psyhology, Psyhopatology), H. Jihn Kim (GS, May 26, High energy physics, astro-particle physics, scintillator), Ronald C. Kessler (GS, May 24, Psychiatric Epidemiology), Albert Einstein (GS, May 21, Physics), Paul Erdös (GS, April 20, Number Theory, Combinatorics, Probability, Set Theory, Mathematical Analysis), Terence Tao (GS, April 20, Mathematical analysis, Combinatorics, Random Matrix Theory, PDE, Number theory), L. Leydesdorf (GS, May 24, Bibliometrics and Scientometrics), Christian G. Meyer (GS, Infectious diseases), Terence Tao (S, May 20, Mathematical analysis, Combinatorics, Random matrix theory, PDE, Number theory), George Andrews (GS, May 05, Mathematical analysis, Combinatorics, Number theory), Michael McAleer (GS, May 27, Economics, Econometrics, Financial economics, Financial econometrics), Jorge E. Hirsch (S, May 23, Physics, Supercodutivity, Ferromagnetism), Paul Erdös (S, April 21, Number theory, Combinatorics, Probability, Set theory, Mathematical analysis), Frank. T. Edelmann (May 27, GS, Inorganic chemistry), Carl Friedrich Gauss (GS, May 28, Mathematics), George Andrews (S, May 23, Mathematical analysis, Combinatorics, Number theory), Stratos Papadimitrou (GS, May 27, Transportation, Maritime logistics), Doron Zeilberger (S, May 23, Combinatorics, Hypergeometric identities, Symbolic computation), Irena Orović (GS, May 27, Signal processing, Multimedia signals and systems, Compressive sensing, Time-frequency signal analysis), Carla. D. Savage (GS, June 3, Algorithms, Combinatorics), Velibor Spalević (GS, May 25, Sediment yield, Soil erosion modeling, Erosion, Watershed management, IntErO), Parisa Ziarati (GS, June 26, Toxicology, Green chemistry, Enviromental chemistry, Inorganic chemistry), Alexander Yong (GS, May 23, Combinatorics, Algebraic geometry, Representation theory), David Kalaj, (GS, May 23, Complex analysis), Petar Vukoslavčević (GS, May 28, Mechanical engineering, Measurement, Turbulent flow, Boundary layer, Applied Mathematics), Katarina Monkova (S, May 23, Mechanical engineering), Ljuben R. Mutafchiev (GS, June 3, Random combinatorial structures).

Table A2. Data of 23 reseacher profiles: the all citations in domain $D$ of asymptotic normality of citations (bases on Google Scholar (GS) and Scopus (S), April, May and June 2023).

| | A. Einstein (GS) | P. Erdös (GS) | T. Tao (GS) | L. Leydesdorf (GS) | C.G. Meyer (GS) | G. Andrews (GS) | M. McAleer (GS) | J. Hirsch (S) | P. Erdös (S) | F.T. Edelmann (GS) | C.F. Gauss (GS) |
|---|---|---|---|---|---|---|---|---|---|---|---|
| 1 | 23200 | 23995 | 18709 | 12829 | 11973 | 6906 | 1108 | 7159 | 890 | 602 | 1953 |
| 2 | 20218 | 2070 | 12037 | 3505 | 3537 | 5653 | 617 | 2378 | 610 | 529 | 1417 |
| 3 | 7572 | 1955 | 8541 | 3269 | 1509 | 1352 | 553 | 881 | 519 | 474 | 534 |
| 4 | 7216 | 1769 | 8325 | 1738 | 697 | 1299 | 465 | 710 | 498 | 471 | 533 |
| 5 | 7029 | 1268 | 2452 | 1395 | 690 | 1108 | 434 | 686 | 437 | 366 | 510 |
| 6 | 6124 | 1254 | 1966 | 1346 | 524 | 731 | 429 | 659 | 353 | 346 | 453 |
| 7 | 4674 | 1201 | 1902 | 1170 | 512 | 587 | 414 | 504 | 344 | 331 | 368 |
| 8 | 4309 | 1131 | 1687 | 1089 | 476 | 585 | 378 | 440 | 335 | 236 | 323 |
| 9 | 4085 | 1121 | 1233 | 841 | 468 | 381 | 371 | 369 | 330 | 232 | 313 |
| 10 | 4007 | 931 | 1023 | 765 | 412 | 329 | 365 | 364 | 301 | 175 | 220 |
| 11 | 3787 | 911 | 710 | 725 | 404 | 315 | 355 | 361 | 273 | 147 | 205 |
| 12 | 3505 | 888 | 599 | 704 | 393 | 271 | 353 | 275 | 257 | 138 | 176 |
| 13 | 2536 | 885 | 590 | 699 | 390 | 266 | 336 | 263 | 257 | 137 | 147 |
| 14 | 2323 | 878 | 567 | 697 | 361 | 249 | 335 | 249 | 243 | 130 | 146 |
| 15 | 1813 | 762 | 529 | 695 | 354 | 221 | 334 | 224 | 230 | 129 | 144 |
| 16 | 1795 | 691 | 522 | 655 | 350 | 213 | 331 | 215 | 220 | 129 | 141 |
| 17 | 1687 | 653 | 473 | 606 | 348 | 208 | 287 | 209 | 209 | 128 | 131 |
| 18 | 1531 | 609 | 453 | 553 | 344 | 193 | 264 | 205 | 196 | 96 | 127 |
| 19 | 1516 | 586 | 450 | 544 | 329 | 172 | 262 | 202 | 195 | 96 | 116 |
| 20 | 1398 | 586 | 408 | 531 | 320 | 157 | 241 | 191 | 191 | 93 | 116 |
| 21 | 1389 | 575 | 401 | 525 | 318 | 149 | 222 | 186 | 177 | 90 | 115 |
| 22 | 1349 | 561 | 371 | 507 | 285 | 148 | 220 | 175 | 177 | 87 | 114 |
| 23 | 1279 | 560 | 370 | 504 | 276 | 145 | 213 | 159 | 175 | 84 | 108 |
| 24 | 1185 | 558 | 370 | 485 | 271 | 144 | 211 | 154 | 150 | 83 | 106 |
| 25 | 1153 | 545 | 334 | 427 | 264 | 141 | 197 | 153 | 149 | 83 | 97 |
| 26 | 1148 | 541 | 296 | 378 | 256 | 131 | 191 | 139 | 146 | 78 | 96 |
| 27 | 1119 | 488 | 289 | 342 | 253 | 127 | 189 | 129 | 138 | 78 | 96 |
| 28 | 1091 | 478 | 289 | 334 | 251 | 121 | 181 | 126 | 137 | 76 | 90 |
| 29 | 1026 | 477 | 278 | 329 | 244 | 107 | 180 | 121 | 137 | 76 | 90 |
| 30 | 995 | 473 | 276 | 329 | 244 | 105 | 174 | 119 | 135 | 74 | 83 |
| 31 | 984 | 468 | 276 | 328 | 239 | 104 | 171 | 119 | 135 | 73 | 81 |
| 32 | 974 | 463 | 276 | 328 | 236 | 102 | 170 | 107 | 125 | 73 | 75 |
| 33 | 863 | 451 | 270 | 323 | 232 | 101 | 167 | 106 | 120 | 72 | 74 |
| 34 | 860 | 450 | 269 | 318 | 232 | 100 | 163 | 105 | 115 | 72 | 68 |
| 35 | 832 | 429 | 256 | 312 | 230 | 100 | 157 | 100 | 113 | 71 | 66 |
| 36 | 746 | 407 | 254 | 302 | 225 | 99 | 154 | 99 | 110 | 70 | 62 |
| 37 | 668 | 399 | 253 | 300 | 224 | 98 | 149 | 94 | 110 | 69 | 61 |
| 38 | 567 | 399 | 249 | 299 | 222 | 97 | 146 | 93 | 108 | 69 | 60 |
| 39 | 558 | 388 | 247 | 297 | 215 | 97 | 145 | 93 | 105 | 69 | 58 |
| 40 | 520 | 377 | 242 | 288 | 200 | 95 | 143 | 90 | 101 | 66 | 57 |
| 41 | 518 | 373 | 239 | 287 | 198 | 95 | 143 | 86 | 99 | 66 | 53 |
| 42 | 507 | 363 | 239 | 271 | 198 | 93 | 139 | 86 | 97 | 64 | 47 |
| 43 | 500 | 358 | 235 | 270 | 193 | 91 | 135 | 85 | 88 | 64 | 45 |
| 44 | 487 | 357 | 233 | 265 | 190 | 91 | 134 | 84 | 86 | 64 | **45** |
| 45 | 483 | 354 | 233 | 248 | 190 | 91 | 133 | 82 | 86 | 64 | 42 |
| 46 | 471 | 352 | 232 | 243 | 178 | 90 | 131 | 82 | 85 | 64 | 41 |
| 47 | 440 | 332 | 228 | 242 | 175 | 84 | 130 | 81 | 83 | 63 | 39 |
| 48 | 436 | 326 | 227 | 240 | 173 | 84 | 130 | 80 | 83 | 59 | 39 |
| 49 | 436 | 323 | 226 | 236 | 171 | 81 | 125 | 78 | 81 | 58 | 38 |
| 50 | 431 | 306 | 225 | 235 | 170 | 80 | 124 | 78 | 78 | 57 | 36 |
| 51 | 431 | 302 | 223 | 234 | 165 | 80 | 124 | 77 | 78 | 56 | 36 |
| 52 | 396 | 294 | 220 | 232 | 164 | 80 | 121 | 72 | 75 | 56 | 36 |
| 53 | 384 | 292 | 217 | 231 | 164 | 77 | 120 | 71 | 75 | 55 | 34 |
| 54 | 382 | 287 | 213 | 220 | 163 | 77 | 119 | 69 | 75 | **55** | |
| 55 | 380 | 286 | 212 | 212 | 163 | 77 | 111 | 68 | 70 | 54 | |
| 56 | 367 | 279 | 211 | 211 | 162 | 75 | 109 | 68 | 70 | 54 | |
| 57 | 361 | 274 | 207 | 211 | 157 | 74 | 108 | 65 | 69 | 54 | |
| 58 | 357 | 268 | 205 | 211 | 156 | 74 | 107 | 62 | 67 | 54 | |
| 59 | 356 | 268 | 203 | 209 | 154 | 73 | 101 | 62 | 67 | | |

| | | | | | | | | | |
|---|---|---|---|---|---|---|---|---|---|
| 60 | 355 | 261 | 203 | 209 | 147 | 73 | 98 | 60 | 66 |
| 61 | 330 | 261 | 202 | 208 | 146 | 70 | 95 | 60 | 64 |
| 62 | 311 | 255 | 201 | 206 | 139 | 70 | 94 | 59 | **64** |
| 63 | 310 | 250 | 200 | 206 | 137 | 69 | 92 | 59 | 61 |
| 64 | 304 | 246 | 176 | 205 | 135 | 67 | 92 | 59 | 61 |
| 65 | 294 | 245 | 174 | 203 | 127 | 67 | 91 | 59 | 60 |
| 66 | 285 | 245 | 174 | 202 | 125 | **67** | 91 | 57 | 60 |
| 67 | 280 | 244 | 174 | 199 | 125 | 65 | 90 | 56 | 59 |
| 68 | 276 | 242 | 167 | 197 | 122 | 64 | 88 | | 58 |
| 69 | 275 | 242 | 166 | 196 | 121 | 63 | 88 | | 58 |
| 70 | 264 | 236 | 160 | 196 | 121 | 63 | 88 | | |
| 71 | 258 | 232 | 160 | 194 | 121 | | 87 | | |
| 72 | 257 | 229 | 158 | 191 | 121 | | 86 | | |
| 73 | 254 | 221 | 154 | 190 | 120 | | 86 | | |
| 74 | 252 | 221 | 154 | 185 | 118 | | 84 | | |
| 75 | 251 | 219 | 154 | 183 | 117 | | 84 | | |
| 76 | 244 | 219 | 153 | 177 | 117 | | 83 | | |
| 77 | 233 | 216 | 148 | 173 | 117 | | 81 | | |
| 78 | 232 | 211 | 148 | 173 | 116 | | 81 | | |
| 79 | 228 | 210 | 148 | 172 | 116 | | **81** | | |
| 80 | 222 | 206 | 148 | 171 | 116 | | 77 | | |
| 81 | 218 | 204 | 145 | 167 | 115 | | 76 | | |
| 82 | 216 | 203 | 143 | 166 | 113 | | 74 | | |
| 83 | 216 | 195 | 142 | 164 | 110 | | 73 | | |
| 84 | 212 | 194 | 141 | 163 | 106 | | 72 | | |
| 85 | 212 | 191 | 139 | 163 | 106 | | | | |
| 86 | 209 | 191 | 139 | 161 | 105 | | | | |
| 87 | 208 | 191 | 138 | 156 | 101 | | | | |
| 88 | 201 | 186 | 137 | 154 | 100 | | | | |
| 89 | 199 | 178 | 136 | 150 | 100 | | | | |
| 90 | 194 | 177 | 133 | 148 | 99 | | | | |
| 91 | 185 | 175 | 133 | 148 | 97 | | | | |
| 92 | 182 | 174 | 132 | 147 | 96 | | | | |
| 93 | 181 | 173 | 131 | 145 | 95 | | | | |
| 94 | 176 | 171 | 127 | 145 | 95 | | | | |
| 95 | 176 | 168 | 125 | 145 | **95** | | | | |
| 96 | 172 | 166 | 125 | 143 | 93 | | | | |
| 97 | 171 | 165 | 124 | 142 | 93 | | | | |
| 98 | 170 | 165 | 123 | 140 | 92 | | | | |
| 99 | 169 | 165 | 122 | 137 | 92 | | | | |
| 100 | 167 | 164 | 121 | 137 | 91 | | | | |
| 101 | 166 | 164 | 119 | 136 | 91 | | | | |
| 102 | 163 | 163 | 118 | 136 | | | | | |
| 103 | 157 | 162 | 115 | 135 | | | | | |
| 104 | 147 | 161 | 115 | 133 | | | | | |
| 105 | 147 | 157 | 111 | 131 | | | | | |
| 106 | 147 | 157 | **108** | 128 | | | | | |
| 107 | 146 | 157 | 106 | 128 | | | | | |
| 108 | 145 | 155 | 105 | 126 | | | | | |
| 109 | 145 | 147 | 105 | 125 | | | | | |
| 110 | 144 | 144 | 103 | 123 | | | | | |
| 111 | 144 | 141 | 103 | 122 | | | | | |
| 112 | 144 | 140 | 102 | 121 | | | | | |
| 113 | 143 | 140 | 102 | 119 | | | | | |
| 114 | 141 | 139 | 101 | 119 | | | | | |
| 115 | 136 | 139 | 101 | 118 | | | | | |
| 116 | 135 | 138 | 100 | 117 | | | | | |
| 117 | 134 | 138 | | **117** | | | | | |
| 118 | 134 | 136 | | 117 | | | | | |
| 119 | 133 | 135 | | 116 | | | | | |
| 120 | 130 | 134 | | 115 | | | | | |
| 121 | 125 | 133 | | 115 | | | | | |
| 122 | 123 | 132 | | 113 | | | | | |
| 123 | **123** | 131 | | 113 | | | | | |
| 124 | 119 | 131 | | 113 | | | | | |

| | | | | |
|---|---|---|---|---|
| 125 | 117 | 130 | | 110 |
| 126 | 116 | 129 | | |
| 127 | 115 | 129 | | |
| 128 | 110 | **128** | | |
| 129 | 108 | 127 | | |
| 130 | 107 | 126 | | |
| 131 | 106 | 125 | | |
| 132 | 106 | 124 | | |
| 133 | 103 | 124 | | |
| 134 | 101 | 124 | | |
| 135 | 100 | 124 | | |
| 136 | 100 | 123 | | |
| 137 | 98 | 122 | | |
| 138 | 98 | 122 | | |
| 139 | 98 | 122 | | |
| 140 | 92 | 122 | | |
| 141 | 91 | 121 | | |
| 142 | 90 | | | |
| 143 | 90 | | | |
| 144 | 89 | | | |
| 145 | 88 | | | |
| 146 | 87 | | | |

**Table A2-Continued.** Data of 23 researcher profiles: the all citations in domain $D$ of asymptotic normality of citations (bases on Google Scholar (GS) and Scopus (S), April, May and June 2023).

| | G. Andrews (S) | S. Papadimitriou (GS) | D. Zeilberger (S) | I. Orović (GS) | C.D. Savage (GS) | V. Spalević (GS) | P. Ziarati (GS) | A. Yong (GS) | D. Kalaj (GS) | P. Vukoslavčević (GS) | K. Monkova (GS) | Lj.R. Mutafchiev (GS) |
|---|---|---|---|---|---|---|---|---|---|---|---|---|
| 1 | 667 | 754 | 295 | 211 | 581 | 164 | 70 | 111 | 95 | 213 | 62 | 29 |
| 2 | 286 | 509 | 251 | 169 | 171 | 64 | 66 | 86 | 95 | 125 | 41 | 22 |
| 3 | 227 | 486 | 224 | 116 | 127 | 57 | 63 | 84 | 92 | 122 | 35 | 21 |
| 4 | 169 | 450 | 164 | 114 | 123 | 51 | 62 | 82 | 91 | 73 | 35 | 18 |
| 5 | 165 | 440 | 136 | 105 | 92 | 51 | 57 | 81 | 69 | 63 | 32 | 17 |
| 6 | 165 | 354 | 98 | 102 | 88 | 49 | 53 | 79 | 54 | 44 | 30 | 16 |
| 7 | 131 | 284 | 95 | 87 | 82 | 46 | 52 | 67 | 41 | 34 | 27 | 14 |
| 8 | 127 | 263 | 88 | 84 | 70 | 43 | 51 | 61 | 40 | 33 | 25 | 12 |
| 9 | 113 | 237 | 86 | 67 | 64 | 41 | 42 | 49 | 33 | 25 | 24 | 12 |
| 10 | 98 | 211 | 84 | 64 | 63 | 40 | 40 | 44 | 32 | 22 | 22 | 11 |
| 11 | 96 | 186 | 81 | 62 | 60 | 39 | 39 | 44 | 31 | 16 | 20 | 10 |
| 12 | 92 | 153 | 78 | 56 | 58 | 38 | 39 | 42 | 30 | 15 | 18 | |
| 13 | 83 | 116 | 67 | 53 | 55 | 37 | 34 | 37 | 29 | 15 | 18 | |
| 14 | 77 | 101 | 61 | 53 | 53 | 37 | 33 | 37 | 28 | **14** | 17 | |
| 15 | 77 | 86 | 57 | 53 | 52 | 36 | 29 | 37 | 27 | 11 | 17 | |
| 16 | 75 | 66 | 54 | 51 | 50 | 35 | 28 | 35 | **25** | 9 | 17 | |
| 17 | 71 | 54 | 53 | 47 | 49 | 34 | 28 | 24 | 25 | | 16 | |
| 18 | 65 | 47 | 38 | 45 | 47 | 34 | 28 | 32 | 22 | | | |
| 19 | 64 | 45 | 35 | 41 | 43 | 34 | 27 | 32 | 22 | | | |
| 20 | 63 | 34 | 35 | 41 | 42 | 33 | 27 | 31 | 22 | | | |
| 21 | 58 | 34 | 32 | 39 | 42 | 32 | 27 | 28 | 22 | | | |
| 22 | 57 | 33 | 31 | 38 | 38 | 32 | 27 | 25 | 22 | | | |
| 23 | 54 | 32 | 28 | 36 | 37 | 31 | 26 | 24 | | | | |
| 24 | 52 | 30 | 26 | 36 | 36 | 31 | 26 | 24 | | | | |
| 25 | 52 | 29 | **25** | 34 | 35 | 31 | 26 | | | | | |
| 26 | 52 | **26** | 24 | 32 | 35 | 31 | 26 | | | | | |
| 27 | 51 | 25 | 24 | 32 | 35 | 30 | 26 | | | | | |
| 28 | 51 | 21 | 23 | 32 | 35 | 29 | | | | | | |
| 29 | 50 | 21 | 23 | 31 | **32** | 29 | | | | | | |
| 30 | 48 | 21 | 21 | 31 | 29 | 28 | | | | | | |
| 31 | 47 | 20 | 20 | 31 | 28 | | | | | | | |
| 32 | 47 | 20 | | 31 | | | | | | | | |
| 33 | 47 | | | | | | | | | | | |
| 34 | 46 | | | | | | | | | | | |

| | |
|---|---|
| 35 | 45 |
| 36 | 43 |
| 37 | 42 |
| 38 | 42 |
| 39 | 41 |
| 40 | 41 |
| 41 | 41 |
| 42 | 39 |

Note:

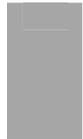 $h$-index

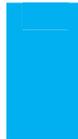 The bounds of domain $D$ of asymptotic normality of citations

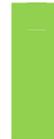 The values for which $h_d + d = h$

**Table A3.** The all citations of Sigmund Freud in domain $D$ of asymptotic normality of citations (Google Scholar, June 26 2023)

| # | cit | # | cit | # | cit | # | cit | # | cit | # | cit |
|---|---|---|---|---|---|---|---|---|---|---|---|
| 1 | 83366 | 61 | 1760 | 121 | 803 | 181 | 525 | 241 | 377 | 301 | 273 |
| 2 | 28948 | 62 | 1751 | 122 | 798 | 182 | 519 | 242 | 376 | 302 | 271 |
| 3 | 22542 | 63 | 1732 | 123 | 793 | 183 | 515 | 243 | 375 | 303 | 271 |
| 4 | 19196 | 64 | 1617 | 124 | 778 | 184 | 514 | 244 | 372 | 304 | 268 |
| 5 | 15892 | 65 | 1615 | 125 | 774 | 185 | 511 | 245 | 372 | 305 | 265 |
| 6 | 15401 | 66 | 1534 | 126 | 771 | 186 | 507 | 246 | 371 | 306 | 265 |
| 7 | 15321 | 67 | 1521 | 127 | 758 | 187 | 503 | 247 | 368 | 307 | 265 |
| 8 | 9847 | 68 | 1518 | 128 | 754 | 188 | 500 | 248 | 364 | 308 | 265 |
| 9 | 9602 | 69 | 1510 | 129 | 752 | 189 | 498 | 249 | 359 | 309 | 262 |
| 10 | 9297 | 70 | 1451 | 130 | 748 | 190 | 497 | 250 | 356 | 310 | 261 |
| 11 | 8686 | 71 | 1449 | 131 | 747 | 191 | 493 | 251 | 356 | 311 | 260 |
| 12 | 8653 | 72 | 1352 | 132 | 732 | 192 | 491 | 252 | 355 | 312 | 259 |
| 13 | 7870 | 73 | 1336 | 133 | 730 | 193 | 488 | 253 | 353 | 313 | 259 |
| 14 | 7466 | 74 | 1321 | 134 | 728 | 194 | 474 | 254 | 350 | 314 | 258 |
| 15 | 7281 | 75 | 1317 | 135 | 723 | 195 | 474 | 255 | 347 | 315 | 258 |
| 16 | 7248 | 76 | 1306 | 136 | 718 | 196 | 473 | 256 | 346 | 316 | 255 |
| 17 | 7094 | 77 | 1298 | 137 | 715 | 197 | 467 | 257 | 344 | 317 | 254 |
| 18 | 6435 | 78 | 1263 | 138 | 710 | 198 | 464 | 258 | 342 | 318 | 253 |
| 19 | 5391 | 79 | 1241 | 139 | 701 | 199 | 457 | 259 | 338 | 319 | 253 |
| 20 | 5356 | 80 | 1236 | 140 | 696 | 200 | 454 | 260 | 338 | 320 | 252 |
| 21 | 5137 | 81 | 1199 | 141 | 693 | 201 | 453 | 261 | 332 | 321 | 252 |
| 22 | 5125 | 82 | 1191 | 142 | 690 | 202 | 449 | 262 | 324 | 322 | 250 |
| 23 | 4419 | 83 | 1190 | 143 | 684 | 203 | 445 | 263 | 323 | 323 | 249 |
| 24 | 4239 | 84 | 1187 | 144 | 684 | 204 | 442 | 264 | 322 | 324 | 248 |
| 25 | 4227 | 85 | 1177 | 145 | 662 | 205 | 441 | 265 | 321 | 325 | 247 |
| 26 | 3999 | 86 | 1165 | 146 | 661 | 206 | 436 | 266 | 320 | 326 | 247 |
| 27 | 3249 | 87 | 1157 | 147 | 656 | 207 | 434 | 267 | 319 | 327 | 245 |
| 28 | 3197 | 88 | 1156 | 148 | 645 | 208 | 431 | 268 | 316 | 328 | 245 |
| 29 | 3146 | 89 | 1150 | 149 | 632 | 209 | 430 | 269 | 316 | 329 | 244 |
| 30 | 3057 | 90 | 1149 | 150 | 631 | 210 | 429 | 270 | 316 | 330 | 244 |
| 31 | 2989 | 91 | 1118 | 151 | 630 | 211 | 427 | 271 | 313 | 331 | 243 |
| 32 | 2974 | 92 | 1115 | 152 | 626 | 212 | 424 | 272 | 310 | 332 | 241 |
| 33 | 2944 | 93 | 1098 | 153 | 626 | 213 | 423 | 273 | 310 | 333 | 239 |
| 34 | 2919 | 94 | 1091 | 154 | 613 | 214 | 422 | 274 | 310 | 334 | 239 |
| 35 | 2716 | 95 | 1081 | 155 | 611 | 215 | 417 | 275 | 309 | | |
| 36 | 2671 | 96 | 1080 | 156 | 607 | 216 | 416 | 276 | 308 | | |
| 37 | 2652 | 97 | 1075 | 157 | 604 | 217 | 412 | 277 | 306 | | |
| 38 | 2639 | 98 | 1068 | 158 | 596 | 218 | 411 | 278 | 305 | | |
| 39 | 2631 | 99 | 1040 | 159 | 595 | 219 | 409 | 279 | 305 | | |
| 40 | 2625 | 100 | 1024 | 160 | 592 | 220 | 409 | 280 | 304 | | |
| 41 | 2511 | 101 | 1022 | 161 | 588 | 221 | 405 | 281 | 303 | | |
| 42 | 2497 | 102 | 996 | 162 | 573 | 222 | 403 | 282 | 296 | | |
| 43 | 2469 | 103 | 987 | 163 | 566 | 223 | 403 | 283 | 293 | | |
| 44 | 2373 | 104 | 968 | 164 | 562 | 224 | 403 | 284 | 293 | | |
| 45 | 2351 | 105 | 967 | 165 | 561 | 225 | 403 | 285 | 292 | | |
| 46 | 2278 | 106 | 924 | 166 | 557 | 226 | 402 | 286 | 291 | | |
| 47 | 2262 | 107 | 921 | 167 | 555 | 227 | 402 | 287 | 288 | | |
| 48 | 2170 | 108 | 906 | 168 | 552 | 228 | 400 | 288 | 288 | | |
| 49 | 2105 | 109 | 897 | 169 | 551 | 229 | 399 | 289 | 287 | | |
| 50 | 2046 | 110 | 873 | 170 | 549 | 230 | 399 | 290 | 284 | | |
| 51 | 2032 | 111 | 867 | 171 | 548 | 231 | 394 | 291 | 284 | | |
| 52 | 1994 | 112 | 853 | 172 | 547 | 232 | 389 | 292 | 283 | | |
| 53 | 1992 | 113 | 850 | 173 | 546 | 233 | 388 | 293 | 281 | | |
| 54 | 1948 | 114 | 843 | 174 | 540 | 234 | 388 | 294 | 280 | | |
| 55 | 1935 | 115 | 838 | 175 | 539 | 235 | 387 | 295 | 280 | | |
| 56 | 1934 | 116 | 827 | 176 | 536 | 236 | 386 | 296 | 278 | | |
| 57 | 1893 | 117 | 827 | 177 | 533 | 237 | 385 | 297 | 276 | | |
| 58 | 1873 | 118 | 809 | 178 | 531 | 238 | 384 | 298 | 276 | | |
| 59 | 1868 | 119 | 806 | 179 | 527 | 239 | 383 | 299 | 276 | | |
| 60 | 1800 | 120 | 805 | 180 | 526 | 240 | 381 | 300 | 275 | | |

Note:

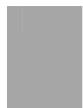 $h$-index ($h = 288$)

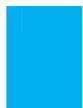 The bounds of domain $D$ of asymptotic normality of citations

**Table A4.** The all citations of Terence Tao and the bounds of domain *D* of asymptotic normality of citations (Scopus, April 17 2023)

| | | | | | | | | | |
|---|---|---|---|---|---|---|---|---|---|
| 1 | 12745 | 61 | 80 | 121 | 33 | 181 | 17 | 241 | 5 |
| 2 | 5548 | 62 | 80 | 122 | 32 | 182 | 17 | 242 | 5 |
| 3 | 5408 | 63 | 78 | 123 | 32 | 183 | 16 | 243 | 5 |
| 4 | 5407 | 64 | 78 | 124 | 32 | 184 | 16 | 244 | 5 |
| 5 | 2053 | 65 | 77 | 125 | 32 | 185 | 16 | 245 | 5 |
| 6 | 1391 | 66 | 76 | 126 | 31 | 186 | 16 | 246 | 4 |
| 7 | 1099 | 67 | 75 | 127 | 31 | 187 | 16 | 247 | 4 |
| 8 | 450 | 68 | 72 | 128 | 31 | 188 | 16 | 248 | 4 |
| 9 | 356 | 69 | 70 | 129 | 30 | 189 | 15 | 249 | 4 |
| 10 | 278 | **70** | **70** | 130 | 30 | 190 | 15 | 250 | 4 |
| 11 | 268 | 71 | 67 | 131 | 30 | 191 | 15 | 251 | 4 |
| 12 | 267 | 72 | 67 | 132 | 30 | 192 | 15 | 252 | 4 |
| 13 | 255 | 73 | 67 | 133 | 30 | 193 | 15 | 253 | 3 |
| 14 | 233 | 74 | 65 | 134 | 29 | 194 | 15 | 254 | 3 |
| 15 | 231 | 75 | 65 | 135 | 29 | 195 | 15 | 255 | 3 |
| 16 | 220 | 76 | 63 | 136 | 29 | 196 | 14 | 256 | 3 |
| 17 | 212 | 77 | 63 | 137 | 28 | 197 | 14 | 257 | 3 |
| 18 | 194 | 78 | 63 | 138 | 28 | 198 | 13 | 258 | 3 |
| 19 | 184 | 79 | 61 | 139 | 28 | 199 | 13 | 259 | 3 |
| 20 | 163 | 80 | 61 | 140 | 28 | 200 | 13 | 260 | 3 |
| 21 | 159 | 81 | 60 | 141 | 28 | 201 | 13 | 261 | 3 |
| 22 | 150 | 82 | 57 | 142 | 27 | 202 | 12 | 262 | 2 |
| 23 | 145 | 83 | 57 | 143 | 27 | 203 | 12 | 263 | 2 |
| 24 | 141 | 84 | 57 | 144 | 27 | 204 | 12 | 264 | 2 |
| 25 | 141 | 85 | 56 | 145 | 26 | 205 | 12 | 265 | 2 |
| 26 | 137 | 86 | 55 | 146 | 25 | 206 | 12 | 266 | 2 |
| 27 | 136 | 87 | 54 | 147 | 25 | 207 | 12 | 267 | 2 |
| 28 | 135 | 88 | 53 | 148 | 25 | 208 | 12 | 268 | 2 |
| 29 | 133 | 89 | 52 | 149 | 24 | 209 | 11 | 269 | 2 |
| 30 | 133 | 90 | 50 | 150 | 24 | 210 | 11 | 270 | 2 |
| 31 | 127 | 91 | 48 | 151 | 23 | 211 | 11 | 271 | 2 |
| 32 | 125 | 92 | 47 | 152 | 23 | 212 | 11 | 272 | 1 |
| 33 | 125 | 93 | 47 | 153 | 23 | 213 | 11 | 273 | 1 |
| 34 | 120 | 94 | 47 | 154 | 23 | 214 | 11 | 274 | 1 |
| 35 | 120 | 95 | 46 | 155 | 23 | 215 | 11 | 275 | 1 |
| 36 | 117 | 96 | 46 | 156 | 22 | 216 | 11 | 276 | 1 |
| 37 | 115 | 97 | 45 | 157 | 22 | 217 | 10 | 277 | 1 |
| 38 | 115 | 98 | 45 | 158 | 22 | 218 | 10 | 278 | 1 |
| 39 | 114 | 99 | 44 | 159 | 21 | 219 | 9 | 279 | |
| 40 | 112 | 100 | 44 | 160 | 21 | 220 | 9 | 280 | |
| 41 | 112 | 101 | 43 | 161 | 21 | 221 | 9 | | |
| 42 | 108 | 102 | 42 | 162 | 20 | 222 | 9 | | |
| 43 | 106 | 103 | 42 | 163 | 20 | 223 | 9 | | |
| 44 | 105 | 104 | 42 | 164 | 20 | 224 | 8 | | |
| 45 | 104 | 105 | 41 | 165 | 20 | 225 | 8 | | |
| 46 | 101 | 106 | 40 | 166 | 20 | 226 | 8 | | |
| 47 | 99 | 107 | 40 | 167 | 19 | 227 | 8 | | |
| 48 | 95 | 108 | 40 | 168 | 19 | 228 | 8 | | |
| 49 | 94 | 109 | 39 | 169 | 19 | 229 | 8 | | |
| 50 | 93 | 110 | 39 | 170 | 18 | 230 | 7 | | |
| 51 | 92 | 111 | 39 | 171 | 18 | 231 | 7 | | |
| 52 | 90 | 112 | 38 | 172 | 18 | 232 | 7 | | |
| 53 | 90 | 113 | 37 | 173 | 18 | 233 | 7 | | |
| 54 | 88 | 114 | 36 | 174 | 18 | 234 | 7 | | |
| 55 | 85 | 115 | 36 | 175 | 18 | 235 | 7 | | |
| 56 | 84 | 116 | 36 | 176 | 18 | 236 | 7 | | |
| 57 | 84 | 117 | 36 | 177 | 18 | 237 | 6 | | |
| 58 | 83 | 118 | 36 | 178 | 17 | 238 | 6 | | |
| 59 | 81 | 119 | 35 | 179 | 17 | 239 | 6 | | |
| 60 | 80 | 120 | 34 | 180 | 17 | 240 | 6 | | |

**Table A5.** The all citations of H.J. Kim in domain *D* of asymptotic normality of citations (Google Scholar, May 26 2023)

| | | | | | | | | | | | |
|---|---|---|---|---|---|---|---|---|---|---|---|
| 1 | 20949 | 61 | 1034 | 121 | 631 | 181 | 496 | 241 | 422 | 301 | 368 |
| 2 | 16571 | 62 | 1030 | 122 | 625 | 182 | 495 | 242 | 421 | 302 | 367 |
| 3 | 11853 | 63 | 1017 | 123 | 625 | 183 | 493 | 243 | 420 | 303 | 367 |
| 4 | 7813 | 64 | 1015 | 124 | 617 | 184 | 490 | 244 | 419 | 304 | 367 |
| 5 | 5326 | 65 | 1014 | 125 | 615 | 185 | 488 | 245 | 419 | 305 | 365 |
| 6 | 4773 | 66 | 1005 | 126 | 605 | 186 | 487 | 246 | 418 | 306 | 365 |
| 7 | 4683 | 67 | 988 | 127 | 604 | 187 | 484 | 247 | 415 | 307 | 363 |
| 8 | 4674 | 68 | 975 | 128 | 602 | 188 | 483 | 248 | 415 | 308 | 363 |
| 9 | 4605 | 69 | 973 | 129 | 598 | 189 | 482 | 249 | 412 | 309 | 362 |
| 10 | 2882 | 70 | 959 | 130 | 595 | 190 | 482 | 250 | 411 | 310 | 362 |
| 11 | 2863 | 71 | 944 | 131 | 595 | 191 | 481 | 251 | 411 | 311 | 362 |
| 12 | 2800 | 72 | 944 | 132 | 592 | 192 | 480 | 252 | 410 | 312 | 362 |
| 13 | 2727 | 73 | 931 | 133 | 591 | 193 | 479 | 253 | 408 | 313 | 359 |
| 14 | 2494 | 74 | 928 | 134 | 589 | 194 | 479 | 254 | 407 | 314 | 358 |
| 15 | 2344 | 75 | 922 | 135 | 587 | 195 | 479 | 255 | 407 | 315 | 358 |
| 16 | 2288 | 76 | 899 | 136 | 586 | 196 | 479 | 256 | 406 | 316 | 358 |
| 17 | 2263 | 77 | 896 | 137 | 585 | 197 | 476 | 257 | 405 | 317 | 356 |
| 18 | 2252 | 78 | 896 | 138 | 581 | 198 | 475 | 258 | 403 | 318 | 354 |
| 19 | 2218 | 79 | 894 | 139 | 580 | 199 | 473 | 259 | 400 | 319 | 354 |
| 20 | 2179 | 80 | 894 | 140 | 577 | 200 | 472 | 260 | 399 | 320 | 354 |
| 21 | 2005 | 81 | 875 | 141 | 575 | 201 | 471 | 261 | 397 | 321 | 351 |
| 22 | 1966 | 82 | 862 | 142 | 565 | 202 | 468 | 262 | 395 | 322 | 350 |
| 23 | 1946 | 83 | 857 | 143 | 565 | 203 | 466 | 263 | 395 | 323 | 349 |
| 24 | 1934 | 84 | 856 | 144 | 565 | 204 | 466 | 264 | 395 | 324 | 348 |
| 25 | 1884 | 85 | 855 | 145 | 563 | 205 | 465 | 265 | 392 | 325 | 348 |
| 26 | 1878 | 86 | 849 | 146 | 563 | 206 | 464 | 266 | 391 | 326 | 348 |
| 27 | 1795 | 87 | 827 | 147 | 561 | 207 | 464 | 267 | 391 | 327 | 347 |
| 28 | 1787 | 88 | 819 | 148 | 548 | 208 | 463 | 268 | 391 | 328 | 347 |
| 29 | 1750 | 89 | 811 | 149 | 545 | 209 | 462 | 269 | 391 | 329 | 346 |
| 30 | 1730 | 90 | 809 | 150 | 544 | 210 | 462 | 270 | 389 | 330 | 345 |
| 31 | 1669 | 91 | 807 | 151 | 540 | 211 | 458 | 271 | 389 | 331 | 344 |
| 32 | 1634 | 92 | 807 | 152 | 540 | 212 | 458 | 272 | 388 | 332 | 344 |
| 33 | 1591 | 93 | 806 | 153 | 540 | 213 | 457 | 273 | 387 | 333 | 344 |
| 34 | 1578 | 94 | 805 | 154 | 539 | 214 | 455 | 274 | 386 | 334 | 342 |
| 35 | 1567 | 95 | 793 | 155 | 539 | 215 | 455 | 275 | 385 | 335 | 340 |
| 36 | 1505 | 96 | 793 | 156 | 537 | 216 | 454 | 276 | 385 | 336 | 339 |
| 37 | 1502 | 97 | 776 | 157 | 534 | 217 | 453 | 277 | 384 | 337 | 338 |
| 38 | 1498 | 98 | 776 | 158 | 531 | 218 | 453 | 278 | 384 | 338 | **338** |
| 39 | 1437 | 99 | 757 | 159 | 529 | 219 | 451 | 279 | 384 | 339 | 337 |
| 40 | 1383 | 100 | 751 | 160 | 529 | 220 | 448 | 280 | 384 | 340 | 336 |
| 41 | 1352 | 101 | 744 | 161 | 528 | 221 | 448 | 281 | 384 | 341 | 336 |
| 42 | 1330 | 102 | 741 | 162 | 528 | 222 | 447 | 282 | 383 | 342 | 336 |
| 43 | 1267 | 103 | 738 | 163 | 525 | 223 | 446 | 283 | 383 | 343 | 335 |
| 44 | 1253 | 104 | 726 | 164 | 524 | 224 | 443 | 284 | 382 | 344 | 334 |
| 45 | 1211 | 105 | 722 | 165 | 523 | 225 | 439 | 285 | 380 | 345 | 334 |
| 46 | 1208 | 106 | 721 | 166 | 522 | 226 | 438 | 286 | 380 | 346 | 332 |
| 47 | 1194 | 107 | 717 | 167 | 521 | 227 | 436 | 287 | 380 | | 330 |
| 48 | 1185 | 108 | 715 | 168 | 520 | 228 | 432 | 288 | 378 | | |
| 49 | 1174 | 109 | 714 | 169 | 518 | 229 | 431 | 289 | 377 | | |
| 50 | 1167 | 110 | 710 | 170 | 517 | 230 | 430 | 290 | 377 | | |
| 51 | 1158 | 111 | 680 | 171 | 515 | 231 | 429 | 291 | 376 | | |
| 52 | 1150 | 112 | 678 | 172 | 514 | 232 | 429 | 292 | 375 | | |
| 53 | 1123 | 113 | 670 | 173 | 512 | 233 | 429 | 293 | 374 | | |
| 54 | 1086 | 114 | 662 | 174 | 511 | 234 | 426 | 294 | 373 | | |
| 55 | 1081 | 115 | 657 | 175 | 508 | 235 | 425 | 295 | 373 | | |
| 56 | 1080 | 116 | 657 | 176 | 506 | 236 | 425 | 296 | 372 | | |
| 57 | 1078 | 117 | 654 | 177 | 504 | 237 | 424 | 297 | 370 | | |
| 58 | 1074 | 118 | 645 | 178 | 504 | 238 | 424 | 298 | 369 | | |
| 59 | 1050 | 119 | 635 | 179 | 504 | 239 | 424 | 299 | 369 | | |
| 60 | 1042 | 120 | 631 | 180 | 498 | 240 | 424 | 300 | 369 | | |

Note:

$h$-index ($h = 70$, signed in bold) 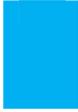 The bounds of domain $D$ of asymptotic normality of citations

**Table A6.** The all citations of R. Kessler in domain *D* of asymptotic normality of citations (Google Scholar, May 24 2023)

| | | | | | | | | | | | | | | | |
|---|---|---|---|---|---|---|---|---|---|---|---|---|---|---|---|
| 1 | 34554 | 61 | 1564 | 121 | 877 | 181 | 627 | 241 | 458 | 301 | 365 | 361 | 273 |
| 2 | 20715 | 62 | 1560 | 122 | 876 | 182 | 623 | 242 | 457 | 302 | 364 | 362 | 273 |
| 3 | 15855 | 63 | 1558 | 123 | 875 | 183 | 623 | 243 | 457 | 303 | 363 | 363 | 272 |
| 4 | 15375 | 64 | 1552 | 124 | 872 | 184 | 618 | 244 | 456 | 304 | 363 | 364 | 272 |
| 5 | 10173 | 65 | 1550 | 125 | 872 | 185 | 612 | 245 | 440 | 305 | 362 | | |
| 6 | 9435 | 66 | 1528 | 126 | 869 | 186 | 606 | 246 | 440 | 306 | 358 | | |
| 7 | 6639 | 67 | 1520 | 127 | 869 | 187 | 606 | 247 | 440 | 307 | 358 | | |
| 8 | 6320 | 68 | 1506 | 128 | 867 | 188 | 602 | 248 | 438 | 308 | 355 | | |
| 9 | 5141 | 69 | 1477 | 129 | 866 | 189 | 586 | 249 | 437 | 309 | 353 | | |
| 10 | 5125 | 70 | 1446 | 130 | 853 | 190 | 585 | 250 | 436 | 310 | 353 | | |
| 11 | 4787 | 71 | 1425 | 131 | 850 | 191 | 585 | 251 | 434 | 311 | 350 | | |
| 12 | 4491 | 72 | 1413 | 132 | 850 | 192 | 581 | 252 | 432 | 312 | 346 | | |
| 13 | 3669 | 73 | 1410 | 133 | 844 | 193 | 580 | 253 | 428 | 313 | 345 | | |
| 14 | 3431 | 74 | 1410 | 134 | 844 | 194 | 580 | 254 | 427 | 314 | 345 | | |
| 15 | 3419 | 75 | 1369 | 135 | 844 | 195 | 580 | 255 | 421 | 315 | 344 | | |
| 16 | 3370 | 76 | 1352 | 136 | 825 | 196 | 577 | 256 | 421 | 316 | 342 | | |
| 17 | 3068 | 77 | 1345 | 137 | 818 | 197 | 566 | 257 | 421 | 317 | 340 | | |
| 18 | 3044 | 78 | 1328 | 138 | 817 | 198 | 563 | 258 | 419 | 318 | 340 | | |
| 19 | 3038 | 79 | 1315 | 139 | 796 | 199 | 563 | 259 | 416 | 319 | 338 | | |
| 20 | 2995 | 80 | 1309 | 140 | 793 | 200 | 562 | 260 | 413 | 320 | 338 | | |
| 21 | 2960 | 81 | 1263 | 141 | 792 | 201 | 553 | 261 | 413 | 321 | 337 | | |
| 22 | 2951 | 82 | 1259 | 142 | 789 | 202 | 551 | 262 | 413 | 322 | 335 | | |
| 23 | 2938 | 83 | 1257 | 143 | 784 | 203 | 550 | 263 | 409 | 323 | 333 | | |
| 24 | 2919 | 84 | 1257 | 144 | 784 | 204 | 547 | 264 | 408 | 324 | 330 | | |
| 25 | 2910 | 85 | 1250 | 145 | 780 | 205 | 546 | 265 | 408 | 325 | 330 | | |
| 26 | 2772 | 86 | 1214 | 146 | 778 | 206 | 546 | 266 | 407 | 326 | 330 | | |
| 27 | 2771 | 87 | 1194 | 147 | 774 | 207 | 545 | 267 | 406 | 327 | **329** | | |
| 28 | 2763 | 88 | 1185 | 148 | 774 | 208 | 544 | 268 | 406 | 328 | 326 | | |
| 29 | 2752 | 89 | 1176 | 149 | 772 | 209 | 540 | 269 | 405 | 329 | 325 | | |
| 30 | 2727 | 90 | 1164 | 150 | 759 | 210 | 539 | 270 | 404 | 330 | 324 | | |
| 31 | 2630 | 91 | 1162 | 151 | 754 | 211 | 538 | 271 | 402 | 331 | 323 | | |
| 32 | 2536 | 92 | 1147 | 152 | 754 | 212 | 535 | 272 | 398 | 332 | 321 | | |
| 33 | 2470 | 93 | 1134 | 153 | 753 | 213 | 533 | 273 | 397 | 333 | 320 | | |
| 34 | 2424 | 94 | 1121 | 154 | 751 | 214 | 532 | 274 | 397 | 334 | 318 | | |
| 35 | 2395 | 95 | 1114 | 155 | 745 | 215 | 529 | 275 | 395 | 335 | 317 | | |
| 36 | 2304 | 96 | 1111 | 156 | 740 | 216 | 526 | 276 | 392 | 336 | 315 | | |
| 37 | 2261 | 97 | 1106 | 157 | 738 | 217 | 524 | 277 | 390 | 337 | 314 | | |
| 38 | 2213 | 98 | 1101 | 158 | 731 | 218 | 518 | 278 | 390 | 338 | 314 | | |
| 39 | 2205 | 99 | 1094 | 159 | 727 | 219 | 514 | 279 | 390 | 339 | 313 | | |
| 40 | 2082 | 100 | 1093 | 160 | 724 | 220 | 513 | 280 | 388 | 340 | 311 | | |
| 41 | 2045 | 101 | 1089 | 161 | 723 | 221 | 511 | 281 | 388 | 341 | 310 | | |
| 42 | 2016 | 102 | 1076 | 162 | 703 | 222 | 509 | 282 | 387 | 342 | 305 | | |
| 43 | 1987 | 103 | 1075 | 163 | 692 | 223 | 508 | 283 | 386 | 343 | 305 | | |
| 44 | 1945 | 104 | 1073 | 164 | 692 | 224 | 503 | 284 | 386 | 344 | 304 | | |
| 45 | 1931 | 105 | 1066 | 165 | 676 | 225 | 497 | 285 | 385 | 345 | 302 | | |
| 46 | 1868 | 106 | 1053 | 166 | 676 | 226 | 493 | 286 | 384 | 346 | 302 | | |
| 47 | 1853 | 107 | 1031 | 167 | 674 | 227 | 490 | 287 | 382 | 347 | 300 | | |
| 48 | 1787 | 108 | 1019 | 168 | 671 | 228 | 485 | 288 | 380 | 348 | 299 | | |
| 49 | 1750 | 109 | 1018 | 169 | 670 | 229 | 482 | 289 | 379 | 349 | 298 | | |
| 50 | 1749 | 110 | 1014 | 170 | 668 | 230 | 482 | 290 | 379 | 350 | 295 | | |
| 51 | 1747 | 111 | 992 | 171 | 663 | 231 | 481 | 291 | 379 | 351 | 294 | | |
| 52 | 1733 | 112 | 990 | 172 | 662 | 232 | 480 | 292 | 378 | 352 | 291 | | |
| 53 | 1711 | 113 | 985 | 173 | 656 | 233 | 474 | 293 | 378 | 353 | 290 | | |
| 54 | 1706 | 114 | 972 | 174 | 645 | 234 | 472 | 294 | 377 | 354 | 290 | | |
| 55 | 1696 | 115 | 967 | 175 | 642 | 235 | 467 | 295 | 377 | 355 | 289 | | |
| 56 | 1668 | 116 | 937 | 176 | 641 | 236 | 466 | 296 | 377 | 356 | 288 | | |
| 57 | 1637 | 117 | 927 | 177 | 637 | 237 | 462 | 297 | 376 | 357 | 286 | | |
| 58 | 1612 | 118 | 901 | 178 | 633 | 238 | 461 | 298 | 373 | 358 | 280 | | |
| 59 | 1566 | 119 | 893 | 179 | 630 | 239 | 460 | 299 | 370 | 359 | 280 | | |
| 60 | 1566 | 120 | 889 | 180 | 629 | 240 | 460 | 300 | 370 | 360 | 279 | | |

Note:

$h$-index ($h = 327$, signed in bold)

▮ The bounds of domain $D$ of asymptotic normality of citations

**We present below the calculations via software *Wolfram Mathematica* 11 of all bibliometric indicators and the indices from Table 5 of A. Einstein, as well as the calculations of the h-defect, the e_d and e_{d + 1} indices for H.J. Kim, R. Kessler, L. Leydesdorf, J. Hirsch, C.F. Gauss and C.G. Meyer.**

**1) Albert Einstein (Google Scholar, May 21 2023)**

```
e_59

(150565-
(23200+20218+7572+7216+7029+6124+4674+4309+4085+4007+3787+3505+2536+2323+1813
+1795+1687+1531+1516+1398+1389+1349+1279+1185+1153+1148+1119+1091+1026+995+98
4+974+863+860+832+746+668+567+558+520+518+507+500+487+483+471+440+436+
436+431+431+396+384+382+380+367+361+357+356)+123+119+117+116+115+110+108+107+
106+106+103+101+100+100+98+98+98+98+92+91+90+90+89+88+87-87^2)^0.5

88.27230596285564
```

**e_59=88.27230596285564**

**h_59=87**

```
e_60

(150565-
(23200+20218+7572+7216+7029+6124+4674+4309+4085+4007+3787+3505+2536+2323+1813
+1795+1687+1531+1516+1398+1389+1349+1279+1185+1153+1148+1119+1091+1026+995+98
4+974+863+860+832+746+668+567+558+520+518+507+500+487+483+471+440+436+
436+431+431+396+384+382+380+367+361+357+356+355)+123+119+117+116+115+110+108+
107+106+106+103+101+100+100+98+98+98+98+92+91+90+90+89+88+87-87^2)^0.5

86.23804264940155
```

**e_60=86.23804264940155**

```
As e_59> h_59 and e_60< h_60, d(the h- defect of approximative normality)=59
```
**59**

```
Ncit_59(h)

150565-
(23200+20218+7572+7216+7029+6124+4674+4309+4085+4007+3787+3505+2536+2323+1813
+1795+1687+1531+1516+1398+1389+1349+1279+1185+1153+1148+1119+1091+1026+995+98
4+974+863+860+832+746+668+567+558+520+518+507+500+487+483+471+440+436+436+431
+431+396+384+382+380+367+361+357+356)+123+119+117+116+115+110+108+107+106+106
+103+101+100+100+98+98+98+98+92+91+90+90+89+88+87
```

**15361**

```
Ncit_59
```

```
161009-
(23200+20218+7572+7216+7029+6124+4674+4309+4085+4007+3787+3505+2536+2323+1813
+1795+1687+1531+1516+1398+1389+1349+1279+1185+1153+1148+1119+1091+1026+995+98
4+974+863+860+832+746+668+567+558+520+518+507+500+487+483+471+440+436+436+431
+431+396+384+382+380+367+361+357+356)
```

**23255**

```
q_59=2*N_cit^(59)(h)/h_59^2-1
2*15361/87^2-1.0
```

**3.05892**

```
q_ 59/e_59
3.058924560708152/88.27230596285564
```

**0.0346533**

```
I_ 59=(a_59, b_59)=(((1- q_ 59/e_59)*h_59/0.5404446394667307)^2, ((1+ q_ 59/e_59)*h_59/0.5404446394667307)^2

a_ 59
((1-0.03465327576233622`)*87/0.5404446394667307)^2
```

**24149.2**

```
b_ 59
((1+0.03465327576233622`)*87/0.5404446394667307)^2
```

**27741.2**
```
So I_ 59=(24149.20,27741.23)

J_ 59=(c_ 59,d_ 59)=(24149.20+137754, 27741.23+137754)=
```
**(161903.20,165495.23)**

**Remark. Observe that the lower bound of the interval J_ 59 is 161903.20 which is very close to the total number of Einstein's citations, Namely c_ 59-N_cit=161903.20-161009=894.2, and the relative error \delta(c_ 59)= 894.2/161009=0.005554**

```
A'=\bar{J_ 59}=(ç_ 59+d_ 59)/2=(161903.20+165495.23)/2 =
```
**163699.22**

```
A'-N_cit=163699.22-161009 = 
```
**2690.22**

```
\delta{A'}=(A'-N_cit)/N_cit=(163699.22-161009)/161009=
```
**0.0167085**

```
Since e_59- h_59 = 88.27230596285564-87=1.27230596285564>1, we define \beta= 1 and \alpha=0, and so

B'=\alpha*c_59+\beta*d_59=0*161903.20+1*165495.23=
```
**165495.23**

```
\delta{B'}=(B'-N_cit)/N_cit=(165495.23-161009)/161009=
```
**0.02786322503710979**

```
Now we calculate the intervals I_60 and J_60.

e_60
```

(150565-(23200+20218+7572+7216+7029+6124+4674+4309+4085+4007+3787+3505+2536+2323+1813+1795+1687+1531+1516+1398+1389+1349+1279+1185+1153+1148+1119+1091+1026+995+984+974+863+860+832+746+668+567+558+520+518+507+500+487+483+471+440+436+436+431+431+396+384+382+380+367+361+357+356+355)+123+119+117+116+115+110+108+107+106+106+103+101+100+100+98+98+98+98+92+91+90+90+89+88+87−87^2)^0.5

**86.23804264940155**

Ncit_60(h)

150565-(23200+20218+7572+7216+7029+6124+4674+4309+4085+4007+3787+3505+2536+2323+1813+1795+1687+1531+1516+1398+1389+1349+1279+1185+1153+1148+1119+1091+1026+995+984+974+863+860+832+746+668+567+558+520+518+507+500+487+483+471+440+436+436+431+431+396+384+382+380+367+361+357+356+355)+123+119+117+116+115+110+108+107+106+106+103+101+100+100+98+98+98+98+92+91+90+90+89+88+87

**15006**

Ncit_60

161009-(23200+20218+7572+7216+7029+6124+4674+4309+4085+4007+3787+3505+2536+2323+1813+1795+1687+1531+1516+1398+1389+1349+1279+1185+1153+1148+1119+1091+1026+995+984+974+863+860+832+746+668+567+558+520+518+507+500+487+483+471+440+436+436+431+431+396+384+382+380+367+361+357+356+355)

**22900**

$q\_60 = 2 \cdot N\_cit^{(60)}(h)/h\_60^2 - 1 = 2 \cdot 15006/87.0^2 - 1$

**2.9651208878319464**

$q\_60 / e\_60$

2.9651208878319464/86.23804264940155

**0.03438297990930251**

$I\_60 = (a\_60, b\_60) = (((1 - q\_60/e\_60) \cdot h\_60 / 0.5404446394667307)^2, ((1 + q\_60/e\_60) \cdot h\_60 / 0.5404446394667307)^2$

a_60
((1-0.03438297990930251)*87/0.5404446394667307)^2

**24162.72**

b_60
((1+0.03438297990930251)*87/0.5404446394667307)^2

**27726.74**

So $I\_60 = (24162.72, 27726.74)$

$J\_60 = (c\_60, d\_60) = (24162.72 + 138109, 27726.74 + 138109) =$ **(162271.72, 165835.74)**

$A'' = \bar{J\_60} = (ç\_60 + d\_60)/2 = (162271.72 + 165835.74)/2 =$ **164054.73**

A''-N_cit=164054.73-161009 = 3045.73

\delta{A''}=(A''-N_cit)/N_cit=(164054.73-161009)/161009= 0.0189165

A=(A'+ A'')/2= (163699.22+164054.73)=**163876.98**
  A-N_cit=161009-163876.98= **-2867.98**

\delta{A}=(A-N_cit)/N_cit=(163876.98-161009)/161009=**0.0178125**

Since h_60- e_60 = 87-86.23804264940155=0.2380426494015495<1, we define \beta'=(the fractional part of h_60- e_60)=0.2380426494015495 and \alpha'=0.7619573505984505, and so

B''=\beta'*c_60+\alpha'*d_60=0.2380426494015495*161903.20+0.7619573505984505*165495.23=**164640.17**

B''-N_cit=164640.17-161009=**3631.17**

\delta{B''}=(B''-N_cit)/N_cit=(**164640.17**-161009)/161009=**0.0225526**

**Finally,**

**B=** (B'+ B'')/2=(165495.23+164640.17)/2=**165067.70**

\delta{B}=(B-N_cit)/N_cit=(**165067.70**-161009)/161009= **0.0252079**

**Notice that \delta{B}=0.0252079>0.0178125=\delta{A}**

**2) H.J. KIM (Google Scholar, May 26 2023) h=338, N_cit=518589**

e_13
($\overline{312483-}$
(20949+16571+11853+7813+5326+4773+4683+4674+4605+2882+2863+2800+2727)+(338+337+336+336+335+334+334)-333^2)^0.5
**e_13= 333.80383460949037>333= h_13**

e_14
($\overline{312483-}$
(20949+16571+11853+7813+5326+4773+4683+4674+4605+2882+2863+2800+2727+2494)+(338+337+336+336+335+334+334)-332^2)^0.5
**e_14= 331.052865868882<332= h_14**

**So d (the $h$- defect)=333**

**The first 200 publications**

(20949+16571+11853+7813+5326+4773+4683+4674+4605+2882+2863+2800+2727+2494+2344+2288+2263+2252+2218+2179+2005+1966+1946+1934+1884+1878+1795+1787+1750+1730+1669+1634+1591+1578+1567+1505+1502+1498+1437+1383+1352+1330+1267+1253+1211+1208+1194+1185+1174+1167+1158+1150+1123+1086+1081+1080+1078+1074+1050+1042+1034+1030+1017+1015+1014+1005+988+975+973+959+944+944+931+928+922+899+896+896+894+894+875+862+857+856+855+849+827+819+811+809+807+807+806+805+793+793+776+776+757+751)+(744+741+738+726+722+721+717+715+714+710+680+678+670+662+657+657+654+645+635+631+631+625+625+617+615+605+604+602+598+595+595+592+591+589+587+586+585+581+580+577+575+565+565+565+563+563+561+548+545+544+540+540+540+539+539+5

37+534+531+529+529+528+528+525+524+523+522+521+520+518+517+515+514+512+511+508+506+504+504+504+498+496+495+493+490+488+487+484+483+482+482+481+480+479+479+479+479+476+475+473+472)

**3) Ronald Kessler (Google Scholar, May 24 2023)**

e_94

(432109-(34554+20715+15855+15375+10173+9435+6639+6320+5141+5125+4787+4491+3969+3431+3419+3370+3068+3044+3038+2995+2960+2951+2919+2916+2772+2771+2763+2752+2727+2630+2563+2470+2424+2395+2304+2261+2213+2205+2082+2045+2016+1987+1945+1931+1868+1853+1787+1750+1749+1749+1747+1733+1711+1706+1696+1668+1637+1612+1566+1566+1564+1560+1558+1552+1550+1528+1520+1506+1477+1446+1425+1413+1410+1369+1352+1345+1328+1315+1309+1263+1259+1257+1257+1250+1214+1194+1185+1176+1164+1162+1147+1147+1134+1121)+326+325+324+323+321+320+318+317+315+314+314+313+311+310+305+305+304+302+302+300+299+298+295+294+291+290+290+289+288+286+280+280+279+273+273+272+272**-270**^2)^0.5

**e_94=271.157>270=h_94**

e_95

(432109-(34554+20715+15855+15375+10173+9435+6639+6320+5141+5125+4787+4491+3969+3431+3419+3370+3068+3044+3038+2995+2960+2951+2919+2916+2772+2771+2763+2752+2727+2630+2563+2470+2424+2395+2304+2261+2213+2205+2082+2045+2016+1987+1945+1931+1868+1853+1787+1750+1749+1749+1747+1733+1711+1706+1696+1668+1637+1612+1566+1566+1564+1560+1558+1552+1550+1528+1520+1506+1477+1446+1425+1413+1410+1369+1352+1345+1328+1315+1309+1263+1259+1257+1257+1250+1214+1194+1185+1176+1164+1162+1147+1147+1134+1121+1114)+326+325+324+323+321+320+318+317+315+314+314+313+311+310+305+305+304+302+302+300+299+298+295+294+291+290+290+289+288+286+280+280+279+273+273+272+272+270**-270**^2)^0.5

**e_ 95=269.596<270=h_95**

**So d(the $h$-defect)=270**

**4) L. Leydesdorf (Google Scholar, May 24 2023)**

e_14

(54898-(12829+3505+3269+1738+1395+1346+1170+1089+841+765+725+704+699+697)+117+117+116+115+113+113**-111**^2)^0.5

**e_14=111.786**

e_15
(54898-(12829+3505+3269+1738+1395+1346+1170+1089+841+765+725+704+699+697+655)+117+117+116+115+113+113**-110**^2)^0.5

**e_15=109.827**

**So d(the $h$-defect)=111**

**5) J. Hirsch (Google Scholar May 24 2023)**

e_12
(20507-
(7159+2378+881+710+686+659+504+440+369+364+361+275)+(60+59+59+59+59+57+56)-
55^2)^0.5
**e_12=55.7225 >55=h_ 12**

e_13
(20507-
(7159+2378+881+710+686+659+504+440+369+364+361+275+263)+(60+59+59+59+59+57+56
+55)-55^2)^0.5
**e_13=53.8238<55= h_ 13**

**So d(the $h$- defect)=55**

**6) C.F. Gauss (GS, May 28 2023,  h=44 Nc=11606**
N_cit(h)
1953+1417+534+533+510+453+368+323+313+220+205+176+147+146+144+141+131+127+116
+116+115+114+108+106+97+96+96+90+90+83+81+75+74+68+66+62+61+60+58+57+53+47+45
+45
9920

**N_cit(h)=9920**

e
(9920-44^2)^0.5
89.3532
**e=89.3532**

q
2*9920/44^2-1.0
9.24793
**q=9.24793**

q/e
9.258264462809917/89.3532316147547`
0.103614
**q/e=0.103614**

h+q
9.24793388429752+44
53.2479
**h+q=53.2479**

h_AN
0.5404446394667307*(11606)^0.5
**h_AN=58.2227**
lower bound a
the interval I=(a,b)=(((1-
q/e)*(h+q)/0.5404446394667307))^2,((1+q/e)*(h+q)/0.5404446394667307))^2

((1-0.10361420953107553`)*53.247933884297524`/0.5404446394667307)^2
7799.97
**a=7799.97**

uper bound b
((1+0.10361420953107553`)*53.247933884297524`/0.5404446394667307)^2
11823.3

**b=11823.3**

**I=(a,b)=(7799.97, 11823.3)**
Since e-h=45.41>1(even >h=44), we assume the uper bound b=11826.47 for the estimate of N_cit=11606, i.e., \beta=1, \alpha=0 and **B=\alpha *a+\beta *b=11826.47.**
Hence, the absolute error of N_cit is
**B-N_cit=11826.47-11606=220.46999999999935**
and the relative error of N_cit is equal
**\delta(B)=(B-N_cit)/N_cit=0.01899620885748745**
uper bound b
((1+0.10361420953107553`)*53.247933884297524`/0.5404446394667307)^2
11823.3
**I=(a,b)=(7799.97,11823.27)**

Since e-h=45.35>1(even >h=44), we assume \beta=1; \alpha=0 and the uper bound b=11823.27 is the estimate of N_cit=11606, i.e.,
**B=(\alpha*a+\beta*b=b)=11823.27.** Hence,
**B-N_cit=11823.27-11606=217.27000000000044**
and the relative error of N_cit is equal
**\delta(A)=(A-N_cit)/N_cit=0.003912560036966091**

e_ 18
(9920-(1953+1417+534+533+510+453+368+323+313+220+205+176+147+146+144+141+131+127)+42+41+39+39+38+36+36+36-34^2)^0.5
35.0714
**e_ 18=35.0714>35= h_ 18**

e_ 19
(9930-(1953+1417+534+533+510+453+368+323+313+220+205+176+147+146+144+141+131+127+126)+42+41+39+39+38+36+36+36+34-34^2)^0.5
33.8821

**e_ 19=33.8821<34= h_ 19**

**So d (the $h$-defect)=18**

**7) Christian G Meyer (Google Scholar May 28 2023), N_cit=49100, h=95**
e_ 9
(36684-(11973+3573+1509+697+590+524+512+476+468)+93+93+92+92+91-91^2)^0.5
92.4229
**e_ 9=92.4229>91=h_9**

e_ 10
(36684-(11973+3573+1509+697+590+524+512+476+468+412)+93+93+92+92+91-91^2)^0.5
90.1665
**e_ 10=90.1665<91=h_91**

**So d(the $h$-defect)=91**